\documentclass{article}
\usepackage{graphicx}
\usepackage{amssymb,amsfonts,amsthm}
\usepackage{amsmath,amsthm,amssymb}
\usepackage{amsfonts}

\newtheorem{intro}{Theorem}[section]

\newtheorem{cintro}[intro]{Corollary}
\newtheorem{theorem}{Theorem}[section]
\newtheorem{lemma}[theorem]{Lemma}

\newtheorem{cy}[theorem]{Corollary}
\newtheorem{prop}[theorem]{Proposition}
\theoremstyle{definition}
\newtheorem{df}[theorem]{Definition}

\newtheorem{rk}[theorem]{Remark}

\newcounter{ppp}

\newcommand{\me}{\medskip}

\newcommand{\Lab}{\phi}

\newcommand{\tool}{\stackrel{\ell}{\too} }

\newcommand{\ttt}{{\cal T}}

\newcommand{\bb}{{\cal B}}
\newcommand{\topp}{{\bf top}}
\newcommand{\ttopp}{{\bf ttop}}
\newcommand{\tbott}{{\bf tbot}}
\newcommand{\bott}{{\bf bot}}
\newcommand{\vk}{van Kampen }

\newcommand{\iv}{^{-1}}

\newcommand{\too}{\to }

\begin{document}

\renewcommand{\theequation}{\thesection.\arabic{equation}}

\title{Space functions of groups.}
\author{A.Yu.Olshanskii \thanks{The
author was supported in part by the NSF grant DMS 0700811 and by the Russian Fund for Basic Research grant 08-01-00573}}
\maketitle

\begin{abstract} We consider space functions $s(n)$ of finitely presented groups 
$G =\langle A\mid R\rangle .$ (These functions have a natural geometric analog.) To define $s(n)$ we start with a word $w$ over $A$ of length at most $n$ equal to $1$ in $G$ and use relations from $R$ for elementary transformations to obtain the empty word; $s(n)$ bounds from above the tape space (or computer memory) one needs to transform any word of length at most $n$ vanishing in $G$ to the empty word. 
One of the main results obtained is the following criterion: 
A finitely generated group $H$ has decidable word problem of polynomial space complexity if and only if $H$ is a subgroup of a finitely presented group $G$ with a polynomial space function. 
\end{abstract}


\medskip


\section{Introduction}\label{intro}

 Time and space
complexities are the main properties of algorithms. Their counterparts 
in Group
Theory are the Dehn and filling length (or space) functions of finitely presented groups.
In this paper, we study the interrelation of space functions of groups and the space complexity
of the algorithmic word problem in groups.

Let $G=\langle A\mid R\rangle $ be a group presentation, where $A$ is a set of generators
and $R$ is a set of defining relators. Recall that relators belong to the free group with basis $A$,
and a group word $w$ in generators $A$ (i.e., a word over $A^{\pm 1}$) represents the identity
of $G$ iff there is a {\it derivation}

\begin{equation}\label{rewri}
w\equiv w_0\to w_1\to \dots\to w_{t-1}\to w_t\equiv 1
\end{equation}
where $1$ is the empty word, the sign $\equiv$ is used for letter-by-letter equality of words, and for every $i=1,\dots,t$, the word $w_i$ results from $w_{i-1}$
after application of one of the elementary $R$-transformations. As such transformations one
can take free reductions of subwords $aa^{-1}\to 1$ ($a\in A^{\pm 1}$), removing subwords
$r^{\pm 1}$, where $r\in R$, and the inverse transformations.

The minimal non-decreasing function $f(n)\colon \mathbb{N}\to \mathbb{N}$ such that
for every word $w$ vanishing in $G$ and having length $||w||\le n,$ there exists
a derivation (\ref{rewri}) with $t\le f(n),$ is called the Dehn function of the presentation 
$G=\langle A\mid R\rangle$ \cite{Gr}. For {\it finitely presented} groups (i.e., both sets $A$ and $R$ are finite)
Dehn functions are usually taken up to equivalence to get rid of the dependence  on
a finite presentation for $G$ (see \cite{MO}). To introduce this equivalence $\sim,$ we write $f\preceq g$
if there is a positive integer $c$ such that
\begin{equation}\label{prec}
f(n)\le cg(cn)+cn\;\;\; for \;\; any \;\;n\in \mathbb{N} 
\end{equation}

For example, we say that a function $f$ is {\it polynomial} if $f\preceq g$ for a polynomial $g.$
From now on we use the following equivalence for non-decreasing functions $f$ and $g$ on $\mathbb{N}.$

\begin{equation}\label{equiv}
f\sim g\;\;\; if\;\; both\;\; f\preceq g\;\; and\;\; g\preceq f   
\end{equation}

It is not difficult to see that the Dehn function $f(n)$ of a finitely presented group $G$ is recursive (or bounded from above by a recursive function) iff the word problem is algorithmically decidable for $G$ (see \cite{Ger1}, \cite{BRS}). In this case, the word
problem can be solved by a primitive algorithm that, given a word $w$ of length $n,$ just checks if
there exists a  derivation (\ref{rewri}) of length $\le f(n).$  Therefore the nondeterministic {\it time}  complexity of the word problem in $G$ is bounded from above by $f(n).$ Moreover if $H$ is a finitely generated subgroup of $G,$
then one can use the rewriting procedure (\ref{rewri}) for $H,$ and so the nondeterministic time complexity of the word
problem for $H$ is also bounded by $f(n).$  

It turns out that a converse statement is also true. Assume that the word problem can be solved in
a {\it finitely generated} group $H$ by a nondeterministic Turing machine ($NTM$) with  time complexity $T(n).$
Then $H$ is a subgroup of a {\it finitely presented} group $G$ with Dehn function equivalent to $n^2T(n^2)^4$
\cite{BORS}.
As the main corollary, one concludes that the word problem of a finitely generated group $H$ has time
complexity of class $NP$ (i.e., there {\it exists} a non-deterministic  algorithm of polynomial time
complexity, that solves the word problem for $H$) iff $H$ is a subgroup of a finitely presented group
with polynomial Dehn function.

We want to obtain similar statements for space functions. It is clear that to perform
the rewriting (\ref{rewri}) one needs space equal to $\max_{0\le i\le t}||w_i||,$ and this observation 
leads to the definition of space function of a finitely presented group. However when handling
groups, one can enlarge the set of elementary $R$-transformations and obtain different definitions
of space functions. One can either (1) consider only the transformations we defined above and obtain
the filling length function introduced by Gromov \cite{Gr} (also see  \cite{GR}, \cite{Bir1}); or (2) also 
allow replacement of words by their cyclic
permutations as was suggested by Bridson and Riley \cite{BR}; or (3) additionaly use the replacement of a word
$w\equiv uv$ by the pair $(u,v)$ if both $u$ and $v$ are trivial in $G$ (see \cite{BR} again).
Starting with these different sets of transformations one comes to different space functions
called in \cite{BR}, respectively, filling length function ($FL$), free filling length functions
($FFL$), and fragmenting free filling length functions ($FFFL$). Each of these functions has a
visual geometric interpretation in terms of the transformations of loops in the Cayley complex of $G$ 
using, respectively,  null-homotopy, free null-homotopy, and  free null-homotopy with bifurcations.
It is proved in \cite{BR} that these functions behave differently for the same finitely presented
group  $G$, for instance, $FFFL$ can grow linearly while $FL$ and $FFL$ have exponential growth.
There are many other features of these functions presented in \cite{BR} to justify their
``inclusion in the pantheon of filling invariants''.

 In this paper, we choose the third version ($FFFL$), and this choice is justified by the theorems on the connections between such functions 
and the space complexity of the word problem for groups.\footnote{An embedding statement for
the $FL$-functions was conjectured by J.-C. Birget in \cite{Bir1}.}  Thus we operate with finite sequences of words
$W=(w_1,\dots,w_s)$ over a group alphabet $A.$ Given a finitely presented group $G,$ we say that
a finite sequence $W'=(w'_1,\dots,w'_{s'})$ results from $W$ after application of an elementary $R$-transformation if $s'\in\{s-1,s,s+1\}$ and one of the following is done for some $w_i$ ($i=1,\dots, s$):
\begin{itemize}

\item a subword $aa^{-1}$ is removed from or inserted to  $w_i$ ($a\in A^{\pm 1}$);
\item a subword $r$ or $r^{-1}$ is removed from or inserted to  $w_i$ ($r\in R$);
\item $w_i$ is replaced by a cyclic conjugate;
\item $w_i\equiv uv$, and $w_i$ is replaced by the pair $u,v$, i.e. $W'=(w_1,\dots, w_{i-1},u,v,w_{i+1},\dots, w_s)$;
\item $w_i$ is removed if it is empty, i.e. $W'= (w_1,\dots, w_{i-1},w_{i+1},\dots, w_s).$

\end{itemize}

Clearly, we have $w=1$ in the group $G$ iff there exists an $R$-rewriting starting with $(w)$
and ending with the empty string $(\;)$.

For every finite sequence $W=(w_1,\dots, w_s)$ we set $||W||=\sum_{i=1}^s ||w_i||.$ By definition,
the space of a rewriting $W_0\to\dots\to W_t$ is $\max_{j=0}^t ||W_j||.$ If a word vanishes in $G,$
then $space(w)=space_G(w)$ is the minimum of spaces of all rewritings starting with $(w)$ and ending
with the empty string. The {\it space function} of the group presentation $G=\langle A\mid R\rangle$
(or briefly, of the group $G$) is the function  
$$S_G(n) = \max(space(w),\;\; where\;\; w=1 \;\;in \;\; G\;\; and \;\; ||w||\le n)$$
The space functions of finitely presented groups will be regarded up to the equivalence defined
by (\ref{prec}) and (\ref{equiv}), and so their growth will be at least linear. It is observed
in \cite{BR} that the equivalence class of $S_G$ does not depend on a finite presentation of
the group $G$, moreover this class is invariant under quasi-isometries. 

An accurate definition of the space complexity (function) $f(n)$ for a Turing machine ($TM$) will be recalled in Subsection \ref{definitions}.
Now we just note that it is conventional that for a multi-tape $TM,$ the function $f(n)$ counts only the space of work tapes. However since the space complexities of machines are taken here up to the same equivalence as  the space functions of groups, the adding of the space of the input
tape does not change the equivalence class. 

The sequence $W_0\to\dots\to W_t$ can be easily produced by an $NTM$ such that the
computation needs at most $2\max_{i=1}^t ||W_i||+const$ tape squares. (See also Section 3 of \cite{SBR}
or Remark 2.4 in \cite{BR}.) This immediately implies

\begin{prop} \label{propos} The space function of a finitely presented group $G$ is equivalent to the space complexity
of a non-deterministic two-tape $TM$. The language accepted by this machine coincides with 
the set of words equal to $1$ in the group.
\end{prop} 

In particular, the non-deterministic space complexity of the word problem in a finitely
presented group $G$ does not exceed the space function of $G.$ It follows from \cite{MO}, \cite{CMO} that
the literal converse statement fails. Moreover, a counter-example can be given by Baumslag's
$1$-relator group \cite{B} $G=\langle a,b\mid (aba^{-1})b(aba^{-1})^{-1}=b^2\rangle$
because the space function of $G$ is not bounded from above by any multi-exponential function
(see papers of S.Gersten \cite{Ger} and A.Platonov \cite{Pl}) while the space complexity of the
word problem for $G$ is at most exponential. This was proved by M. Kapovich and Schupp (unpublished), moreover, it is polynomial (announced by A. G. Myasnikov, A. Ushakov, and Dong Wook Won).
Therefore the correct formulation has to take into consideration that the algorithm from Proposition \ref{propos}
solves the word problem not only for $G$ but also for every finitely generated subgroup of the group $G.$ In the deterministic case we get a sharper formulation:

\begin{theorem} \label{main1} Let $H$ be a finitely generated group such that the word problem for $H$ is
decidable by a deterministic $TM$ ($DTM$) with space complexity $f(n)$. Then $H$ is a subgroup of a finitely
presented group $G$ with space function equivalent to $f(n).$ 

One can choose the group $G$ so  that $f(n)$ is also equivalent to the logarithm of the Dehn function of $G$.

\end{theorem}

The main corollary of this theorem applies to polynomial space complexity. We say that a finitely 
generated group $G$ belongs to the class $PSPACE$ (to $NPSPACE$) if 
the word problem for $G$ is decidable by some $DTM$ (some $NTM$) with polynomial space
complexity. But $NPSPACE = PSPACE$ by the remarkable theorem of Savitch (see \cite{DK}, Corollary 1.31)
(Since there is no similar equality for deterministic  {\it time} complexity, it therefore
has no natural algebraic counterpart.)
Therefore in contrast to the main result of \cite{BORS},  Proposition \ref{propos} and 
Theorem \ref{main1} eliminate any non-determinism in

\begin{cy} \label{pspace} A {\em finitely generated} group  $H$ belongs to $PSPACE$ iff
$H$ is a subgroup of a {\em finitely presented} group $G$ having polynomial
space function.
\end{cy}

Thus, given a 'good' algorithm solving 
the word problem in $H$ (e.g., using a matrix representation of $H$, etc.), it is
possible to find a bigger group $G$ whose deterministically modified ('silly')
natural algorithm solves the word problem for both $H$ and $G,$ and whose
space complexity is not much worse than the space complexity of the original algorithm.
Note that most common finitely generated groups belong to $PSPACE,$ with a wide range 
from linear
groups (two-tape  machines need just logspace on the work tape to solve their word problem \cite{LZ}) to free Burnside groups
of large odd exponents (N.Boatman, unpublished).

Another natural question raised in this paper is the realization problem: Which functions
$f(n)\colon \mathbb {N}\to\mathbb {N}$ are, up to equivalence, the space functions of finitely presented
groups? There are not many examples; linear and exponential ones can be found in \cite{BR}, but it is not easy even to specify a group with
space function $n^2.$

\begin{theorem} \label{realiz} The space complexity $f(n)$ of arbitrary $DTM$ $M$ is equivalent
to the space function of some finitely presented group $G$.

   In addition, one can choose the group $G$ so that $f(n)$ is also equivalent to the logarithm of the Dehn function of $G,$ and there is a one-to-one linear time reduction of
the decision problem of $M$ to the word problem of $G.$

\end{theorem}

This theorem gives a tremendous class of space functions of groups, including functions equivalent to 
$[\exp{\sqrt n}],$ $[n^k]$ ($k\in \mathbb{N}$), $[n^k\log^l n],$  $[n^k\log^l (\log\log n)^m],$ etc.  Note that we do not assume in the formulation of Theorem \ref{realiz} that the function
$f(n)$ is superadditive (i.e., $f(m+n)\ge f(m)+f(n)$) or grows sufficiently fast. (Compare with Theorem 1.2 \cite{SBR}
on Dehn functions of groups.) It follows, in particular, that there exists a finitely presented
group whose space function is not equivalent to any superadditive function. Recall that 
it is unknown if the Dehn function of arbitrary finitely presented group
 is equivalent to a superadditive function; see  \cite{GS}. 
 
 J.-C. Birget observed (email communication) that  in addition, Theorem \ref{realiz}  implies 
 a statement that does not involve space functions of groups (see also Remark \ref{extend}.)
 The first claim of Corollary \ref{complete} below follows, in fact, from each of
 the papers  \cite{T} (by B.A.Trakhtenbrot) and \cite{V} (by M.K. Valiev).  
 Let us start with a 
$DTM$ $M$ solving a $PSPACE$ complete problem (for the definition and 
the existence see \cite{DK}, Theorem 3.23). Then choose a finitely presented 
group $G=G(M)$ in accordance with Theorem \ref{realiz}, i.e., $G$ has a 
polynomial space function. By Proposition \ref{propos}, 
$G$ belongs to $NPSPACE = PSPACE,$ and so we have 

\begin{cy} \label{complete} There is a finitely presented group $G$ with  $PSPACE$ complete word problem. Moreover the group $G$ can be chosen with polynomial space function.
\end{cy}

Theorem \ref{main1} implies a non-deterministic corollary. To formulate it we recall that a function $f: \mathbb{N}\to \mathbb{N}$ is called
{\it fully space-constructible} ($FSC$) if there exists a two-tape $DTM$ that halts on any input $x$
of length $n$ after visiting exactly $f(n)$ tape squares of the work tape. Most common functions
are $FSC$ (see \cite{DK}).

\begin{cy} \label{main} Let $H$ be a finitely generated group such that the word problem for $H$ is
decidable by an $NTM$ having $FSC$ space complexity $f(n)$. Then $H$ is a subgroup of a finitely
presented group $G$ with space function equivalent to $f(n)^2.$
\end{cy}

Finally, we describe the functions $n^{\alpha}$ which are (up to equivalence)  space functions of groups.
Our approach is to modify the proof of Savitch's theorem from \cite{DK} and the proof from \cite{SBR}, where the similar problem was
considered for Dehn functions if $\alpha \ge 4,$ and close necessary and sufficient
conditions were obtained. (See also a dense series of examples with $\alpha\ge 2$ presented by Brady and Bridson \cite{BB}.)  Now we have $\alpha\ge 1$ in Corollary \ref{alpha} below. 
Also it is remarkable that for space functions the necessary and sufficient
conditions just coincide. 
To formulate the criterion, 
we  call a real number $\alpha$ {\it computable with space} $\le f(m)$ if there exists a
$DTM$ which, given a natural number $m,$ computes a binary rational approximation of $\alpha$
with an error  $O(2^{-m}),$ and the space of this computation $\le f(m).$

We have obtained the following
criterion.

\begin{cy} \label{alpha} For a real number $\alpha\ge 1,$ 
the function $[n^{\alpha}]$  is equivalent to the space function of a finitely
presented group iff $\alpha$ is  computable with space $\le 2^{2^m}.$ 

\end{cy}

It follows that functions $[n^\pi],$ $[n^{\sqrt e}]$, and $[n^{\alpha}]$ with any algebraic $\alpha\ge 1$ are  all  space functions of finitely presented groups. 

The space function is defined for a simply connected geodesic metric space under
some weak restrictions, in particular, for the universal cover of any closed connected
Riemannian manifold. (See \cite{BR} for details; we just note here that to define
the space (= $FFFL$) function, one should consider free homotopy with the possibility
of separating a loop into two loops in a bifurcation point.) It is proved in
\cite{BR} (Theorem E) that if a finitely presented group $G$ acts properly and
cocompactly by isometries on such a space $X$, then the space function of $G$ is
equivalent to the space function of $X.$ Since every finitely presented group is
a fundamental group of a connected compact Riemannian manifold, we can use corollaries
\ref{realiz} and \ref{alpha} to formulate one more

\begin{cy} \label{Riem} For every space complexity $f(n)$ of a $DTM,$ 
there exists a closed connected Riemannian manifold $M$ such that the
space function of the universal cover $\tilde M$ is equivalent to $f(n)$. 

In particular,
if  a real number $\alpha\ge 1$ is computable with space $\le 2^{2^m},$
then there exists such a universal cover $\tilde M$ with  space function equivalent to $n^{\alpha}.$
\end{cy}   

\medskip

To some extent, our constructions can be traced back to the works of P.Novikov, Boon, Britton and
other authors who invented group-theoretical interpretations of $TM$-s (see \cite{R}, ch. 12).
The hub relator is a multiple copy of the accept configuration of a machine. Then we
use  the language of van Kampen diagrams and construct a disc diagram  
for arbitrary accepting computation. A computational disc has the hub cell in the center,
and this hub is surrounded by a number of similar sectors.
At first we have to estimate the sizes
of computational discs, and to do this we should know the {\it generalized} space complexity of a machine.
It estimates the space of computations starting with {\it arbitrary} accept configuration, not only
with input ones. So we should modify the initial machine to be able to control the generalized
space complexity. (The modification from \cite{SBR} helps to control the time complexity but corrupts
the space complexity.) 

The next modification is due to the symmetry of algebraic relations: since $u=v$ always implies $v=u,$
the algebraic version of a machine $M$ always interprets the symmetrization of $M.$ Thus
we are concerned that the symmetrization preserves the basic characteristics, e.g. the accepted
language and the space functions. We are able to do this only if the initial machine is deterministic
or can be transformed to a deterministic machine under the control of basic properties. (The known symmetrization trick  from \cite{Bir} or \cite{SBR} does not work here since it does not preserve
the space complexity.) This causes the
restrictions in the formulations of Theorem \ref{main1}  and Corollary \ref{main}. 

The interpretation problem for groups remains much harder than for semigroups even after
modifying the machine because the group theoretic simulation can execute unforeseen computations with non-positive
words. Boon and Novikov secured the positiveness of admissible configurations in discs with the help of an
additional 'quadratic letter' (see \cite{R}, ch.12), but this involves a difficult control of
parameters for the constructed group. A new approach was suggested in \cite{SBR}. Invented by
Sapir, $S$-machines can work with non-positive words on the tapes. Here we also construct an
$S$-machine which is a somewhat modified composition of a convenient Turing machine with an 'adding machine' $Z(A)$
introduced in \cite{OS}. Fortunately, $Z(A)$ does not change the space of computations but
controls positiveness of configurations.

Recall that we should not just simulate the work of a machine but construct an embedding of
given group $H$ into a finitely presented group $G$ in the spirit of the Higman Embedding Theorem.
(See \cite{R}, ch.13). To this end we use a version of the two-disc scheme 
applied in \cite{OS2}. One can say more simply that the
configurations on the boundary of discs of the first type are longer than the words written on the 
boundary of corresponding discs of the second type, and the 
surpluses are relations of the group $H.$ Since every relation of $H$ holds in $G,$ we obtain  a 
homomorphism $H\to G$ that turns out to be injective.

To estimate the space function of the group $G$ we use van Kampen diagrams and
induct on the number of hubs. The basic problem is to cut a diagram $\Delta$ having at least two
hubs into two subdiagrams with hubs, so that
the perimeters of the subdiagrams do not exceed the perimeter of $\Delta.$ We have to introduce a new metric where the length of a word
depends on its syllable factorization. To find a shortcut we use the mirror symmetry of sectors
in discs, but unfortunately for the two-disc scheme, one of the sectors has no (mirror) copies, which creates
technical obstacles.\footnote{Note that one can give shorter proofs of Theorem \ref{realiz} and 
Corollary \ref{alpha} which do not need the two-disc scheme; but here we just obtain these results after Theorem \ref{main1} is proved.}
We study both exact and non-accurate copying for various types of complete and incomplete sectors.
A number of concepts, e.g. bands, trapezia, discs were introduced in algorithmic group theory 
early on, in papers \cite{SBR}, \cite{O}, \cite{BORS} and subsequent ones, and we reproduce
them in Section \ref{gd} (and partly in Sections \ref{mach} and \ref{compare}) in a modern form. Some others
(in particular, replica, unfinished diagram) are new.

In future work, we mean to consider similar problems for semigroups where the
simulating of machines is easier, and where the definition of space function
does not need  cyclic shifts and fragmentation of words, i.e.,
the space functions are $FL$-functions. For groups, an $FL$-analogue
of the results from this paper is still unachievable.

\section{Machines}\label{mach}

Suppose a $DTM$ $M$ of space complexity $S(n)$ solves the word problem in 
a finitely generated group $H.$ Then one can introduce a finitely
presented group $G$ whose relations simulate the work of $M$
and construct a natural homomorphism $H\to G.$ This will be done in Section $3.$
But as mentioned in the Introduction, one can prove neither that this homomorphism is injective nor that the space function of $G$ is bounded by $S(n),$ unless
$M$ enjoys some additional properties. Therefore  the machine $M$ 
has to be amended at the beginning. In this section we substitute $M$ for a non-deterministic 
S-machine $\cal S$ which is needed for the injectivity of the mapping $H\to G,$
and prove that $\cal S$ accepts the same language and has the same space 
complexity as the original $M.$ Moreover, the main Lemma \ref{MtocalS} of the
section claims that the generalized space complexity of $\cal S$ does not
exceed its space complexity; this property will be exploited in Subsection \ref{vse} 
in order to get an upper estimate of the space function of $G.$

\subsection{Definitions}\label{definitions}
We will use a model of recognizing $TM$-s which is close to the model from \cite{SBR}.
Recall that a {\it (multi-tape) $TM$ with $k$ tapes and $k$
heads} is  a tuple $$M= \langle X, Y, Q,
\Theta, \vec s_1, \vec s_0 \rangle$$ where $X$ is the input alphabet, 
$Y=\sqcup_{i=1}^k Y_i$ is the tape alphabet, $Y_1 \supset X,$
$Q=\sqcup_{i=1}^k Q_i$ is the set of states of the heads of the
machine, $\Theta$ is a set of transitions (commands), $\vec s_1$ is
the $k$-vector of start states, $\vec s_0$ is the $k$-vector of
accept states. ($\sqcup$ denotes the disjoint union.) The sets $X,Y, Q, \Theta$
are finite.

We
assume that the machine normally starts working with
states of the heads forming the vector $\vec s_1$,  with the head
placed at the right end of each tape, and accepts if it reaches the
state vector $\vec s_0$.  In general, the machine can be turned on
in any configuration and turned off at any time.

A {\em configuration} of tape number $i$ of a $TM$ is a word
$u q v$ where $q\in Q_i$ is the current state of the head,
$u$ is the word to the left of the head, and $v$ is the word to the
right of the head. 
A tape is {\em empty} if $u$, $v$ are empty words.

A {\em configuration} $U$ of the machine $M$ is a word 
$$\alpha_1U_1\omega_1\alpha_2U_2\omega_2... \alpha_kU_k\omega_k$$
where $U_i$ is the configuration of tape $i$, and $\alpha_i, \omega_i$ 
are special separating symbols.  For unification of notation, we
shall treat $\alpha_i, \omega_i$ as heads of the machine too.
These heads correspond to tapes that are always empty and do not
change during a computation.

An {\it input configuration} is a configuration where all tapes
except the first one are empty, the configuration of the first
tape (let us call it the {\it input tape}) is of the form $uq$,
$q\in Q_1$, $u$ is a word in the alphabet $X$, 
and the states form the start vector $\vec s_1$. The
{\em accept configuration} is the configuration where the state
vector is $\vec s_0$, the accept vector of the machine, and all
tapes are empty.  (The requirement that the tapes must be
empty is often removed for auxiliary machines which are used 
in construction of bigger machines.)  

To every $\theta\in \Theta,$ there corresponds a command (marked by the
same letter $\theta$), i.e., a pair of sequences of words
 $[V_1,...,V_k]$ and $[V'_1,...,V'_k]$ 
 such that for each $j\le k,$ either both $V_j$ and $V'_j$ are configurations
 of the tape number $j,$ or $V_j \equiv\alpha_j q$ and $V'_j \equiv \alpha_j q',$
 or $V_j \equiv q\omega_j$ and $V'_j \equiv q' \omega_j ,$ or $V_j \equiv \alpha_j q\omega_j$ and $V'_j \equiv \alpha_j q' \omega_j $ ($q,q'\in Q_j$ ).

In order to execute this command, the machine checks if $V_i$
is a subword of the configuration of tape $i$ for each $i\le k,$
and if this condition holds the machine replaces $V_i$ by $V'_i$
for all $i=1,\dots,k.$ Therefore we also use the notation:
$\theta: [V_1\to V'_1,\dots, V_k\to V'_k].$

 Suppose we have a sequence of configurations $w_0,...,w_t$ and a word
 $h\equiv \theta_1\dots\theta_{t}$ in the alphabet $\Theta,$
such that for every $i=1,..., t$ the machine passes from $w_{i-1}$ to
$w_i$ by applying the command $\theta_i$. Then the sequence
$(w_0\to w_1\to\dots\to w_t)$ is said to be a {\it computation with
history} $h.$ In this case we shall write $w_0\cdot h=w_t.$
 The number $t$ will be called the {\em time} or {\em length} of the computation.

A configuration $w$ is called {\em accepted} by a machine $M$ if
there exists at least one computation which starts with $w$ and
ends with the accept configuration. We do not only consider 
deterministic $TM$s, for example, we allow several 
transitions with the same left side. Moreover, for non-deterministic 
$TM$s, one may correspond identically equal 
executions to different symbols $\theta, \theta'\in \Theta$.

A word $u$ in the input alphabet $X$ is said to be {\em accepted} by the machine if the
corresponding input configuration is accepted.  The set of all
accepted words over the alphabet $X$ is called the {\em language  ${\cal L}_M$
recognized by the machine $M$}.

  Let $|w_i|_a$ ($i=0,...,t$) be the number of tape letters (or tape squares) in 
  the configuration $w_i$. (As in \cite{SBR}, the tape letters are called $a$-letters.) 
  Then the maximum of all $|w_i|_a$ will be called the {\em space
of computation} $C: w_0\to w_1\to\dots\to w_t$
and will be denoted by $space_M(C)$. By $space_M(w)$, we denote the minimal natural number $s$
such that there is an accepted computation of space at most $s,$ starting with the 
configuration $w.$  If $u\in {\cal L}_M,$ then, by definition, $space_M(u)$ is the space of the
corresponding input configuration $w.$

The number $S(n)=S_M(n)$ is the minimum of the numbers $space(u)$ over all
words $u\in {\cal L}_M,$ with $||u||\le n.$
The function $S(n)$ will be called the {\it space complexity} of the Turing
machine. 


The definition of the {\it generalized space complexity}
$S'(n)=S'_M(n) $ is  similar to the definition of  space complexity
but we consider arbitrary accepted
configurations $w$ with $|w|_a=n$, not just input
configurations as in the definition of  $S(n)$.  It is clear that
$S(n)\leq S'(n).$

To obtain the definitions of  $time_M(w)$, $time_M(u),$ {\it time complexity} $T_M(n)$ and
{\it generalized time complexity} $T'_M(n),$  one should replace 'space' by 'time' in the
previous definitions.

One may assume that only input configurations
involve the state letters from $\vec s_1$ and only one command is applicable to the
input configurations. (Indeed, given an NTM $M,$ one can add
additional states and add new commands changing the states from the new
vector $\vec s_1$.)
Similarly, one may assume that  there is a unique
accept configuration $\vec s_0$ and a unique accepting command.  
Under both of these assumptions, we will say that the
machine satisfies the $\vec s_{10}$-{\it condition}. These assumptions change
neither the language $\cal L$ nor the functions $S_{M}(n)$ and $S'_{M}(n).$

\subsection{Machines with equivalent space and generalized space complexities}\label{equal}
  
In this subsection, we construct an NTM $M_2$ which depends on 
an NTM $M_1,$ and prove Lemma \ref{M1M2} that allows us to replace
the original machine $M_1$ by a machine $M_2$ inheritting  the basic
characteristics of $M_1$ and having equivalent generalized space complexity
and space complexity.  

Let an NTM $M_1$ have $k$ tapes, and let its first tape  be the input tape. Then 
we add a tape numbered $k+1,$
which is empty for input configurations, and we
organize the work of the 3-stage machine $M_2$ as sequential work
of the following machines $M_{21}$, $M_{22},$ and $M_{23}.$ 

The machine $M_{21}$ uses only one
command $\theta_*$ that does not change states and adds one square
with an auxiliary letter  $*$ to the $(k+1)$-st tape, i.e the command $\theta_*$ has the form
$$[q_1\to q_1, \dots, q_k\to q_k, q_{k+1}\to *\;q_{k+1}]$$ The machine $M_{21}$ can execute this 
command arbitrarily many times while the tapes numbered $1,\dots,k$ keep 
the copy of an input configuration of $M_1$ unchanged.  Then  a connecting rule 
$\theta_{12}: [q_1\to q'_1, \dots, q_k\to q'_k, q_{k+1}\to q'_{k+1}]$ 
changes all states of the heads and switches on the machine $M_{22}.$

The work of $M_{22}$ on the tapes with numbers $1,\dots,k$ copies 
the work of $M_1.$ But the extension $\theta'$ of  every command $\theta$ of $M_1$ to the $(k+1)$-st tape is defined 
so that its application does not change the current 
space. More precisely,  if the application of $\theta$ inserts $m_1$ tape squares and deletes $m_2$
tape squares on $k$ tapes, then $\theta'$ inserts $m_2-m_1$ (deletes $m_1-m_2$) squares with 
letter $*$ on the $(k+1)$-st tape if $m_1 - m_2\le 0$ (if $m_1 - m_2\ge 0$).
That is the $(k+1)$-st component of $\theta'$ has the form $q'_{k+1}\to *^{m_2-m_1}q'_{k+1}$
(resp., $*^{m_1-m_2}q'_{k+1}\to q'_{k+1}$).
Note that one cannot apply $\theta'$ if $m_1-m_2$ exceeds the current number of
squares on the tape numbered $k+1.$ 

The connecting command $\theta_{23}$ is applicable when $M_{22}$ obtains 
the accept configuration on the first $k$ tapes. It changes the states and 
switches on the machine $M_{23}$ erasing all squares on the $(k+1)$-st tape (one by one), 
and then $M_2$ accepts. (The tape alphabet of $M_{23}$ has only one letter $*.$)

Let $w$ be a configuration of the machine $M_2$ such that $w\cdot \theta_*$ is defined, 
or such that $w$ is obtained after an application of the connecting command $\theta_{12}.$ Then we have
an input configuration on the tapes with numbers $1,\dots,k$ (plus several $*$-s on the $(k+1)$-st
tape). We will denote by $u(w)$ the input word $u$ written on the first tape.  It is an 
input word for the machine $M_1$ as well, and if it is accepted by $M_1$, the expression 
$space_{M_1} u(w)$ makes sense. 

The connecting commands $\theta_{12}$ and $\theta_{23}$ are not invertible in $M_2$
by definition. Therefore every non-empty accepting computation of $M_2$ has history of the
form $h_1h_2h_3,$ or $h_2h_3$, or $h_3$, where $h_l$ is the history for $M_{2l},$
($l=1,2,3$).
(To simplify notation
we attribute the command $\theta_{12}$ (the command $\theta_{23}$) to $h_2$ (to $h_3$).)

\begin{lemma} \label{M1M2}

 (a) The machines $M_1$ and $M_2$  recognize the same language $\cal L$. 
(b) The   space complexity $S_{M_2}(n)$ and the generalized space complexity $S'_{M_2}(n)$ of $M_2$ 
are both equivalent to $S_{M_1}(n)$. (c) If $w$ is an accepted configuration for $M_2,$ and
the command $\theta_*$ is applicable to $w,$ then $space_{M_2}(w)=\max (space_{M_1} (u(w)), |w|_a).$

\end{lemma}

\proof Assume that $u\in {\cal L}={\cal L}_{M_1}.$ Then $u\in {\cal L}' = {\cal L}_{M_2}$ because the machine $M_{21}$ can 
insert sufficiently many squares 
(equal to  $space_{M_1}(u)-||u||$) so that the accepting computation of $M_1$ can be simulated 
by $M_{22}.$ Also it is clear from the definition of $M_2$, that every accepting 
computation  for $M_2$ having a history $h_1h_2h_3$ as above, 
simulates, at stage 2, an accepting computation of $M_1$ with history $h_2.$ Therefore 
${\cal L}'=\cal L$ and $ S_{M_1}(n) = S_{M_2}(n).$

Assume now that $C: w=w_0\to\dots\to w_n$ is an accepting computation of $M_2$ with $space_{M_2}(C)=space_{M_2}(w)$ and  
$h\equiv h_1h_2h_3$ is the history with the above factorization ($h_1$ or $h_1h_2$ can be empty here).  If the word $h_1$ is empty, then $||w_0||\ge\dots\ge ||w_n||$ by
the definition of the machines $M_{22}$ and $M_{32}$. Hence the space of this computation
is equal to $|w|_a$  \footnote{ Here and in what follows we keep in mind that the difference $||w_i||-|w_i|_a$ is a constant
for any computation.}.
Then let $h_1$ be non-empty. It follows that the machine $M_2$
starts working with a copy of an input configuration of the machine $M_1,$ i.e., the input
tape of this configuration contains an input word $u=u(w),$ and the additional $(k+1)$-st
tape has $m$ squares for some $m\ge 0.$ Moreover  $u\in \cal L$ 
since the computation of $M_2$ is accepting. We consider two cases.

{\bf Case 1.} Suppose $m \ge space_{M_1}(u)-||u||.$ This inequality says that the additional tape
has enough squares to enable $M_{22}$  to simulate the accepting computation
of $M_1$ with the input word $u.$ Hence there is an $M_2$-computation  $w_0\to\dots\to w_{n'}$
with history of the form $h'_2h'_3$, and so its space, as well as the space of our original
accepting computation, is  $|w|_a.$

{\bf Case 2.} Suppose $m < space_{M_1}(u)-||u||.$ Then there is a computation $w_0\to\dots\to w_{n'}$
such that the commands of its $M_{21}$-stage insert squares until the total number of 
squares of the $(k+1)$-st tape becomes
equal to $space_{M_1}(u)-||u||,$ and then the machines $M_{21}$ and $M_{23}$ work in their standard manner.
The space of this (and the original) computation is $space_{M_1}(u).$

The estimates obtained in cases 1 and 2 prove statement (c) of the lemma. 
They also show that $$S_{M_1}(n)=S_{M_2}(n) \le S'_{M_2}(n)\le\max( S_{M_1}(n), n))\sim S_{M_1}(n) $$ 
and statement (b) is completely proved too.
\endproof

\subsection{Symmetric machines}\label{symm}

A simulation of the work of a machine $M$ by algebraic relations leads,
in fact, to the simulation of the machine $M^{sym}$ capable of inverting every command of $M$.
Therefore we must control the basic properties of such symmetrization, and this will be done
in Lemma \ref{generspace}.

For every command $\theta$ of a $TM$, given by the vector $[V_1\to V'_1,\dots,V_k\to V'_k]$,
the vector $[V'_1\to V_1,\dots,V'_k\to V_k]$ also gives  a command
of some $TM$. These two commands $\theta$ and $\theta^{-1}$ are called 
{\em mutually inverse}.

From now on we will assume that the machine $M_1$ we started with in Subsection \ref{equal} is a DTM and satisfies the $\vec s_{10}$-condition. 

Since the machine $M_1$ is deterministic, the machine $M_2$ has no invertible commands at all.
The definition of the {\it symmetric} machine $M_3= M_2^{sym}$ is the following. Suppose $M_2= \langle X, Y, Q,
\Theta, \vec s_1, \vec s_0 \rangle.$ 
Then by definition, $M_2^{sym} = \langle X, Y, Q,
\Theta^{sym}, \vec s_1, \vec s_0 \rangle,$ where $\Theta^{sym}$
is the minimal {\it symmetric} set containing $\Theta,$ that is,
with every command $[V_1\to V'_1,\dots, V_{k+1}\to V'_{k+1}]$ it contains the inverse command 
$[V'_1\to V_1,\dots, V'_{k+1}\to V_{k+1}]$; in other words,
$\Theta^{sym}=\Theta^+\sqcup\Theta^-,$ where 
$\Theta^+=\Theta$ and $\Theta^-=\{\theta^{-1}\mid \theta\in \Theta\}.$

A computation $w_0\to\dots\to w_t$ of $M_3$ (or of another machine) is called {\it reduced} if its history 
is a reduced word. If the history $h\equiv\theta_1\dots\theta_t$ contains a subword $\theta_{i}\theta_{i+1},$ where the commands
$\theta_{i}$ and $\theta_{i+1}$ are mutually inverse, then obviously there is a shorter computation
$w_0\to\dots \to w_{i-1}\equiv w_{i+1}\to\dots\to w_t$ whose space does not exceed the space of the original
computation.

\begin{lemma} \label{M1M3} 
Let $w\equiv w_0\to w_1\to\dots\to w_t$ be an accepted reduced computation of the machine $M_3$, and the command $\theta_*$
be applicable to $w.$ Then

(a) the word $u(w)$ belongs to the language $\cal L$ recognized by $M_1$ and

(b) the  space of this computation is at least $space_{M_1} (u(w)).$

\end{lemma}

\proof Let $h\equiv\theta(1)\dots\theta(t)$ be the history of the computation. If for some $i$, $w_{i+1}=w_i\cdot \theta(i+1)$ where $\theta(i+1)\equiv\theta_{23}$ or $\theta(i+1)^{\pm 1}$ is a command of $M_{23},$
then one can modify our accepted computation so that, for $j>i+1,$ every command $\theta(j)$ is a command of $M_{23}$
and $||w_j||\le ||w_{j-1}||.$ Hence we may assume that $h$  has exactly one letter
$\theta_{23}$ followed by the commands of $M_{23}$, and 
$h\equiv h_0\tau_1h_1\tau_2\dots\tau_s h_s,$ where  
$\tau_s\equiv\theta_{23},$ $\tau_i\equiv\theta_{12}^{(-1)^{s-i-1}}$ for $i<s,$ and the subwords $h_i$-s contain
no connecting commands. We may also assume that the subword $h_0$ is
empty since the command $\theta_*$ does not change the subword $u(w).$

Since, by the ${\vec s}_{10}$-condition, only one command of the machine $M_1$ (and of its analog $M_{22}$) accepts, the last command of $h_{s-1}$ is this unique command of $M_{22},$
and so this last command is positive.  
Therefore if $h_{s-1}$ 
contains a letter $\theta^{-1},$ where $\theta$ is a command of $M_{22},$  then $h$ has a 2-letter
subword $\theta_1^{-1}\theta_2$, where both  $\theta_1$ and $\theta_2$ are commands
of $M_{22}.$ Hence there is a configuration $w_i$  such that both
$\theta_1$ and $\theta_2$ are applicable to $w_i.$ This is impossible since the
machine $M_1$ is deterministic and the history $h_{s-1}$ is a reduced word.
Therefore $h_{s-1}$ is entirely the history of a computation of $M_{22},$ the computation
$w_{i_{s-1}}\to\dots\to w_t$ with history $\tau_{s-1}h_{s-1}\tau_s h_s \equiv \theta_{12}h_{s-1}\theta_{23} h_s$ 
is an accepting computation of $M_2,$ the word $u(w_{i_{s-1}})$ 
belongs to the language $\cal L,$ and the  space of the computation  $w_{i_{s-1}}\to\dots\to w_t$
is at least $space_{M_1} (u(w_{i_{s-1}}))$ by Lemma \ref{M1M2} (c).

Assume, by induction on $j,$ that $u(w_{i_{s-j}})$ belongs to $\cal L$ for $j\ge 1$, 
where the computation $w_{i_{s-j}}\to\dots\to w_t$  has history $\tau_{s-j}h_{s-j}\dots\tau_s h_s,$
and the  space of this computation is at least $space_{M_1} (u(w_{i_{s-j}})).$

Then the word $w_{i_{s-j-1}}$ has similar properties if $h_{s-j-1}$ consists of the commands of $M_{21}$ or
their inverses since these commands do not change the content of the tapes numbered $1,\dots, k.$
Otherwise the commands of $h_{s-j-1}$ are commands of $M_{22}$ (and their inverses), and since this 
machine is deterministic, the word $h_{s-j-1}$ has no subwords $\theta_1^{-1}\theta_2$ with positive $\theta_1$
and $\theta_2.$ Therefore we have $h_{s-j-1}\equiv g'g''^{-1},$ where both $g'$ and $g''$ are (positive) histories 
of $M_{22}$-computations. 
This implies the equality $ (w_{i_{s-j-1}}\cdot \tau_{s-j-1})\cdot g'=  w_{i_{s-j}}\cdot g''.$
Since the commands $\theta_{12}^{\pm 1}$  
do not change $u(w_i)$-s, we have $W_{i_{s-j-1}}\cdot g' = W_{i_{s-j}}\cdot g'',$ where $W_{i_{s-j-1}}$ and $W_{i_{s-j}}$ 
are the input configurations for the machine $M_1$ with inputs $u(w_{i_{s-j-1}})$ and $u(w_{i_{s-j}}),$ 
respectively.  (Here we use identical letters for the corresponding commands of $M_1$ and $M_{22}.$)

The machine $M_1$ is deterministic, and so the accepted computation for $W_{i_{s-j}}$  must look like
$W_{i_{s-j}} \to\dots \to W_{i_{s-j}}\cdot g''\to\dots $ , and consequently, the configuration 
$W_{i_{s-j}}\cdot g''$ is accepted by $M_1.$
Therefore we can construct the accepted computation $W_{i_{s-j-1}} \to\dots \to W_{i_{s-j-1}}\cdot g'
=W_{i_{s-j}}\cdot g'' \to\dots$ for $M_1$, and so the word $u(w_{i_{s-j-1}})$ belongs to $\cal L$,  
as desired. 

The constructed accepted computation of $M_1$ is decomposed in two parts. It follows from the definition
of $M_{22}$ that the  space of the first part
is majorized by the space of the $M_3$-computation $w_{i_{s-j-1}} \to\dots \to  w_{i_{s-j-1}}\cdot g'$
which is a part of the computation $w_{i_{s-j-1}}\to\dots\to w_t.$ The second part occurs 
in a deterministic accepted $M_1$-computation with input $u(w_{i_{s-j}})$, and so, by the inductive hypothesis, the  space of 
this part does not exceed the  space of the $M_3$-computation  $w_{i_{s-j}}\to\dots\to w_t.$ 
Hence $space_{M_1} (u(w_{i_{s-j-1}}))$ does not exceed the  space of the $M_3$-computation $w_{i_{s-j-1}}\to\dots\to w_t.$

Since $w \equiv w_{i_1},$ the lemma is proved by induction on $j.$ \endproof

\begin{lemma} \label{M3} The machines $M_1$ and $M_3$ recognize the same language. The generalized space complexities $S'_{M_2}(n)$ and $S'_{M_3}(n)$  are equivalent.
\end{lemma}

\proof We recall that every computation of the machine $M_2$ is also a computation of $M_3$. Therefore the first
statement follows from Lemmas \ref{M1M2} (a) and \ref{M1M3} (a). 

To prove the second part, it suffices to prove that for every accepted configuration $w$ of $M_2$ (of $M_3$), 
there is an accepted configuration $w'$ of $M_3$ (of $M_2$)
such that $||w'||\le ||w||$ but  $space_{M_3}(w')\ge space_{M_2}(w)$ (respectively,  $space_{M_2}(w')\ge space_{M_3}(w)$).

(1) Consider an accepted computation
$w\equiv w_0\to\dots\to w_t$ of $M_2$ whose  space is equal to $space_{M_2}(w).$ If the first command of this
computation is not a command of $M_{21},$ then $||w_0||\ge ||w_1||\ge\dots\ge ||w_t||,$ and therefore
$space_{M_2}(w) = |w|_a \le space_{M_3}(w)$, so one can choose $w'\equiv w.$ If the first command is a command
of $M_{21},$ then by Lemmas \ref{M1M2}(c) and \ref{M1M3} (b), we have 
$space_{M_2}(w) = \max (space_{M_1}(u(w)), |w|_a) \le space_{M_3}(w),$ and again $w' \equiv w.$

(2) Now consider a reduced accepted computation
$w\equiv w_0\to\dots\to w_t$ of $M_3$ whose  space is equal to $space_{M_3}(w).$ If the first command (or its inverse)
is a command of $M_{23},$ then the commands of the shortest accepted computation with minimal  space just erase squares.
Hence $space_{M_3}(w) = |w|_a=space_{M_2}(w)$, and we can choose $w'$ equal to $w.$ If the first command is
a command of $M_{21}^{sym}$, then by Lemma  \ref{M1M3} (a), the word $w$ is accepted by $M_2.$ Since every command of
$M_2$ is a command of $M_3$, we have
$space_{M_3}(w) \le space_{M_2}(w),$ and again it suffices to set $w'\equiv w.$

Thus, we may assume that the first command of our computation (or its inverse) is a command of $M_{22}.$ Therefore
the history $h$ of this computation has a prefix $h'h''$ with non-empty $h'',$ where every command of $h'$ 
is a command of $M_{22}^{sym}$ and either every command of $h''$ (or the inverse command) is a command of $M_{21}$ or every command of $h''$ 
is a command of $M_{23}.$ In the latter case, we may assume that $h\equiv h'h''$ and $||w_0||\ge ||w_1||\ge\dots \ge ||w_t||,$
and so $space_{M_3}(w) = |w|_a\le space_{M_2}(w),$ and $w'\equiv w.$ In the former case, we set $w'=w\cdot h'$ and note
that $||w||=||w_1||=\dots = ||w'||$ since the commands of the computation $w\to\dots\to w'$ with history $h'$ do not change 
the number of tape squares. In particular, we have $space_{M_3}(w) \le space_{M_3}(w').$ Since the command $\theta_*$ is 
applicable to $w'$ in this case,
we have   $space_{M_3}(w') \le space_{M_2}(w')$ as was observed in the previous paragraph. 
Therefore $space_{M_3}(w) \le space_{M_3} (w')\le space_{M_2}(w'),$ as desired; and the lemma is proved.
   
 \endproof  

\begin{lemma} \label{generspace}
For every DTM $M$ recognizing a language $\cal L$ and having space complexity $S(n)$ (for every NTM $M$ 
recognizing a language $\cal L$ and having FSC space
complexity $f(n)$), there exists an NTM $M'$ with the following properties.
\begin{enumerate}

 \item The machine $M'$ recognizes the language $\cal L$. 

\item $M'$ is symmetric. 

\item The   space and the generalized space complexities of $M'$ are equivalent
 to $S(n)$ (respectively, are equivalent to $f(n)^2$).
 
\item For every command  $[V_1\to V'_1,\dots, V_k\to V'_k]$ of $M'$, we have $\sum |V_i|_a +\sum |V'_i|_a\le 1,$ i.e., at most one tape letter is involved in the command.  
 \end{enumerate}
 \end{lemma}
   
\proof Assume that the machine $M$ is deterministic. Then starting with $M_1=M$,
we construct the machine $M_2$ described in  Subsection \ref{equal} and the machine
$M_3=M_2^{sym}.$ Now Lemma \ref{M3} implies statement 1,  and statement
2 is true since $M_3=M_2^{sym}.$ Then, by Lemmas \ref{M3} and \ref{M1M2} (b),  we have
$$S'_{M_3}(n) \sim S'_{M_2}(n)\sim S_{M_2}(n)\sim S_{M_1}(n)=S(n),$$ and $S(n) \le S_{M_3}(n)\le S'_{M_3}(n)$ by Lemma \ref{M1M3} (b), and so all these functions are equivalent.
Finally, we modify $M_3$ to obtain property $4$. For example, if for a one-tape machine we have  a command
$aq\to bq',$ then we introduce a new state letter $q''$ and replace this command by
two commands $aq\to q''$ and $q''\to bq'.$ It is easy to see that the obtained
machine $M'$ satisfies property $4$ and inherits properties $1 - 3$ from $M_3.$

If $M$ is non-deterministic, then we first use  that the
function $f(n)$ is $FSC,$ and therefore, by Savitch's theorem
(\cite{DK}, Theorem 1.30), there exists a $DTM$ $M_1$ accepting
the same language $L$ with space complexity equivalent to $f(n)^2$.
So the replacement of $S(n)$ by $f(n)^2$ in the previous paragraph provides the proof of the
non-deterministic version of the lemma. \endproof

\subsection{S-machines}\label{Smach}

By Lemma \ref{simul}, computations of a machine are faithfully represented
by  special van Kampen diagrams (trapezia) over the group $G$. However this statement
is true for $S$-machines but it is false  for standard $TM$-s.  

Ordinary Turing machines work with positive words and they can see some letters on the tape
near the position where the head is. The
command executed by the machine depends not only on the state of the head but also on the letter(s)
observed by the head. In contrast, S-machines introduced in [SBR] work with words in group alphabets
and they are almost  ``blind'', i.e., the heads do not observe the tape letters. But the heads can ``see''
each other if there are no tape letters between them. We will use the following precise definition
of an $S$-machine  as a rewriting system.

Let $k$ be a natural number. Consider a {\it language of admissible words}. It
consists of words of the form $$q_1u_1q_2\dots u_k q_{k+1},$$ where  $q_i$ are letters
from disjoint sets $Q_i$ ($i=1,\dots,k+1$),
$u_i$ are reduced words in a group alphabet $Y_i,$
(i.e. every letter
$a$ belongs to it iff the inverse letter $a^{-1}$ does) and
the sets $Y=\sqcup Y_i$ and $Q=\sqcup Q_i$ 
are finite.
The letters from $Q$ are
called {\it state} letters, and the letters from $Y$ are {\it tape} letters.
Notice that in every admissible word, there is exactly one representative
of each $Q_i$ and these representatives appear in this word in the order
of the indicees of $Q_i.$ (i.e., unlike \cite{OS}, we consider
only the regular order of $Q_i$-s in admissible words).

There is a finite set of {\it commands} (or {\it rules}) $\Theta.$
To every $\theta\in \Theta$, we associate two sequences of reduced
words from the free group $F(Q\cup Y)$:
$B(\theta)=[U_1,...,U_{k+1}]$, $T(\theta)=[V_1,...,V_{k+1}]$, and
a subset $Y(\theta)=\sqcup Y_i(\theta)$ of $Y$, where
$Y_i(\theta)\subseteq Y_i$.

The words $U_i, V_i$ satisfy the following restriction:

\begin{itemize}
\item[(*)] For every $i=1,...,k+1$, the words $U_i$ and $V_i$ have the form
$$U_i\equiv v_{i-1}q_iu_i, \quad V_i\equiv v_{i-1}'q_{i}'u_{i}'$$ where $q_{i}, q_{i}'\in
Q_{i}$, $u_{i}$ and $u_{i}'$ are words in the alphabet $Y_i(\theta)$, $v_{i-1}$ and $v_{i-1}'$ are words in the alphabet
$Y_{i-1}(\theta)$. The words $v_0, v'_0, u_{k+1}, u'_{k+1}$ are empty.

\end{itemize}

Sometimes we will denote the rule $\theta$ by
$[U_1\to V_1,...,U_{k+1}\to V_{k+1}]$. This notation contains no information 
about the sets
$Y_i(\theta)$. In most cases it will be clear what these sets are.
In the $S$-machines used in this paper, the sets $Y_i(\theta)$
will mostly be  equal to either $Y_i$ or $\emptyset$. By default
$Y_i(\theta)=Y_i$.

We will use the notation
$v_iq_iu_i\tool v_i'q_i'u_i'$ for a part of a rule when the
corresponding $Y_i(\theta)$ is empty (a similar notation has been
used in \cite{OS3}). In particular, $u_i$ and $u'_i$ are empty
words in this case.

Every $S$-rule $\theta=[U_1\to V_1,...,U_{k+1}\to V_{k+1}]$ has an
inverse $\theta\iv=[V_1\to U_1,...,V_{k+1}\to U_{k+1}]$; we set
$Y_i(\theta\iv)=Y_i(\theta)$. We always divide the set of rules
$\Theta$ of an $S$-machine into two disjoint parts, $\Theta^+$ and
$\Theta^-$ such that for every $\theta\in \Theta^+$, $\theta\iv\in
\Theta^-$ and for every $\theta\in\Theta^-$,
$\theta\iv\in\Theta^+$. The rules from $\Theta^+$ (resp.
$\Theta^-$) are called {\em positive} (resp. {\em negative}).

To apply an $S$-rule $\theta$ to an admissible word
$W\equiv q_1w_1q_2\dots w_k q_{k+1}$
means to check if every $w_i$ is a word in the alphabet $Y_i(\theta)$ and then, if $W$ satisfies
 this condition, to  simultaneously replace subwords $U_i$ by subwords $V_i$ ($i=1,\dots,k+1$).
 This replacement can be performed in the form $q_i \to
 v'_{i-1}v_i^{-1}q'_iu_i^{-1}u'_{i}$ followed by reducing 
 the resulted word. The following convention is important in the definition of
$S$-machine: {\it After every application of a rewriting rule, the word is automatically 
reduced. The reducing is not considered a separate step of an $S$-machine.}

   The definitions of computation, its history, input admissible words, the accept word,
   the language of admissible words,  space of a computation,  space and generalized complexities, time and generalized time complexities of an $S$-machine are similar to those for a $TM$. (One should
   replace the word ``configuration" by ``admissible word" in the definitions.)

   Although $S$-machines are usually highly non-deterministic, they are better adapted for
   simulating by finitely presented groups than
   ordinary $TM$s.  On the other hand, it is mentioned in \cite{SBR} that 
   every symmetric NTM  $M$ can be viewed as an $S$-machine $S(M)$: just interpret the commands of the Turing 
   machine as $S$-rules. (For example, the part of a rule of the form $aq\to bq'$ is interpreted as
   $q\to a^{-1}bq'$, and the part of the form $\alpha_j q\to \alpha_j q'$ is interpreted by the pair
   $\alpha_j\tool\alpha_j$, $q\to q'.$)
   Unfortunately, the language recognized by $S(M)$ is in general much bigger than the 
   language recognized by $M$ since $M$ works with a {\it positive} tape alphabet only. 
   Nevertheless the following statement is true:

\begin{lemma} \label{tmtostm} (Compare with Prop. 4.1\cite{SBR}.) Every computation of a symmetric $NTM$ $M$ is a computation of  $S(M)$ with the same history. If $M$ satisfies property 
4 from Lemma \ref{generspace}, then every {\em positive} computation of $S(M)$, i.e.,  
a computation consisting of positive words, is a computation of $M$ with the same history.
\end{lemma}

\proof The first statement follows from the definition of the machine $S(M).$

Every positive admissible word $W$ of $S(M)$  is a configuration of the Turing machine $M$. Assume that a rule  $\bar\theta$ of $S(M)$ corresponding to a command $\theta$ of $M$ is applicable
to this word $W$ and the word $W\cdot \bar\theta$ is positive. Recall that by property 4 of Lemma \ref{generspace},  $\theta$ involves at most one tape letter (e.g., it cannot replace a
tape letter by a tape letter or have a part of the form $aq\to aq'$). Therefore the positiveness of both $W$ and $W\cdot \bar\theta$
implies that the application of $\bar\theta$ just coincides with the application of $\theta.$
The second statement of the lemma follows.
\endproof

\subsection{Composition with an adding machine}\label{compos}

One the one hand, $S$-machines much better fit for simulating their work
by group relations than ordinary Turing machines (and moreover, $S$-machines are themselves treated in \cite{OS} as HNN-extensions of a free group with basis $Y\cup Q$). On the other
hand, when a symmetric Turing mashine $M$ starts working as an S-machine (see
Subsection \ref{Smach}), the language of (positive) accepted words can enlarge uncontrolably.  Therefore we will consider a composition $\cal S$ of a symmetric Turing machine $M$
and an auxiliarly S-machine $Z(A)$ from \cite{OS}. Fortunately, this composition is an S-machine
which inherits the language, the space function and the generalized space function of $M.$

In \cite{OS}, the main duty of $Z(A)$
was the exponential slowing down of basic computations, while now we will mainly use the capacity of $Z(A)$
(observed in Lemma 3.25 (2) \cite{OS})  to check whether an admissible word is positive or not. 

The tape alphabet of $Z(A)$ consists of an alphabet $A^{\pm 1}$ and two copies $A_1^{\pm 1}$ and $A_2^{\pm 1}$ of $A^{\pm 1}$
while the input alphabet is $A^{\pm 1}$. The admissible words are of the form $LupvR$, where $u$ is a reduced
word in the alphabet $A^{\pm 1} \cup A_1^{\pm 1},$ $v$ is a reduced word in $A_2^{\pm 1}$, the symbols $L, p, R,$ 
are state letters, the commands do not change $L$ and $R$, and $p\in \{p(1),p(2),p(3)\}$. The input
words have the form $Lup(1)R$ and the accept ones are of the form $Lup(3)R.$ 

For every
letter $a\in A$ let $a_1, a_2$ denote its copies in $A_1, A_2$.


Thus, the set of state letters is $Q_1\cup Q_2\cup Q_3$ where $$Q_1=\{L\},
Q_2=\{p(1),p(2), p(3)\}, Q_3=\{R\}.$$ The set of tape letters is
$Y_1\cup Y_2$ where $Y_1=A\cup A_1$ and $Y_2=A_2$.

The machine $Z(A)$ has the following positive rules
($a$ is an
arbitrary letter from $A$ below). The following comments explain the meanings of
these rules.

\begin{itemize}

\item $r_1(a)=[L\to L, p(1)\to a_1\iv p(1)a_2, R\to R]$.
\me

{\em Comment.} The state letter $p(1)$ moves left searching for a
letter from $A$ and replacing letters from $A_1$ by their copies
in $A_2$.

\me

\item $r_{12}(a)=[L\to L, p(1)\to a\iv a_1p(2), R\to R]$.

{\em Comment.} When the first letter $a$ of $A$ is found, it is
replaced by $a_1$, and $p(1)$ turns into $p(2)$.

\me

\item $r_2(a)=[L\to L, p(2)\to ap(2)a_2\iv, R\to R]$.

\me

{\em Comment.} The state letter $p(2)$ moves toward $R$.

\me

\item $r_{21}=[L\to L, p(2)\tool p(1), R\to R]$, $Y_1(r_{21})=Y_1,
Y_2(r_{21})=\emptyset$.

\me

{\em Comment.} $p(2)$ and $R$ meet, and the cycle starts again.

\me

\item $r_{13}=[L\tool L, p(1)\to p(3), R\to R]$, $Y_1(r_{13})=\emptyset,
Y_2(r_{13})=A_2$.

\me

{\em Comment.} If $p(1)$ never finds a letter from $A$, then the cycle
ends, and $p(1)$ turns into $p(3)$; $p$ and $L$ must stay next to each
other in order for this rule to be executable.

\item $r_{3}(a)=[L\to L, p(3)\to ap(3)a_2\iv, R\to R]$,
$Y_1(r_3(a))=A, Y_2(r_3(a))=A_2$

 \me

{\em Comment.} The letter $p(3)$ returns to $R$.

\end{itemize}


Here we will not explicitly use the form of these rules and will rather formulate, 
in Lemma \ref{Z}, the required properties of $Z(A)$ obtained in \cite{OS}. 
Nevertheless the reader can see how this simple $S$-machine works if it `agrees' 
to apply the rules in the order recommended in the comments. In particular, if the input word 
is $La^np(1)R$ ($n>0$) for some $a\in A,$ then every maximal
subcomputation ending with $r_{21}$ and having no other rules $r_{21}$ in its history,
changes the word over $Y_1$ in the manner that
one changes a binary $n$-digit number when one adds $1$ to it. (Therefore the
machine $Z(A)$ is called {\it adding}.) However Lemma \ref{Z} is not so obvious
since (1) the machine $Z(A)$ is non-deterministic and (2) it works with arbitrary reduced tape
words (while Turing machines deal with positive words only).

If $w$ a word in the alphabet $A\cup A_1^{\pm 1} \cup A_2^{\pm 1},$ then its {\it projection} onto $A^{\pm 1}$
takes every letter to its copy in $A.$

\begin{lemma} \label{Z} The following properties of the machine $Z(A)$ hold.

(1) Every positive input word $u$ in the alphabet $A$ is accepted by a canonical computation of
$Z(A)$ with a positive history and such that all words appearing in this computation
have equal lengths. 

(2) For every computation $LupvR\equiv w\to\dots\to w'\equiv Lu'p'v'R$ of $Z(A),$   
the projections of the words $uv$ and  $u'v'$ onto $A$ are freely equal. In particular,
$u\equiv u'$ if the words $v$ and $v'$ are empty and the words $u$ and $u'$ contain no letters
from $A_1^{\pm 1}$.

(3) If $w_0\to\dots\to w_t$ is a reduced computation of $Z(A)$ and $||w_0||<||w_1||,$
then $||w_1||\le ||w_2||\le\dots\le ||w_t||.$

(4) For every reduced computation $w_0\to\dots\to w_t,$ we have $||w_i\|\le \max(||w_0||,||w_t||)$ ($i=0,\dots,t$).

(5) If $w\equiv LupR,$ where $p\equiv p(1)$ (or $p\equiv p(3)$), $w\equiv w_0\to\dots\to w_t$ is a reduced computation, $w_t$
contains the subword $p(3)R$ (respectively, $p(1)R$), and all $a$-letters of $w_0$ and $w_t$ are from $A^{\pm 1}$, then $u$ is a positive word and all words of the computation have the same length. 
The length of this computation is at least $2^{||u||}.$ 

(6) There is no reduced computation $w\equiv w_0\to\dots\to w_t$ of length $t\ge 1$ such that both $w_0$ and $w_t$
contain $p(1)R$ or both of them contain $p(3)R$ and all $a$-letters of $w_0$ and $w_t$ belong to $A^{\pm 1}.$ 

\end{lemma}

\proof (1) To construct the canonical computation one should just use the above comments.

(2) This claim  is the statement of Lemma 3.18 of \cite{OS}.

(3) This statement is true by Lemma 3.24 of \cite{OS} 

(4), (5) These statements are contained in Lemma 3.25 of \cite{OS}.

(6) This statement is contained in Lemma 3.27 of \cite{OS}.

\endproof

Roughly speaking, the composition of machines $M$ and $Z(A)$ was constructed in \cite{OS} as follows.
After an application of (the copy of) a positive command from $M,$ the copies of
$Z(A)$ successivly implement their accepting computations on each of the tapes,
and then (the copy of) the next command of $M$ is applied. Thus $Z(A)$ just (exponentially) slows down the computations of $M.$  

Here we somewhat modify the definition of the composition of a symmetric Turing machine $M$ and the
adding machine $Z(A).$ The difference is that  machines of the form $Z(A)$ 
will  work not only after the application of every command of $M$ but before applications of commands from $M$ as well. This
makes it possible to simulate the work of any symmetric NTM, not only $S$-machines as  was done in \cite{OS}. 
The aim of following construction is to obtain an $S$-machine $\cal S$ which
recognizes the same language and has the same  space and generalized space complexity
as the symmetric Turing machine $M.$  The $S$-machine constructed in [SBR] cannot
serve in the present paper since the  space and the generalized space complexities of that
machine are equivalent to its time complexity.

Consider a symmetric NTM  $M= \langle X, Y, Q,
\Theta, \vec s_1, \vec s_0 \rangle$ with $Y=\sqcup_{i=1}^l Y_i,$ and with the ${\vec s}_{10}$-condition. 
The set $\Theta$ is a disjoint union of positive and negative commands:
$\Theta=\Theta^+\sqcup\Theta^-$. Let $S(M)$ be the associated $S$-machine defined
before Lemma \ref{tmtostm}. We will assume that the admissible words of $S(M)$ are
of the form $k_1u_1k_2\dots u_l k_{l+1},$ where  $k_i$ are letters
from disjoint sets $Q_i$ ($i=1,\dots,l+1$), $u_i$ is a reduced word in the alphabet $Y_i$ ($i\le l$).

To define the composition  ${\cal S} =M \circ Z$ of $M$ and
$Z(A),$ we will insert a $p$-letter between any two
consecutive $q$-letters $k_i, k_{i+1}$ in an admissible word of $S(M)$, so that to treat 
any subword $k_i...p...k_{i+1}$ as an admissible word for a copy of $Z(A)$.
In particular, the {\it start} and the {\it accept} words of $\cal S$ are obtained
from, respectively, the start and accept words of $M.$

First, for every $i=1,...,l$, we make two copies of the alphabet
$Y_i$ of $S(M)$ ($i=1,...,l$): $Y_{i,1}$ and $Y_{i,2}$. The
set of state letters of the new machine is $$K_1\sqcup P_1\sqcup
K_2\sqcup P_2\sqcup...\sqcup P_{l}\sqcup K_{l+1},$$ where $P_i=\{p_i, p_i(\theta,1^-), p_i(\theta,2^-), p_i(\theta,3^-),
p_i(\theta,1^+), p_i(\theta,2^+), p_i(\theta,3^+)\mid
\theta\in\Theta^+\}$, $i=1,...,l$.

The set of tape letters is $$\bar Y=(Y_{1}\sqcup Y_{1,1})\sqcup
Y_{1,2}\sqcup (Y_{2}\sqcup Y_{2,1})\sqcup Y_{2,2}\sqcup...\sqcup
(Y_{l}\sqcup Y_{l,1})\sqcup Y_{l,2};$$ the components of this
union will be denoted by $\bar Y_1,...,\bar Y_{2l}$.

The set of positive rules $\bar\Theta$ of $M\circ Z$ is a union
of the set of modified positive rules of $S(M)$ and
of the positive rules of $Z_i(\theta,-)$ and $Z_i(\theta,+)$ ($\theta\in\Theta,
i=1,...,l$) which are copies of the machines $Z(Y_i)$ (also suitably
modified).

More precisely, suppose a positive command $\theta$ of $S(M)$
differs from the unique start and accept commands of $S(M)$ and has 
 the form
$$[k_1u_1\to k_1'u_1', v_1k_2u_2\to
v_1'k_2'u_2',...,v_{l}k_{l+1}\to v_{l}'k_{l+1}']$$ where $k_i, k_i'\in
K_i$, $u_i$ and $v_i$ are words in $Y_i$. Then its copy in $M\circ Z$  is 

$$\bar\theta=\begin{array}{l}[k_1u_1\to k_1'u_1', v_1p_1(\theta, 3^-)\tool
v_1'p_1(\theta,1^+), k_2u_2\to k_2'u_2', ...,\\
 v_{l}p_{l}(\theta, 3^-)\tool v_l'p_{l}(\theta,1^+), k_{l+1}\to
 k_{l+1}']\end{array}$$ with
$\bar Y_{2i-1}(\bar\theta)=Y_{i}(\theta)$ and $\bar Y_{2i}(\theta)=\emptyset$
for every $i$, in particular, the words $u_i, v_i$ are rewritten here in the alphabet $Y_{i}.$


Each machine $Z_i(\theta,-)$ is a copy of the machine $Z(Y_{i}),$ where
every rule $\tau=[U_1\to V_1, U_2\to V_2, U_3\to V_3]$ is replaced
by rule of the form
$$\bar\tau_i(\theta,-)=\left[
\begin{array}{l}\bar U_1\to \bar V_1, \bar U_2\to \bar V_2, \bar U_3\to \bar V_3,\\
k_j\to k_j, p_j(\theta,3^-)\tool  p_j(\theta, 3^-),
j=1,...,i-1,\\
p_s(\theta,1^-)\tool p_s(\theta,1^-), k_{s+1}\to k_{s+1},
s=i+1,...,l\end{array}\right]$$ where $\bar U_1, \bar U_2, \bar
U_3, \bar V_1, \bar V_2, \bar V_3$ are obtained from $U_1, U_2, U_3,
V_1, V_2, V_3$, respectively, by replacing $p(j)$ with
$p_i(\theta,j^-)$, $L$ with $k_{i}$ and $R$ with $k_{i+1}$, and for
$s\ne i$, $\bar Y_{2s-1}(\bar\tau_i(\theta,-))=Y_{i}$.

Similarly, each machine $Z_i(\theta,+)$ is a copy of the machine $Z(Y_{i}),$ where
every rule $\tau=[U_1\to V_1, U_2\to V_2, U_3\to V_3]$ is replaced
by the rule of the form
$$\bar\tau_i(\theta,+)=\left[
\begin{array}{l}\bar U_1\to \bar V_1, \bar U_2\to \bar V_2, \bar U_3\to \bar V_3,\\
k_j'\to k_j', p_j(\theta,3^+)\tool  p_j(\theta, 3^+),
j=1,...,i-1,\\
p_s(\theta,1^+)\tool p_s(\theta,1^+), k_{s+1}'\to k_{s+1}',
s=i+1,...,l\end{array}\right]$$ where $\bar U_1, \bar U_2, \bar
U_3, \bar V_1, \bar V_2, \bar V_3$ are obtained from $U_1, U_2, U_3,
V_1, V_2, V_3$, respectively, by replacing $p(j)$ with
$p_i(\theta,j^+)$, $L$ with $k_{i}'$ and $R$ with $k_{i+1}'$, and for
$s\ne i$, $\bar Y_{2s-1}(\bar\tau_i(\theta,+))=Y_{i}$.

If $\theta$ is the start (the accept) command of $S(M),$ then we introduce only
machines $Z_i(\theta,+)$ (only $Z_i(\theta,-)$, respectively), and replace the
letters $p_j(\theta,3^-)$ (replace $p_j(\theta, 1^+)$, resp.) by $p_j$ in the above
definition of the command $\bar\theta.$

In addition, we need the following {\em transition} rules
$\zeta(\theta,-)$ and $\zeta(\theta,+)$ 
that transform all $p$-letters from and to their original forms.

$$[k_i\to k_i, p_j\tool
p_j(\theta, 1^-), i=1,...,l+1, j=1,...,l].$$

$$[k_i'\to k_i', p_j(\theta,3^+)\tool
p_j, i=1,...,l+1, j=1,...,l].$$

If $\theta$ is the start (the accept) command of $S(M),$ then we introduce only $\zeta(\theta,+)$
(only $\zeta(\theta,-)$, respectively).


Thus while the machine $Z_i(\theta, -)$ (the machine $Z_i(\theta, +)$) works
all other machines $Z_j(\theta, -)$  (all machines $Z_j(\theta, +)$, $j\ne i$) 
must stay idle (their state letters do not
change and do not move away from the corresponding $k$-letters).
After the machine $Z_i(\theta,-)$ (the machine $Z_i(\theta, +)$) finishes, i.e., 
the state letter $p_{i}(\theta,3^{\pm})$ appears next to $k_{i+1}$ (next to $k'_{i+1}$), 
the next machine
$Z_{i+1}(\theta, -)$ (the machine $Z_{i+1}(\theta, +)$) starts working. 
The transition rule $\zeta(\theta,-)$ switches on the consecutive works of the
machines $Z_1(\theta, -),\dots, Z_{l}(\theta, -).$ After all $p$-letters have the
form $p_j(\theta,3^-)$ we can apply the rule $\bar\theta$  and turn all
$p_j(\theta,3^-)$ into $p_j(\theta, 1^+)$. This switches on the consecutive work of
$Z_1(\theta, +),\dots, Z_{l}(\theta, +),$ followed by the transition rule $\zeta(\theta,+).$

So the modified rule  $\bar\theta\in\bar\Theta$  turns on  the machines $Z_1(\theta,+),\dots, Z_l(\theta,+)$ (and $\bar\theta$
 follows after the machines $Z_1(\theta,-),\dots, Z_l(\theta,-)$ finish their work).

Thus, in order to simulate a computation of the symmetric $TM$ $M$  (and of the $S$-machine $S(M)$)
consisting of a sequence of applications of rules $\theta_1,
\theta_2,...,\theta_s$, we first apply all rules corresponding to
$\theta_1$ (before and after $\bar\theta$), then all rules corresponding to $\theta_2$, then all
rules corresponding to $\theta_3$, etc.  The language ${\cal L}_{\cal S}$ of $\cal S$
consists of some words $u$ in the alphabet $Y_{1}.$ In particular, every
admissible input word of $\cal S$ is of the form $\Sigma(u)\equiv k_1up_1k_2p_2k_3\dots k_{l-1}p_lk_l.$
We  denote the accept word of $\cal S$ by $\Sigma_0=\Sigma_0({\cal S}).$ 

The modified rules $\bar\theta$ of $S(M)$ will be called {\em basic
rules} of the $S$-machine $\cal S.$

\subsection{Computations of machine $\cal S$}\label{comput}

Lemma \ref{MtocalS} will complete the part of our paper dealing with machines. This lemma makes possible
to substitute the Turing machine $M'$ from Lemma \ref{generspace} by the S-machine $\cal S.$
Lemma \ref{spacecalS} (4) helps to obtain the lower bound for the space function of the group
$G$ in Subsection \ref{vse}.

 Given a computation $C$ 
 of $\cal S,$ 
 one obtains a computation $C_{S(M)}$  of  $S(M)$ after omitting all the (copies of) commands of machines  $Z_i(\theta,\pm)$ and deleting the additional state letters (letters of the sets $P_1,\dots, P_l$) from the admissible words of the computation $C.$ 
 (By Lemma \ref{Z} (2),  the computation $C_{S(M)}$ is well-defined, although $C_{S(M)}$ is of length $0$ if $C$ has no basic rules.) 
 
 \begin{lemma} \label{calS} 
 
 (1) The start and the accept rules of $\cal S$ are basic rules which are the copies of the start and the
accept commands of $M,$ respectively. The machine $M\circ Z$ satisfies the ${\vec s}_{10}$-condition.

  (2) If the history $h$ of a reduced computation $C: \;\; w_0\to\dots\to w_t$ of $\cal S$ contains no
  basic commands, then $||w_i||\le \max (||w_0||,||w_t||)$ for every $i$ ($0\le i \le t$). If the only
  basic letter of $h$ is the last one, then $||w_i||\le ||w_0||$ for $i<t.$

  (3) If a computation $C$ of $\cal S$ is reduced then $C_{S(M)}$ is also reduced.

 (4) For any reduced computation $C: \;\; w_0\to\dots\to w_t$ of $\cal S$ starting and ending
 with basic rules we have 
$space_{\cal S}(C) = space_{S(M)}(C_{S(M)}).$

 (5) For every positive reduced computation $C: \;\; w_0\to\dots\to w_t$ of the machine $S(M),$ there is a canonical 
 reduced computation $C_{\cal S}$ of $\cal S$ whose history starts and ends with basic commands, 
 such that $(C_{\cal S})_{S(M)} = C.$ Moreover, we have $space_{\cal S}(C_{\cal S}) = space_{S(M)}(C).$
 
 \end{lemma}
 
 \proof (1) Property (1) follows from the similar property of the machine $M$ and from the definition 
 of the machine $M\circ Z.$  

(2) We start with the first statement. If the history of the whole computation consists of the commands of one machine $Z_i(\theta,\pm)$, then the statement
follows from Lemma \ref{Z}(3). Otherwise it has $s\ge 1$ admissible words $w_{i_1},\dots, w_{i_s},$ such that
the history of every subcomputation 
$$C_0: w_0\to\dots\to w_{i_1},\dots, C_j: w_{i_j}\to\dots\to w_{i_{j+1}},\dots,
C_s: w_{i_s}\to\dots\to w_{t}$$ either  consists of the commands of some  $Z_i(\theta,\pm),$ being a maximal
subcomputation with this property, or is equal to a transition letter $\zeta(\theta,\pm)^{\pm 1}$ for some
$\theta$. In the former case, all the admissible words participating in $C_j$ have the same length 
for $j\in [1, s-1]$ by Lemma \ref{Z} (5,6). The same is clearly true in the latter case. Therefore it 
suffices to prove that the lengths of the admissible words do not decrease in the subcomputations 
$C_0^{-1}$ and $C_s.$  

Let us consider $C_s$ only, assuming that it is a computation of some machine $Z_i(\theta,\pm).$
It corresponds to a computation $LupvR\equiv W_0 \to\dots\to W_m \equiv Lu'pv'R$ of a machine of the form $Z(A)$ 
with $m=t-i_s.$  Since $s\ge 1$ and $C_s$ is a maximal subcomputation corresponding to $Z_i(\theta,\pm),$
we must have $||v||=0$. (Otherwise only commands of $Z_i(\theta,\pm)$ could be applied to $w_{i_s}$,
and so $w_{i_s-1}\to w_{i_s}\to\dots\to w_t$ is a longer computation of the same $Z_i(\theta,\pm).$)
Hence the projection of the word $uv$ onto $A$ is reducible, and so $||W_0||\le ||W_i||$ for every $i\in [1,m]$
by Lemma \ref{Z} (2). Therefore either all $W_i$-s have equal lengths or $||W_0||\le ||W_1||\le\dots \le ||W_m||$
by Lemma \ref{Z} (3). In any case, we have $||W_0||\le\dots\le ||W_m||.$ This implies $||w_{i_s}||\le\dots\le ||w_t||,$
as required.  

The proof of the second claim is similar: The computation $w_{i_1}\to\dots\to w_{i_s}$ is a product of subcomputations
$C_j$-s which preserve the lengths of $w_i$-s, while $C_0$ and $C_s$ cannot increase the space.

(3)Assume that $\tau_1\dots\tau_t$ is a history of a computation $w_0\to\dots\to w_t$ of $\cal S$, where $\tau_1$ corresponds to a positive command $\theta$
of $S(M)$, $\tau_t$ corresponds to $\theta^{-1},$ and the other rules are not basic. Then the history $h'$ of the
computation $w_1\to\dots\to w_{t-1}$ is non-empty and $h' \equiv h_1\dots h_s,$ where each $h_i$ ($i=1,\dots,s$)
is  a maximal subword of $h'$ corresponding to the work of some machine $Z_{j(i)}(\theta_i,\pm).$ Since $h_1$ and $h_s$
correspond to $Z_1(\theta,+),$ either $s=1$ or there is $i$ such that both $h_{i-1}$
and $h_{i+1}$ correspond to the same $Z_{j(i)\pm 1}(\theta',\pm).$ It follows that
the computation with history $h_i$ satisfies the assumption of Lemma \ref{Z} (6), a contradiction.

(4) Property (4) follows from (2) and the definition of the computation $C_{S(M)}.$

(5) Given $C$, the computation $C_{\cal S}$ with the same space and with property
 $(C_{\cal S})_{S(M)} = C$ is briefly described above at the end of subsection \ref{compos}
 and with more details (though with submachines $Z_i(\theta, +)$ but without  $Z_i(\theta, -)$) in subsection 3.7 of \cite{OS}.
\endproof

 \begin{lemma} \label{positive} Let $C: \;\; w_0\to\dots\to w_t$ be a reduced computation of $\cal S$ such that the first command $\bar\theta_1$ and the last commands $\bar\theta_t$ are basic ones, and $C_{S(M)}\equiv W_0\to\dots\to W_s$ with $s\ge 2.$ Then the subcomputation $W_1\to\dots\to W_{s-1}$ is positive and $t-2 \ge 2^{|W_1|_a/l}.$ If $\bar\theta_1$ (if $\bar\theta_t$)
 is a start (is an accept) command, then the word $W_0$ (respectively, $W_t$) is also positive.
 \end{lemma}
 
 \proof To justify the first claim of the lemma, it suffices to prove that the word $W_1$ is positive under the assumption that 
there are no basic commands in the computation $w_1\to\dots\to w_{t-1},$ and the word $W_1$ corresponds
 to $w_1$. Therefore it suffices to prove that the word $w_1$ is positive.
 
 We first assume that the first command $\bar\theta$ of the history of $C$ is a positive basic command.
 Then $\bar\theta$ switches on the machine $Z_1(\theta, +).$ Since $s\ge 2,$ this machine must complete
 its work before the computation $C$ ends. The computation of $Z_1(\theta, +)$ cannot be empty since
 otherwise $\bar\theta$ would be followed by $\bar\theta^{-1}$ because $w_1$ involves the state letter $p_1(\theta,1^+).$ 
 But this would contradict the reducibility of $C.$ 
 
 Hence, by Lemma \ref{Z} (6), the work of $Z_1(\theta, +)$  sooner or later leads to an admissible word 
 $w_i$ containing the state letter $p_1(\theta, 3^+).$ By Lemma \ref{Z}(5), the subword of $w_1$ of the form 
 $k_1u_1 p_1(\theta, 1^+)k_2$ is positive, and the time of the work of $Z_1(\theta, +)$ is
 at least $2^{||u_1||}$. Then the machine $Z_1(\theta, +)$ finishes working and switches on
 the machine $Z_2(\theta, +),$ whose work similarly provides the positiveness of the subword of $w_1$ 
 having form $k_2u_2 p_2(\theta, 1+)k_3.$ Finally we obtain that $w_1$ is covered by positive
 subwords, and so it is positive itself.  The time of the work of all $Z_1(\theta, +),
 \dots, Z_l(\theta, +)$ is at least $2^{|W_1|_a/l}$ since $|W_1|_a=\sum_{j=1}^l ||u_j||.$
 
 If the command $\bar\theta^{-1}$ is positive, then it switches on the machine $Z_{l}(\theta,-),$ and we
 first obtain the positiveness of $k_{l}u_{l}p_{l}(\theta, 3^-)k_{l+1},$ then the positiveness of
 $k_{l-1}u_{l-1}p_{l-1}(\theta, 3^-)k_{l},$ and so on.

 Since start and  accept commands leave tape letters unchanged, the second statement follows
 from the positiveness of the word $W_1$ (of the word $W_{t-1}$). The lemma is proved. \endproof

\begin{lemma} \label{spacecalS}

(1) The S-machine $\cal S$ and the symmetric NTM $M$ recognize the same language $\cal L.$

 (2) The space complexities $S_{M}(n)$ and $S_{\cal S}(n)$ of $M$ and $\cal S,$ respectively, are equivalent.

(3) The generalized space complexities $S'_{M}(n)$ and $S'_{\cal S}(n)$ of $M$ and $\cal S$ are equivalent.

(4) We have $T'(n)\succeq \exp(S'_M(n))$ for the generalized time complexity $T'(n)$ of the machine $\cal S.$

\end{lemma}

\proof 

(1)  Assume that a word $u$ belongs to the language $\cal L$ recognized by $M$. By Lemma \ref{tmtostm},
this word belongs to the language of $S(M)$, and the accepted computation $C$ is positive. By
Lemma \ref{calS}(5), $u$ belongs to the language  ${\cal L}_{\cal S}$ of $\cal S.$

Now suppose $u$ belongs to ${\cal L}_{\cal S},$ and $C$ is the accepted computation. Then the computation
$C_{S(M)}$ is positive by Lemma \ref{positive}, and therefore this computation is also an accepted
computation of the machine $M$ by Lemma \ref{tmtostm}, and so $u\in \cal L.$

(2) The above argument shows that if a reduced computation $C:w_0\to\dots\to w_t$ of $\cal S$ accepts 
an admissible input word $w_0,$ then the computation $C_{S(M)}$ is a positive accepting computation 
of both $S(M)$ and $M,$ and so  $S_{M}(n) \le S_{\cal S}(n)$ by Lemma \ref{calS} (4). 
On the other hand, every accepted input configuration $W$ of $M$ is accepted by $S(M).$ By Lemma \ref{calS}(5), 
it has a copy accepted by $\cal S,$ and moreover, the accepting computations of $M$ and $\cal S$ need the same space. Therefore $S_{M}(n) \ge S_{\cal S}(n).$

(3) Assume that $C: W \equiv W_0\to\dots\to W_t$ is an accepting computation of $M$ such that $space_{M}(C)=space_M(W),$
and $C_{\cal S}: w\equiv w_0\to\dots\to w_s.$ For the given $w$ and $w_s$, we also consider a reduced computation 
$C':w\to\dots\to w_s$ of
$\cal S$ with minimal   space and the computation $(C')_{S(M)}: W'_0\to\dots\to W'_{t'}.$ 
If $t'=0$ (no basic rules), then $W_0\equiv W_{t}$ by Lemma \ref{Z} (2), and so  $space_{M}(C)= space_{\cal S}(C').$
Then we assume that $t'>0$ and note that $||W'_0||=||W_0|| $ and $||W'_{t'}||=||W_t||$ by Lemma \ref{Z} (2) since
the corresponding admissible words of $\cal S$ can be connected by computations without basic rules. 
Since $C'_{S(M)}$ is positive by Lemma \ref{positive}, it is also an accepted computation of the machine $M$
by Lemma \ref{tmtostm}.
Therefore, by Lemma \ref{calS} (4), $space_{M}(C)\le space_{S(M)}(C'_{S(M)})= space_{\cal S}(C').$ Since
$|w|_a=|W|_a,$ the last inequality proves that $S'_M(n)\le S'_{\cal S} (n)$ for every $n.$

Now we consider any accepting computation $C: w\equiv w_0\to\dots\to w_t$ of $\cal S$ with  $|w|\le n$, such that  
$space_{\cal S}(C)=space_{\cal S}(w).$ Without loss of generality, we may also assume that
$space_{\cal S}(C[m])=space_{\cal S}(w_m)$ for every subcomputation $C[m]$ of the form
$w_m\to\dots\to w_t.$ 
Let $\rho_{i_1}\equiv\bar\theta_1,\dots, \rho_t\equiv\rho_{i_s}\equiv\bar\theta_s$ be the basic
commands of the history $h\equiv\rho_1\dots \rho_t$ of $C,$  let $C'$ be the subcomputation of $C$
with history $h'\equiv\rho_1\dots \rho_{i_1-1},$ and let $C''$ have history $h''\equiv\rho_{i_1}\dots\rho_t;$
and so $h\equiv h'h''$ and $C=C'C''.$

We denote by $W_0\to W_1\to\dots\to W_s$ the computation $(C'')_{S(M)}$ and by $(C'')_{S(M)}[1]$ the
subcomputation $W_1\to\dots \to W_s,$ which is positive by Lemma \ref{positive}. 
Note that $space_{S(M)}( (C'')_{S(M)}[1]) = space_{S(M)}(W_1)$ since otherwise the subcomputation of $C[i_1]$
could be replaced by a subcomputation which needs less  space by Lemma \ref{calS} (5). Therefore by Lemma
\ref{calS} (4),
$$space_{\cal S}(C'')=space_{S(M)}((C'')_{S(M)})\le $$ $$space_{S(M)}( (C'')_{S(M)}[1])+|||W_1||-||W_0|||= space_{S(M)}(W_1)+c$$
since $|||W_1||-||W_0|||$ is bounded by a constant $c$ depending on the machine $M$ only.  

By Lemma \ref{calS} (2), we have $||w_j||\le ||w_{0}||$ for $j\le i_1$, and so $space_{\cal S}(C')\le |w_0|_a.$ Now, since $||W_1||\le ||W_0||+c\le ||w_{i_1-1}||+c\le ||w_0||+c$, we get  

$$ space_{\cal S}(C) \le \max(space_{\cal S}(C'), space_{\cal S}(C''))\le $$ 
$$\max(|w_0|_a, space_{S(M)}(W_1)+c)
\le\max (|w_0|_a, S'_{S(M)}(|w_0|_a+c)+c)$$

Hence $S'_{\cal S}(n)\le \max( S'_{S(M)}(n+c)+c, n).$ This inequality together with the inequality
$S'_M(n)\le S'_{\cal S} (n)$ obtained earlier, show that $S'_M(n)\sim S'_{\cal S} (n)$. 

(4) Let again $C: W\equiv W_0\to\dots\to W_t$ be an accepting computation of $M$ such that $space_{M}(C)=space_M(W)=S'_M(|W|_a),$ and $C_{\cal S}: w\equiv w_0\to\dots\to w_s.$
For the given $w$ and $w_s$, we also consider a reduced computation 
$C':w\to\dots\to w_{s'}$ of
$\cal S$ with minimal {\it time} $s'$ and the computation $(C')_{S(M)}: W'_0\to\dots\to W'_{t'},$
where one has $W'_0\equiv W_0$ and $W'_{t'}\equiv W_t.$ Therefore $(C')_{S(M)}$ is a positive computation by Lemma \ref{positive}, and  by the choice of $C$ and by Lemma \ref{tmtostm},  $space_{M}(C)\le space_M((C')_{S(M)}).$
Consider also a maximal subcomputation $C''$ of $C'$ starting and ending with basic commands.
Then $(C')_{S(M)}= (C'')_{S(M)}.$ By Lemma \ref{positive}, $space_{S(M)} (C'')_{S(M)}\le c\log (time_{\cal S}  (C''))$ for a constant $c>0.$ Thus by Lemma \ref{tmtostm}, $$S'_M(|W|_a)= space_M(C)\le space_{M}((C')_{S(M)})
= space_{S(M)}((C')_{S(M)})$$ $$= space_{S(M)}((C'')_{S(M)})
\le c\log (time_{\cal S}  (C''))\le c\log (time_{\cal S}  (C'))\le c\log(T'(|W|_a)),$$
and the lemma is proved.

\endproof

 \begin{lemma} \label{MtocalS}
(a) For every DTM $M$ recognizing a language $\cal L$ and having a  space complexity $S(n),$
(b) for every NTM $M$ recognizing a language $\cal L$ and having an FSC  space
complexity $f(n)$,  there exists an $S$-machine $\cal S$  with the following properties.
\begin{enumerate}

 \item The machine $\cal S$ recognizes the language $\cal L$.

\item Respectively, (a) both the   space and the generalized space complexities of $\cal S$ are equivalent to $S(n)$, 
(b) both the   space 
 and the generalized space complexities of $\cal S$ are equivalent
 to $f(n)^2$.

  \item Every command of $\cal S$ or its inverse inserts/deletes at
most one letter on the left and at most one letter on the right of
every state letter.

\item The machine $\cal S$ satisfies the ${\vec s}_{10}$-property. 

\item The unique start command  is of the form 
$q_1\to q'_1, q_2\tool q'_2,\dots, q_{k+1}\tool q'_{k+1},$ where 
$(q_1,\dots, q_{k+1}) \equiv {\vec s}_1.$

\item Any state letters $q$ from $Q_1$ is {\it passive}, i.e., there are
no commands of $\cal S$ of the form  $q\to q'u$ with non-empty $a$-word $u$.
\end{enumerate}
\end{lemma}

\proof For every DTM $M$ recognizing a language $\cal L$ and having a  space complexity $S(n)$ (for every NTM $M$ 
recognizing a language $\cal L$ and having an FSC  space
complexity $f(n)$) one can construct a symmetric NTM $M',$ as in the
formulation of Lemma \ref{generspace}. 
Then the properties 1 and 2 hold for 
the machine ${\cal S} = M'\circ Z$ by Lemma \ref{spacecalS}.

To provide the third property we use the same trick as for property
$4$ of Lemma \ref{generspace}.
 The properties 1 and 2
are obviously preserved.

To obtain the ${\vec s}_{10}$-condition of $\cal S$,  it suffices to add new state
letters and two new (positive) commands. Such a modification preserves
the other properties of $\cal S$, as this was noticed at the end of subsection 
\ref{definitions}. Now the ${\vec s}_{10}$-condition for $\cal S$ follows from
Lemma \ref{calS} (1). The form of the start command follows from our agreement
that all tapes numbered $2,\dots,k$ are empty for the input configurations of the machine $M.$
Finally, the left-most head of the Turing machine $M$ is passive being equal to the
separating symbol $\alpha_1;$ and the same property is inherited by 
${\cal S}$ as this follows from the definition  of the composition $M'\circ Z.$

\endproof

\section{Groups and diagrams }\label{gd}

Recall that given a finitely generated group $H$,
we want to construct an embedding of $H$ into a finitely presented group $G$
whose space function is equivalent to the deterministic space complexity of the algorithmic word problem in $H.$ The relations of $G$ are presented in this section.

\subsection{Construction of the embedding}\label{conemb}

Let $H$ be a finitely generated group with solvable word problem. To prove Theorem 
\ref{main1}, we will suppose that a Turing machine $M$ solves the word problem in $H.$ This 
implies that $H$ has a finite set of generators $\{a_1,..,a_m\}$, and a
word $w$ in these generators is accepted by $M$ iff $w=1$ in $H.$  Here we consider
only positive words in the generators since $M$ can work with positive words only,
 and so we assume that the set of generators is symmetric: for every $a_i$, there is
a generator $a_j$ such that $a_ia_j=1$ in $H.$ Thus every relation holding in $H$ follows from 
relations in $a_1,\dots,a_m,$ with positive left-hand side.

Further we will assume that one of the hypotheses (a), (b) of Lemma \ref{MtocalS} holds.
Therefore we also have the $S$-machine ${\cal S} =\langle X, Y, Q,
\Theta, \vec s_1, \vec s_0 \rangle $ provided by that lemma. The
input alphabet of ${\cal S}$ is the system of generators $a_1,\dots,a_m$ of the group $H$ together with
the symbols of the inverse letters $\{a_1^{-1},..,a_m^{-1}\}.$
Let $S'_{\cal S}(n)$ be the generalized space complexity of $\cal S.$

We denote by $\hat{\cal S}=\langle \hat X, \hat Y, \hat Q,
\hat\Theta, \hat\vec s_1, \hat\vec s_0 \rangle $ a copy of the $S$-machine 
${\cal S}.$ We will assume that $\hat X=X,$  $Y\cap \hat Y=X,$
$Q\cap \hat Q$ consists of the state letters of the vectors ${\vec s}_1,$ 
and ${\vec s}_0,$ 
and $\Theta\cap\hat\Theta=\emptyset.$ Therefore the machines $\cal S$ and $\hat{\cal S}$
have the same admissible input words.

The copy of a command $\theta$ of $\cal S$ is called
$\hat\theta.$ Similar notation is used for the $a$-letters from $\hat Y$
and $q$-letters from $\hat Q.$ The set of rules of the machine and admissible words of 
${\cal S}\cup \hat{\cal S} $ is, by definition, the union of the corresponding sets for 
$\cal S$ and for $\hat{\cal S}.$ The
following lemma is a clear consequence of these definitions.

\begin{lemma} \label{inH} The sets of  accepted admissible input words of the machines 
${\cal S},$ $\hat{\cal S}$ and ${\cal S}\cup\hat{\cal S}$ coincide, and so these machines 
recognize the same language $\cal L.$ They also have equal   space complexities and equal generalized 
space complexities. For every accepting computation $w_1\to\dots$ of ${\cal S}\cup\hat{\cal S},$ 
there is an accepting computation $w_1\to\dots$ of either ${\cal S}$ or $\hat{\cal S}$
whose length and  space does not exceed the length and  space of the original computation.
 \end{lemma}

We consider a group $G({\cal S},L)$  associated with the machine $\cal S.$
 Furthermore as in \cite{OS2}, we need a very similar group $\hat G({\cal S},L)$ to
produce a group embedding  required for the proof of Theorem  \ref{main1}.  

To define $G({\cal S},L)$ we need  many copies of
every letter used in the work of $\cal S;$ this enables us to apply a kind
of hyperbolic argument for the hub structure of van Kampen diagrams. Moreover, the
copies alternate with ``mirror copies"; this trick will be used in Section \ref{compare}.

Therefore we ``multiply" the machine $\cal S$ as follows. For some even $L\ge 40,$
we introduce
$L/2$ copies ${\cal S} ={\cal S}_1, {\cal S}_3,\dots, {\cal S}_{L-1}$ of the
machine ${\cal S}$ and $L/2$ mirror copies ${\cal S}_2, {\cal S}_4,\dots, {\cal S}_L$  of $\cal S.$
(A mirror copy of an arbitrary word $x_1\dots x_n$ 
is, by definition, $x_n\dots x_1$. For even $i,$ the rules  of ${\cal S}_i(L)$ 
transform the words in the mirror manner in comparison with $\cal S.$) 
We also add auxiliary separating state letters $k_1,\dots k_L.$  
For every admissible word $W$ of the machine $\cal S,$ we
define the words $W_1,W_2,\dots, W_L$, where $W_1,W_3,\dots
W_{L-1}$ are copies of $W$ in disjoint alphabets and $W_2,
W_4,..., W_L$ are mirror copies of $W$ also in disjoint alphabets.
The words of the form
$k_1W_1k_2W_2\dots k_LW_L$ are admissible words of the machine
${\cal S}(L)$. The rules  of ${\cal S}(L)$ are in one-to-one correspondence with the
rules of $\cal S;$ they transform the words $W_1, W_3,\dots$ as the rules of $\cal S$ do, 
they transform the words $W_2, W_4,\dots $ in the mirror manner, and 
do not change the new state letters
$k_1,\dots,k_L.$ We identify the set of rules of ${\cal S}(L)$ with $\Theta.$
Thus, by definition of ${\cal S}(L),$ we
also multiply the input and accept configurations. If $u$ is an
input word for $\cal S$, then we denote by $\Sigma(u,L)$ the
corresponding input configuration for ${\cal S}(L)$ containing $L/2$ copies of $u$
and $L/2$ mirror copies of $u$ as subwords. If the admissible words of $\cal S$ have $K$
state letters, then the admissible words of the machine ${\cal S}(L)$ 
have $N=(K+1)L$ state letters. 
Clearly the machine
${\cal S}(L)$ enjoys the properties (1)--(5) of machine $\cal S$ listed in Lemma \ref{MtocalS}. 

The finite set of generators of the group $G({\cal S},L)$ consists of {\em $q$-letters} corresponding to the states
of ${\cal S}(L)$, {\em $a$-letters} corresponding to the tape letters of ${\cal S}(L),$ and
{\em $\theta$-letters} corresponding to the commands.
Thus the set of generators consists of the set of (state) $q$-letters $Q({\cal S},L)=\sqcup_{i=1}^N Q_i$,  
the set of (tape) $a$-letters $Y=\sqcup_{i=1}^{N}
Y_i$ including the input alphabet $X=\sqcup_{i=1}^L X_i;$ and the $\theta$-letters from $N$ copies of
$\Theta^+,$ i.e., for every $\theta\in \Theta^+$, we have $N$
generators $\theta_1,\dots,\theta_N$.

The relations of the group $G({\cal S},L)$ correspond to the rules of the machine ${\cal S}(L)$;
for every $\theta=[U_1\to V_1,\dots U_N\to V_N]\in \Theta^+$, we have
\begin{equation}\label{rel1}
U_i\theta_{i+1}=\theta_i V_i,\,\,\,\, \qquad \theta_j a=a\theta_j, \,\,\,\, i,j=1,...,N
\end{equation}
for all $a\in \bar Y_j(\theta)$. (Here $\theta_{N+1}\equiv\theta_1.$ ) The first type of relations will be
called $(\theta,q)$-{\em relations}, the second type -
$(\theta,a)$-{\em relations}. 

The definition of the machine ${\hat S}(L)$ is similar to that of ${\cal S}(L)$
but the admissible words are of the form $ k_1 \hat W_1
k_2\hat W_2...k_L\hat W_L$, where every $\hat W_i$ is
obtained from $W_i$ after replacement of every letter $x$ by its
copy $\hat x$, and for $i=1,$ we, in addition, delete all
$a$-letters, i.e., the word $\hat W_1$ has no $a$-letters.
 (In other words, instead of the first copy of $\cal S,$ we use the ``machine''
with the same state letters but having no tape letters.) In particular, 
the word $\hat\Sigma(u,L)$ is obtained from $\Sigma(u,L)$ by omitting  the
 first occurrence of the input word $u.$ Again, it
is obvious that properties (1) - (5) from Lemma \ref{MtocalS} hold for the machine $\hat
{\cal S}(L)$ as well. 
The relations of the group $\hat G({\cal S},L)$ are

\begin{equation} \label{rel2}
\hat U_i\hat\theta_{i+1}=\hat\theta_i \hat V_i,\,\,\,\, i=1,...,N, \qquad \hat\theta_j \hat a=\hat a\hat\theta_j
\end{equation}
for all $\hat a\in Y_j(\hat\theta)$ and $j\in [K+1, N].$

We will also use the combined $S$-machine ${\cal} S(L)\cup\hat {\cal S}(L)$. Its
 admissible words are either the admissible words for ${\cal S}(L)$ or
 the admissible words for $\hat{\cal S}(L)$, the set of rules is the
 union of the rules for ${\cal S}(L)$ and $\hat {\cal S}(L).$ Note that 
 the machine ${\cal} S(L)\cup\hat {\cal S}(L)$ does not satisfy the
 $\vec s_{10}$-condition.

 The subwords of the form $(k_iW_ik_{i+1})^{\pm 1}$ (indices taken
 modulo $L$) of the admissible words of the machine ${\cal S}(L)$ are
 said to be  $i$-sector words. Similarly one define $i$-sector
 words for the machines $\hat{\cal S}(L)$ and for ${\cal S}(L)\cup\hat{\cal S}(L)$.
 (Recall that the $q$-letters of the $1$-sector words of ${\cal S}(L)$, except for $k_1$ and $k_2$, are  identified with the corresponding letters of the original $S$-machine $\cal S$.) 
 The state letters of the $i$-sector of ${\cal S}(L)$ are the letters from $\sqcup_{j=(i-1)(K+1)+1}^{i(K+1)} Q_j ,$
 and the tape letters of the $i$-sector are the letters from $\sqcup_{j=(i-1)(K+1)+1}^{i(K+1)} Y_j .$
 Similarly we have state and tape letters of the $i$-sector for $\hat{\cal S}(L)$ and for
 ${\cal S}(L)\cup\hat{\cal S}(L)$.
 By definition, the $\theta$-letters  $\theta_j$ and $\hat\theta_j$ with subscripts
 $j\in [(i-1)(K+1)+1, i(K+1)]$ are $\theta$-letters of the $i$-sector of ${\cal S}(L)\cup\hat{\cal S}(L)$. 
The state, tape and  the $\theta$-letters of the $i$-sector constitute the
{\it alphabet} ${\cal A}_i$ of the $i$-sector.

 The group $G({\cal S}\cup \hat{\cal S}, L)$ is given by the generators and
 relations of both groups $G({\cal S},L)$ and $\hat G({\cal S},L)$. 

Finally, the required group $G$ is given by the generators and
relations of the group $G({\cal S}\cup\hat{\cal S},L)$ and one more additional
relation, namely the {\it hub}-relation  
\begin{equation}\label{rel3}
\Sigma_0=1,
\end{equation}
where
$\Sigma_0=\Sigma(L)$ is the accept word (of length $N$) of the machine ${\cal S}(L)$ (and of $\hat{\cal S}(L)$ as well). 

Suppose an admissible word $W'$ of ${\cal S} (L)$ is
obtained from an admissible word $W$ by an application of a rule
$\theta: \; [U_1\to V_1,\dots,U_N\to V_N]$. This implies that $W = U_1w_1U_2w_2\dots U_N w_N,$
where  $w_i$ is a word in the alphabet $Y_i(\theta)$, and therefore 
$W'=\theta_1^{-1}W\theta_1$ in $G({\cal S},L)$ by Relations (\ref{rel1}), since
$\theta_1\equiv\theta_{N+1}.$ 

Now suppose $u$ is a positive word in the alphabet $\{a_1,\dots, a_m\}$
vanishing in the group $H.$ Then it is recognized by the machine $\cal S,$  
and so the word $\Sigma(u)$ is accepted by the machine ${\cal S}(L)$ and therefore it is conjugate to $\Sigma_0$
in the group $G({\cal S}\cup\hat{\cal S},L).$ Consequently, we have  that
$\Sigma(u)=1$ in $G$  by (\ref{rel3}). Similarly, $\hat\Sigma(u)=1$ in $G.$

Recall that we have identified the alphabet of $1$-sector words of ${\cal S}(L)$ with
the alphabet of $\cal S.$ Therefore the word $\hat\Sigma(u)$ 
results from $\Sigma(u)$ after deleting the subword $u$ in the alphabet
of generators of $H$. Hence the
relations $\Sigma(u)=\hat\Sigma(u)=1$ imply $u=1$ in $G.$ Since
the language of accepted words for ${\cal S}(L)$ contains all the
defining relations of $H,$ we have obtained

\begin{lemma} \label{homo} The mapping  $a_i \mapsto a_i$ ($i=1,\dots, m$) extends to a
homomorphism of the group $H$ to $G.$
\end{lemma}

In Section \ref{compare} we show that this homomorphism is injective.

\subsection{Minimal diagrams}\label{md}

As in \cite{OS1}, we enlarge the set of defining relations of the
group $G$ by adding some consequences of defining relations.
Taking into account Lemma \ref{homo}, we include all cyclically reduced
relations of the group $H$ generated by the set
$\{a_1,\dots,a_m\},$ i.e., all non-empty cyclically reduced words
in $\{a_1^{\pm 1},\dots,a_m^{\pm 1}\}$ which are equal to $1$ in the group $H.$
These relations will be called $H$-{\it relations}. 

We denote by $G_1$ the group given by all generators of the group $G$, by
all $H$-relations, and by all defining relations of $G$ except for
the hub-relation (\ref{rel3}).

\medskip

Recall that a van Kampen {\it diagram} $\Delta $ over a presentation
$P=\langle B\; | \; \mathcal R\rangle$ (or just over the group $P$)
is a finite oriented connected and simply--connected planar 2--complex endowed with a
labeling function $\Lab : E(\Delta )\to B^{\pm 1}$, where $E(\Delta
) $ denotes the set of oriented edges of $\Delta $, such that $\Lab
(e^{-1})\equiv \Lab (e)^{-1}$. Given a cell (that is a 2-cell) $\Pi $ of $\Delta $,
we denote by $\partial \Pi$ the boundary of $\Pi $; similarly,
$\partial \Delta $ denotes the boundary of $\Delta $. The labels of
$\partial \Pi $ and $\partial \Delta $ are defined up to cyclic
permutations. An additional requirement is that the label of any
cell $\Pi $ of $\Delta $ is equal to (a cyclic permutation of) a
word $R^{\pm 1}$, where $R\in \mathcal R$. Labels and lengths of
paths are defined as for Cayley graphs.

The van Kampen Lemma states that a word $W$ over the alphabet $B^{\pm 1}$
represents the identity in the group $P$ if and only
if there exists a diagram $\Delta
$ over $P$ such that 
$\Lab (\partial \Delta )\equiv W$
(\cite{LS}, Ch. 5, Theorem 1.1).

We will study diagrams over the groups $G$ and $G_1$. The edges labeled by state
letters ( = $q$-{\it letters}) will be called $q$-{\it edges}, the edges labeled by tape
letters (= $a$-{\it letters}) will be called $a$-{\it edges}, and the edges labeled by 
the letters from $\Theta$ and $\hat\Theta$ (= $\theta$-{\it letters}) are $\theta$-{\it edges}.
 The cells corresponding
to Relation (\ref{rel3}) are called {\it hubs}, the cells corresponding
to Relations (\ref{rel1}) and (\ref{rel2}) are called $(\theta,q)$-{\it cells}
if they involve $q$-letters, and they are called $(\theta,a)$-{\it cells} otherwise.
The cells corresponding to arbitrary relations of $H$ are $H$-{\it cells.}

The resulting presentation
and the diagrams over it are {\it graded} by the ranks of defining
words and cells as follows. The hubs are the cells of the highest
rank, the rank of $(\theta,q)$-cells is higher than the rank of
$H$-cells (and in Lemma \ref{simple}, $(\theta,k_i)$-cells, with $i\ne 1,2$ are higher than
other $(\theta,q)$-cells), and the $(\theta,a)$-cells are of the lowest rank.

If $\Delta$ and $\Delta'$ are diagrams over $G,$ then we say that
$\Delta$ has a higher type than $\Delta'$ if $\Delta$ has more
hubs, or if the numbers of hubs are the same, but $\Delta$ has more
cells which are next in the hierarchy, and so on.

Clearly the defined partial order on the set of diagrams satisfies
the descending chain condition, and so there is a diagram having
the smallest type among all diagrams with the same boundary label.
Such a diagram is called {\it minimal}.

\subsection{Bands and trapezia}\label{bands}
\label{bands} From now on, we shall mainly consider minimal van
Kampen diagrams. In particular the diagrams are {\it reduced}, i.e.,
they contain no pair of cells that have a common edge $e$ such that the boundary
paths of these cells starting with $e$ have the same label. 
As in \cite{SBR}, \cite{BORS}, to explore van Kampen diagrams over the groups
$G$ and $G_1$ we shall use their simpler subdiagrams such as bands and trapezia.

Here we repeat some necessary definitions.

\begin{df}Let $\cal Z$ be a subset of the set of generators ${\cal X}$ of the group $G$. A
$\cal Z$-band $\bb$ is a sequence of cells $\pi_1,...,\pi_n$ in a \vk
diagram such that

\begin{itemize}
\item Every two consecutive cells $\pi_i$ and $\pi_{i+1}$ in this
sequence have a common edge $e_i$ labeled by a letter from $\cal Z$.
\item Each cell $\pi_i$, $i=1,...,n$ has exactly two $\cal Z$-edges,
$e_{i-1}$ and $e_i$ (i.e. edges labeled by a letter from $\cal Z$).

\item If $n=0$, then $\bb$ is just a $\cal Z$-edge.
\end{itemize}
\end{df}

The counterclockwise boundary of the subdiagram formed by the
cells $\pi_1,...,\pi_n$ of $\bb$ has the factorization $e\iv q_1f q_2\iv$
where $e=e_0$ is a $\cal Z$-edge of $\pi_1$ and $f=e_n$ is an $\cal Z$-edge of
$\pi_n$. We call $q_1$ the {\em bottom} of $\bb$ and $q_2$ the
{\em top} of $\bb$, denoted $\bott(\bb)$ and $\topp(\bb)$.
Top/bottom paths and their inverses are also called the {\em
sides} of the band. The $\cal Z$-edges $e$
and $f$ are called the {\em start} and {\em end} edges of the
band. If $n\ge 1$ but $e=f,$ then the $\cal Z$-band is called a $\cal Z$-{\it annulus}.

We say that a ${\cal Z}_1$-band and a ${\cal Z}_2$-band {\em cross} if 
they have a common cell and ${\cal Z}_1\cap {\cal Z}_2=\emptyset.$

We shall call a $\cal Z$-band {\em maximal} if it is not contained in
any other $\cal Z$-band.

We will consider $q$-bands where $\cal Z$ is one of the sets $Q_i$ of state letters
for the machine ${\cal S}(L)\cup\hat{\cal S}(L)$,
$\theta$-bands for every $\theta\in\Theta$, and $a$-bands where
$M=\{a\}\subseteq Y$. 
 
The convention is that $a$-bands do not
contain $q$-cells, and so they consist of $(\theta,a)$-cells  only.

The papers \cite{O}, \cite{BORS}, \cite{OS2} contain the proof of the
following lemma in a more general setting. (In contrast to Lemmas 6.1 \cite{O} and
 3.11 \cite{OS2}, we have no $x$-cells here.)

\begin{lemma}\label{NoAnnul}
A minimal van Kampen diagram $\Delta$ over $G_1$ has no
$q$-annuli, no $\theta$-annuli, and no  $a$-annuli.
Every $\theta$-band of $\Delta$ shares at most one cell with any
$q$-band and with any $a$-band. 

\end{lemma}

If $W\equiv x_1...x_n$ is a word in an alphabet $X$, $Y$ is another
alphabet, and $\phi\colon X\to Y\cup\{1\}$ (where $1$ is the empty
word) is a map, then $\phi(W)\equiv\phi(x_1)...\phi(x_n)$ is called the
{\em projection} of $W$ onto $Y$. We shall consider the
projections of words in the generators of $G$ onto
$\Theta\sqcup\hat\Theta$ (all $\theta$-letters map to the
corresponding element of $\Theta\sqcup\hat\Theta$,
all other letters map to $1$), and the projection onto the
alphabet $\{Q_1\sqcup \dots \sqcup Q_N\}$ (every
$q$-letter maps to the corresponding $Q_i$, all other
letters map to $1$).

\begin{df}\label{dfsides}{\rm  The projection of the label
of a side of a $q$-band onto the alphabet $\Theta^{\pm 1}$ is
called the {\em history} of the band. The projection of the label
of a side of a $\theta$-band onto the alphabet $\{Q_1,...,Q_n\}$
is called the {\em base} of the band. Similarly we  define the
history of a word and the base of a word. The base of a word $W$
is denoted by $base(W)$}. It will be convenient  to use representatives of
$Q_1, ...,Q_N$ in base words. For example, if $k\in Q_1$, $q\in Q_2$, we shall
say that the word $kaq$ has base $kq$ instead of
$Q_1Q_2$.
\end{df}

\begin{df}\label{dftrap}
{\rm Let $\Delta$ be a minimal \vk diagram over $G_1$
which has  boundary path of the form $p_1\iv q_1p_2q_2\iv$ where:
\medskip

$(TR_1)$ $p_1$ and $p_2$ are sides of $q$-bands,

\medskip

$(TR_2)$ $q_1$, $q_2$ are maximal parts of the sides of
$\theta$-bands such that $\Lab(q_1)$, $\Lab(q_2)$ start and end
with $q$-letters,

\medskip

$(TR_3)$ for every $\theta$-band $\ttt$ in $\Delta$, the labels of
$\topp(\ttt)$ and $\bott(\ttt)$ are reduced.

\medskip

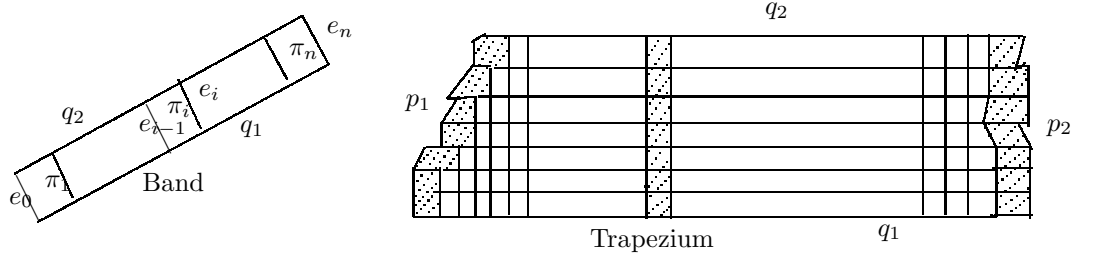
\begin{figure}[h!]
\unitlength 1mm 
\linethickness{0.4pt}
\ifx\plotpoint\undefined\newsavebox{\plotpoint}\fi 
\begin{picture}(148.25,40)(0,115)
\put(76.5,148){\line(1,0){64.25}}
\put(74.75,143.75){\line(1,0){66.25}}
\put(72.75,140){\line(1,0){68.25}}
\put(71,136.5){\line(1,0){69.75}}
\put(69.5,133.25){\line(1,0){71.75}}
\put(68,130.25){\line(1,0){73.25}}
\put(66.75,127.25){\line(1,0){74.5}}
\put(65.5,124){\line(1,0){75.75}}
\put(76.75,147.75){\line(0,-1){23.5}}
\put(74.25,143.75){\line(0,-1){19.75}}
\put(72.25,140){\line(0,-1){16}}
\put(70,133.5){\line(0,-1){9.5}}
\put(67.5,130.5){\line(0,-1){6.5}}
\put(73.5,148){\line(1,0){3.25}}
\multiput(73.5,148.25)(-.0583333,-.0333333){30}{\line(-1,0){.0583333}}
\put(72,147.25){\line(0,-1){2.75}}
\put(72,144.5){\line(0,1){.25}}
\put(72,144.75){\line(0,-1){.75}}
\put(72,144){\line(1,0){2.25}}
\multiput(72,144.25)(-.033653846,-.043269231){104}{\line(0,-1){.043269231}}
\put(68.5,139.75){\line(1,0){1.75}}
\multiput(70,140)(-.03358209,-.05223881){67}{\line(0,-1){.05223881}}
\put(67.75,136.5){\line(1,0){3.25}}
\multiput(69.5,139.75)(.3125,.03125){8}{\line(1,0){.3125}}
\put(67.75,136.5){\line(0,-1){3.25}}
\put(67.75,133.25){\line(1,0){2}}
\put(65.75,133.25){\line(1,0){2}}
\multiput(65.5,133.25)(-.03333333,-.06666667){45}{\line(0,-1){.06666667}}
\put(64,130.25){\line(1,0){3.25}}
\put(64,130.5){\line(0,-1){6.5}}
\put(64,124){\line(0,1){0}}
\put(64,124){\line(1,0){1.75}}
\put(79.25,147.75){\line(0,-1){23.75}}
\multiput(73.68,147.68)(-.375,-.5){5}{{\rule{.4pt}{.4pt}}}
\multiput(74.93,147.93)(-.45,-.75){11}{{\rule{.4pt}{.4pt}}}
\multiput(76.18,147.68)(-.5,-.77941){18}{{\rule{.4pt}{.4pt}}}
\multiput(76.43,145.93)(-.5,-.5833){4}{{\rule{.4pt}{.4pt}}}
\multiput(73.68,141.43)(-.54167,-.75){19}{{\rule{.4pt}{.4pt}}}
\multiput(66.68,133.18)(-.45,-.5){6}{{\rule{.4pt}{.4pt}}}
\multiput(71.93,137.18)(-.55357,-.75){15}{{\rule{.4pt}{.4pt}}}
\multiput(71.93,135.68)(-.375,-.5625){5}{{\rule{.4pt}{.4pt}}}
\multiput(69.68,132.18)(-.4167,-.5833){4}{{\rule{.4pt}{.4pt}}}
\multiput(67.43,128.93)(-.54167,-.58333){7}{{\rule{.4pt}{.4pt}}}
\multiput(66.93,127.18)(-.5,-.55){6}{{\rule{.4pt}{.4pt}}}
\multiput(67.18,125.68)(-.5833,-.5833){4}{{\rule{.4pt}{.4pt}}}
\put(95,148){\line(0,-1){24.25}}
\put(98.25,148){\line(0,-1){24.25}}
\multiput(96.93,147.93)(-.4375,-.625){5}{{\rule{.4pt}{.4pt}}}
\multiput(97.93,146.93)(-.45833,-.66667){7}{{\rule{.4pt}{.4pt}}}
\multiput(97.68,144.18)(-.5,-.7){6}{{\rule{.4pt}{.4pt}}}
\multiput(97.93,142.43)(-.42857,-.57143){8}{{\rule{.4pt}{.4pt}}}
\multiput(97.68,140.43)(-.5,-.7){6}{{\rule{.4pt}{.4pt}}}
\multiput(97.93,138.43)(-.55,-.75){6}{{\rule{.4pt}{.4pt}}}
\multiput(97.68,136.18)(-.5,-.7){6}{{\rule{.4pt}{.4pt}}}
\multiput(97.68,133.93)(-.5,-.75){6}{{\rule{.4pt}{.4pt}}}
\multiput(97.68,131.18)(-.5,-.65){6}{{\rule{.4pt}{.4pt}}}
\multiput(95.18,127.93)(0,0){3}{{\rule{.4pt}{.4pt}}}
\multiput(95.18,127.93)(-.25,0){3}{{\rule{.4pt}{.4pt}}}
\multiput(97.93,129.18)(-.5,-.5){7}{{\rule{.4pt}{.4pt}}}
\multiput(98.18,126.93)(-.6,-.5){6}{{\rule{.4pt}{.4pt}}}
\multiput(97.93,125.18)(-.3333,-.3333){4}{{\rule{.4pt}{.4pt}}}
\put(140.5,148.25){\line(0,-1){8.5}}
\put(140.5,139.75){\line(0,1){0}}
\multiput(140.5,140)(-.0326087,-.1521739){23}{\line(0,-1){.1521739}}
\put(139.75,136.5){\line(0,1){.25}}
\multiput(139.75,136.75)(.03365385,-.06730769){52}{\line(0,-1){.06730769}}
\put(141.5,133.25){\line(0,1){0}}
\put(141.5,133.25){\line(0,-1){9.25}}
\put(141.5,124){\line(-1,0){.25}}
\put(140.25,148){\line(1,0){4.75}}
\multiput(145.25,148)(-.03289474,-.11842105){38}{\line(0,-1){.11842105}}
\put(144.25,144){\line(1,0){1.5}}
\put(145.75,144){\line(0,-1){8}}
\put(145.75,136.5){\line(-1,0){1.5}}
\multiput(144.5,136.25)(.03365385,-.0625){52}{\line(0,-1){.0625}}
\put(146,133.75){\line(0,-1){10}}
\put(146,124.25){\line(1,0){.25}}
\multiput(140.75,124.25)(.0625,-.03125){8}{\line(1,0){.0625}}
\put(141.25,124.25){\line(1,0){4.75}}
\multiput(142.93,147.93)(-.45,-.5){6}{{\rule{.4pt}{.4pt}}}
\multiput(144.43,147.68)(-.5,-.60714){8}{{\rule{.4pt}{.4pt}}}
\multiput(143.68,144.43)(-.5,-.625){7}{{\rule{.4pt}{.4pt}}}
\multiput(145.18,142.93)(.25,.5){3}{{\rule{.4pt}{.4pt}}}
\multiput(144.93,143.93)(.125,.25){3}{{\rule{.4pt}{.4pt}}}
\multiput(145.18,144.43)(-.5,-.65){11}{{\rule{.4pt}{.4pt}}}
\multiput(145.43,142.18)(-.55,-.65){11}{{\rule{.4pt}{.4pt}}}
\multiput(145.18,139.93)(-.5625,-.625){9}{{\rule{.4pt}{.4pt}}}
\multiput(140.68,134.93)(-.125,-.125){3}{{\rule{.4pt}{.4pt}}}
\multiput(145.43,137.43)(-.5,-.53125){9}{{\rule{.4pt}{.4pt}}}
\multiput(144.93,134.68)(-.58333,-.58333){7}{{\rule{.4pt}{.4pt}}}
\multiput(145.68,133.18)(-.60714,-.5){8}{{\rule{.4pt}{.4pt}}}
\multiput(145.68,131.68)(-.57143,-.46429){8}{{\rule{.4pt}{.4pt}}}
\multiput(145.68,129.93)(-.60714,-.53571){8}{{\rule{.4pt}{.4pt}}}
\multiput(145.68,127.68)(-.75,-.5){6}{{\rule{.4pt}{.4pt}}}
\multiput(145.68,126.43)(-.5625,-.4375){5}{{\rule{.4pt}{.4pt}}}
\put(140.25,143.75){\line(1,0){4.5}}
\put(140.25,140){\line(1,0){5.5}}
\put(140,136.5){\line(1,0){4.75}}
\put(141.5,133.25){\line(1,0){4.5}}
\put(141,130.5){\line(1,0){5.25}}
\put(141,127.25){\line(1,0){4.75}}
\put(137.75,148.25){\line(0,-1){24}}
\put(134.75,148){\line(0,-1){23.75}}
\put(131.75,147.75){\line(0,-1){23.5}}
\put(63,138.75){$p_1$}
\put(125.75,121.5){$q_1$}
\put(148.25,135.5){$p_2$}
\put(110.75,151){$q_2$}
\put(87.5,120){Trapezium}
\multiput(11,129.75)(.0613964687,.0337078652){623}{\line(1,0){.0613964687}}
\multiput(49.25,150.75)(.033653846,-.0625){104}{\line(0,-1){.0625}}
\multiput(52.75,144.25)(-.0621990369,-.0337078652){623}{\line(-1,0){.0621990369}}
\put(11.25,129.5){\line(1,-2){3}}
\multiput(16,132.5)(.03353659,-.07317073){82}{\line(0,-1){.07317073}}
\multiput(44.25,147.75)(.03370787,-.06179775){89}{\line(0,-1){.06179775}}
\put(28.5,139){\line(1,-2){3.25}}
\put(31.75,132.5){\line(0,1){0}}
\multiput(32.75,142.25)(.03370787,-.07303371){89}{\line(0,-1){.07303371}}
\put(15,128){$\pi_1$}
\put(47.5,145.5){$\pi_n$}
\put(31.25,137.75){$\pi_i$}
\put(10.25,125.75){$e_0$}
\put(27.5,135.5){$e_{i-1}$}
\put(35.5,140){$e_i$}
\put(41,135.25){$q_1$}
\put(17.25,137){$q_2$}
\put(28,127.5){Band}
\put(52.5,148.5){$e_n$}
\end{picture}
\end{figure}

Then $\Delta$ is called a {\em trapezium}. The path $q_1$ is
called the {\em bottom}, the path $q_2$ is called the {\em top} of
the trapezium, the paths $p_1$ and $p_2$ are called the {\em left
and right sides} of the trapezium. The history of the $q$-band
whose side is $p_2$ is called the {\em history} of the trapezium;
the length of the history is called the {\em height}  of the
trapezium. The base of $q_1$ is called the {\em base} of the
trapezium.}
\end{df}

\begin{rk}{\rm (1) Property $(TR_3)$ is easy to achieve: by folding edges
with the same labels having the same initial vertex, one can make
the boundary label of a subdiagram in a \vk diagram reduced, see
\cite{SBR}.

(2) Notice that the top (bottom) side of a
$\theta$-band $\ttt$ does not necessarily coincide with the top
(bottom) side $q_2$ (side $q_1$) of the corresponding trapezium of height $1$, and $q_2$
($q_1$) is
obtained from $\topp(\ttt)$ (resp. $\bott(\ttt)$) by trimming the
first and the last $a$-edges if these paths start and/or end with $a$-edges.} 
We shall denote the
trimmed top and bottom sides of $\ttt$ by $\ttopp(\ttt)$ and
$\tbott(\ttt)$. By definition, for arbitrary $\theta$-band $\cal T,$ $\ttopp(\cal T)$
is obtained by such a trimming only if $\cal T$ starts and/or ends with a
$(\theta,q)$-cell; otherwise $\ttopp(\cal T)\equiv\topp(\cal T).$
The definition of $\tbott(\cal T)$ is similar. 
\end{rk}

The trapezium $\Delta$ is said to be an $i$-sector if the labels of its
top and bottom paths are $i$-sector words.

\begin{lemma} \label{calAi} Let $\Gamma$ be an $i$-sector, where $i\ne 1$. Then
the sides $p_1$ and $p_2$ are the sides of maximal $k_i$- and $k_{i+1}$-bands
${\cal K}_i$ and ${\cal K}_{i+1} $ of $\Gamma,$ respectively.
If an edge $e$ of $\Gamma,$ belongs to neither ${\cal K}_i$ nor ${\cal K}_{i+1},$ then $\Lab(e)\in {\cal A}_i.$ In particular, $\Gamma$ has no $H$-cells.
\end{lemma}

\proof The first assertion follows from Lemma \ref{NoAnnul}, since the labels
of the top and bottom of an $i$-sector are of the form $k_i\dots k_{i+1}.$ 
Then, by the same lemma, the second assertion is true for the edges of all
$\theta$-cells since every maximal $\theta$-band must connect ${\cal K}_i$ and ${\cal K}_{i+1}.$
Since $i\ne 1$ and the labels of the boundary edges of $H$-cells belong to ${\cal A}_1,$
the $H$-cells of $\Gamma$  have no edges in common either with the $\theta$-cells of
$\Gamma$ or with the boundary $\partial\Gamma.$ Now the minimality of the diagram $\Gamma$
implies that $\Gamma$ has no $H$-cells at all, and so the lemma is proved. \endproof

The following lemma claims that every $i$-sector, $i\ne 1,$ simulates the work
of ${\cal S}(L)\cup\hat{\cal S}(L).$ It summarizes the assertions  of Lemmas 6.1,
6.3, 6.9, and 6.16 from \cite{OS2}. For the formulation (1) below, it is important
that $\cal S$ (and ${\cal S} \cup\hat{\cal S}$) is an $S$-machine. The analog of 
this statement is false for Turing machines. (See \cite{OS1} for a discussion.) 

\begin{lemma}\label{simul} (1) Let $\Delta$ be an $i$-sector for some
$i\ne 1$ with history $\theta_1\dots\theta_d$. Assume that $\Delta$ 
has consecutive maximal $\theta$-bands  ${\cal T}_1,\dots
{\cal T}_d$, and $k_iW_jk_{i+1}$ and $k_iW'_jk_{i+1}$ are the bottom and the 
top labels of ${\cal T}_j,$ ($j=1,\dots,d$).
Let  $U_j$ (resp. $V_j,$ $
i=1,\dots d$) be the copies of $W_j$ (resp. $W'_j$) in the alphabet
of the machine ${\cal S} \cup\hat{\cal S}$. Then $U_j$, $V_j$ are admissible
words for ${\cal S} \cup\hat{\cal S},$ and
$$V_1\equiv U_1\cdot \theta_1, U_2\equiv V_1,\dots, U_d \equiv V_{d-1}, V_d\equiv U_d\cdot \theta_d$$

(2) For every reduced computation $U\cdot h \equiv V$ of $\cal S$ (of $\hat{\cal S}$) 
with $|h|\ge 1$
and for every $i\in [1,N]$ (for every $i\in [2,N]$),
there exists an $i$-sector $\Delta$ with history $h$ and without $H$-cells, whose bottom and top labels
are $k_iU'k_{i+1}$ and $k_iV'k_{i+1}$, where $U'$ (resp. $V'$) is the copy or
the mirror copy of the word $U$ (of $V$) in the alphabet ${\cal A}_i.$ These copies
are mirror copies iff $i$ is even.
\end{lemma}

We call an $i$-sector {\it accepted}
if its top label $k_iV'k_{i+1}$ is the $i$-sector subword of the word $\Sigma_0.$
For $i\ne 1,$  the computation $V_1\to\dots\to V_d$ of the machine $\cal S\cup\hat{\cal S}$ provided
by Lemma \ref{simul} (1) for an accepted $i$-sector is accepting.

\subsection{Replicas} \label{replicas}

Let $\Delta$ be an $i$-sector, where $i\ne 1$. Then by Lemma \ref{calAi}, for every $i'\ne 1$,
one can relabel the edges of $\Delta$ (or of the mirror copy of $\Delta$ if $i-i'$ is odd)
and obtain an $i'$-sector $\Delta'$, which is just
a copy (or a mirror copy) of $\Delta$. 
This symmetry and mirror symmetry are necessary for the diagram surgery
we will utilize in the proof of Lemma \ref{cut} of Section 4.
But one cannot construct such a copy of the $i$-sector if $i'=1$ since there are no
commutativity relations $\hat\theta_j \hat a=\hat a\hat\theta_j$  if
 $\hat\theta_j\in {\cal A}_1$  (see (\ref{rel2})). However one can construct an ersatz-copy of $\Delta$ called a replica, if $\Delta$ is an accepted sector.
 This construction will be used in Sections 4 and 5.
 
 For $i'=1,$
 we construct the {\it replica} of $\Delta$ as follows. (We assume below that $i$ is odd, otherwise
 one first replaces $\Delta$ by its mirror copy.)
 
 At first, the relations (\ref{rel1}) and (\ref{rel2}) make it possible to 
 replace the maximal $k_i$-band ${\cal K}_i$ and $k_{i+1}$-band ${\cal K}_{i+1}$ of
 $\Delta$ by their copies ${\cal K}_1$ and ${\cal K}_2$, respectively. Similarly,
 we replace every maximal $\theta$-band of $\Delta$ by its copy if $\theta$ is
 a command of the machine $\cal S.$ Now let $\cal T$ be a maximal $\theta$-band
 of $\Delta$ and $\theta$ a command of the machine $\hat{\cal S}.$ This $\cal T$
 consists of $(\theta,q)$- and $(\theta,a)$-cells. To construct the replica $\cal T'$
 of $\cal T$ we take the 'copies' of $(\theta,q)$-cells only (but no $a$-edges in these
 'copies', the $a$-edges are contracted to vertices) and build $\cal T'$ from
 them by identifying $\theta$-edges of neighboring cells in the order in which the original $(\theta,q)$-cells appear in $\cal T.$
 
 It remains to close up the holes between $\theta$-bands ${\cal T'}_s$ and ${\cal T'}_{s+1}$ for
 consecutive ${\cal T}_s$ and ${\cal T}_{s+1}.$ If the corresponding letters $\theta_s$ and
 $\theta_{s+1}$ are both commands of $\cal S,$ then we just identify the top of the copy
 ${\cal T'}_s$ and the bottom of the copy ${\cal T'}_{s+1}$. 
 This is possible since $\phi(\topp{\cal T}_s)\equiv
 \phi(\bott{\cal T}_{s+1})$ by Lemma \ref{simul} (1). A similar identification works if $\theta_s$ and $\theta_{s+1}$
 are both the commands of $\hat{\cal S}.$
 
 Assume now that $\theta_s$ is a command of $\cal S$ and $\theta_{s+1}$ is a command of $\hat{\cal S}$
 (or vice versa).  Then $\Lab(\topp({\cal T}_s))\equiv\Lab(\bott{\cal T}_{s+1}))\equiv k_iW_sk_{i+1},$ and the copy $W_s'$ of
 $W_s$ in the alphabet ${\cal A}_1$ is an admissible word by both machines $\cal S$ and $\hat{\cal S}$ since both commands $\theta_s^{-1}$ and $\theta_{s+1}$ are applicable to it. But the only common tape letters  of these machines are the letters of the input alphabet of $\cal S,$ and the only common state letters  of these machines are the letters from the start vector ${\vec s}_1$ and the accept vector ${\vec s}_0$
 of $\cal S$.   Then by property (5) of Lemma \ref{MtocalS}, either $\theta_s^{-1}$ is  the start rule or $\theta_s$ is the accept  rule of $\cal S.$ 
 
In the latter case, we have $\Lab(\topp({\cal T'}_s))\equiv\Lab(\bott{\cal T'}_{s+1}))\equiv k_iW_sk_{i+1}$ 
since the accept words have no $a$-letters at all, and so  identification of 
the top of the copy
 ${\cal T'}_s$ and the bottom of the copy ${\cal T'}_{s+1}$ is possible again.
Let us consider the former case. Then the rule $\theta_s^{-1}$  is of the form 
$q_1\to q'_1, q_2\tool q'_2,\dots, q_{K}\tool q'_{K},$ where 
$(q_1,\dots, q_{K}) \equiv  {\vec s}_1.$
Therefore $W_s'$ can have tape letters only between $q_1$ and $q_2,$ i.e., $W'_s\equiv q_1uq_2q_3\dots q_K,$
where $u$ is a word in the input alphabet. 

Since the sector $\Delta$ is accepted and the rules of $S$-machines are invertible, we have
that $W'_s$ is accepted by the machine ${\cal S}(L)\cup \hat{\cal S}(L),$ and so it is accepted
by $\cal S$ by Lemma \ref{inH}. Hence the word $u$ belongs to the language recognized by the machine
$\cal S$, and therefore $u$ is a word in the generators of $H$, and $u=1$ in $H$ by the 
definition of that language.

Now, to close up the hole between the bands ${\cal T'}_s$ and ${\cal T'}_{s+1}$, it suffices to
paste in an $H$-cell $\pi$ labeled by the cyclically reduced form of the word $u$ between them because 
$\Lab(\topp({\cal T'}_s))\equiv k_1q_ 1uq_2q_3\dots q_K k_2$ and $\Lab(\bott{\cal T'}_{s+1})\equiv k_1q_1q_2q_3\dots q_Kk_2.$

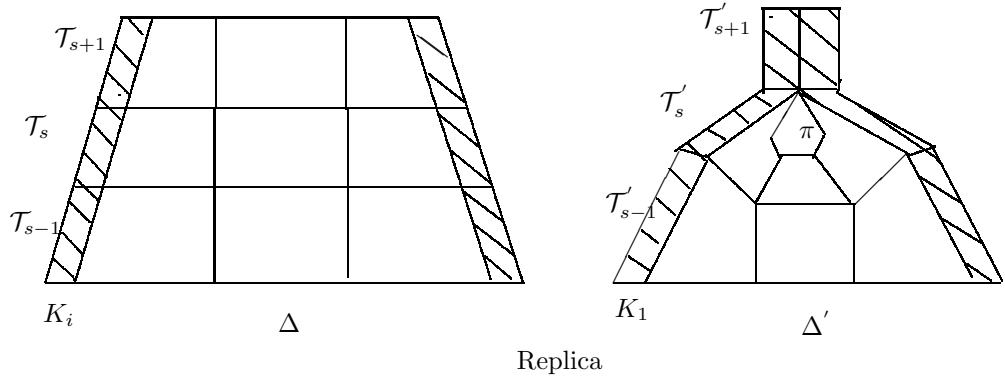
\begin{figure}[h!]

\unitlength 1mm 
\linethickness{0.4pt}
\ifx\plotpoint\undefined\newsavebox{\plotpoint}\fi 
\begin{picture}(139.5,55)(0,80)
\multiput(22.5,126)(-.0337171053,-.1159539474){304}{\line(0,-1){.1159539474}}
\put(12.25,90.75){\line(1,0){63.75}}
\multiput(75.75,90.75)(-.0336826347,.1055389222){334}{\line(0,1){.1055389222}}
\put(22.25,126){\line(1,0){42.25}}
\put(64.5,126){\line(0,1){0}}
\multiput(26.5,126)(-.0337171053,-.1151315789){304}{\line(0,-1){.1151315789}}
\multiput(16.25,91)(-.03125,.03125){8}{\line(0,1){.03125}}
\put(22.25,115.75){\line(-1,0){.25}}
\put(109.25,127){\rule{0\unitlength}{.25\unitlength}}
\put(108.5,126){\rule{.5\unitlength}{.25\unitlength}}
\put(107.75,127.25){\line(0,-1){11.25}}
\put(107.75,116){\line(0,1){0}}
\put(107.5,127.25){\line(1,0){10.25}}
\put(117.75,127.25){\line(1,0){.25}}
\put(117.75,127){\line(0,-1){10.25}}
\multiput(117.75,117.25)(.03125,.0625){8}{\line(0,1){.0625}}
\put(117.75,117.5){\line(0,-1){1}}
\put(117.75,118){\line(0,-1){1.5}}
\multiput(117.75,116.5)(-.03125,-.0625){8}{\line(0,-1){.0625}}
\multiput(117.5,116)(.060465116,-.03372093){215}{\line(1,0){.060465116}}
\put(130.5,108.75){\line(0,1){.25}}
\multiput(107.75,116.5)(-.0479591837,-.0336734694){245}{\line(-1,0){.0479591837}}
\put(112.5,127){\line(0,-1){10.5}}
\put(112.5,116.5){\line(0,1){0}}
\put(107.75,116.5){\line(1,0){9.75}}
\put(117.5,116.5){\line(0,1){0}}
\multiput(112.75,116.5)(-.0468164794,-.0337078652){267}{\line(-1,0){.0468164794}}
\multiput(100.25,107.5)(-.1166667,.0333333){30}{\line(-1,0){.1166667}}
\multiput(112.5,116)(.0602040816,-.0336734694){245}{\line(1,0){.0602040816}}
\multiput(127.25,107.75)(.1,.0333333){30}{\line(1,0){.1}}
\put(112.5,116.25){\line(-3,-5){3.75}}
\multiput(112.5,116)(.033505155,-.054123711){97}{\line(0,-1){.054123711}}
\multiput(115.75,110.75)(.03125,-.03125){8}{\line(0,-1){.03125}}
\multiput(108.75,110)(.03289474,-.07236842){38}{\line(0,-1){.07236842}}
\multiput(115.75,110.25)(-.03333333,-.06666667){45}{\line(0,-1){.06666667}}
\put(109.75,107.75){\line(1,0){4.5}}
\put(112.5,110){$\pi$}
\multiput(100,107.75)(.034974093,-.033678756){193}{\line(1,0){.034974093}}
\put(106.75,101.25){\line(1,0){13.25}}
\put(120,101.25){\line(1,1){6.75}}
\multiput(110.25,107.75)(-.033653846,-.060096154){104}{\line(0,-1){.060096154}}
\multiput(114.5,107.75)(.033687943,-.044326241){141}{\line(0,-1){.044326241}}
\put(96.5,108.25){\line(-1,-2){8.75}}
\put(87.75,90.75){\line(1,0){51.5}}
\multiput(130.25,108.75)(.0336363636,-.0645454545){275}{\line(0,-1){.0645454545}}
\multiput(126.75,108)(.0336538462,-.0634615385){260}{\line(0,-1){.0634615385}}
\multiput(100.25,107)(-.0336734694,-.0663265306){245}{\line(0,-1){.0663265306}}
\multiput(60.5,126)(.0336391437,-.1062691131){327}{\line(0,-1){.1062691131}}
\put(19.25,114){\line(1,0){49}}
\put(16.25,103.5){\line(1,0){55.25}}
\put(34.75,114){\line(0,-1){10.25}}
\put(52.5,113.75){\line(0,-1){10}}
\put(7.5,97.75){${\cal T}_{s-1}$}
\put(34.75,104){\line(0,-1){13}}
\put(52.5,104.5){\line(0,-1){13}}
\put(52.5,91.5){\line(0,1){0}}
\put(106.75,101){\line(0,-1){10}}
\put(119.75,101){\line(0,-1){10}}
\put(35,126){\line(0,-1){11.75}}
\put(52.25,126){\line(0,-1){12.25}}
\put(13.5,122.25){${\cal T}_{s+1}$}
\put(86.75,100.25){${\cal T}_{s-1}^{'}$}
\put(94,113.75){${\cal T}_s^{'}$}
\put(99.25,124.75){${\cal T}_{s+1}^{'}$}
\put(43.25,84.25){$\Delta$}
\put(112.75,83.75){$\Delta^{'}$}
\put(12,85.5){$K_i$}
\put(88,86.25){$K_1$}
\put(75,79.25){Replica}
\put(9.5,110.5){${\cal T}_s$}
\multiput(13.25,94)(.03666667,-.03333333){75}{\line(1,0){.03666667}}
\put(16,91.5){\line(0,1){0}}
\multiput(14.5,97.5)(.03353659,-.03658537){82}{\line(0,-1){.03658537}}
\put(17.25,94.5){\line(0,1){0}}
\multiput(15.5,101.5)(.03353659,-.03353659){82}{\line(0,-1){.03353659}}
\put(18.25,98.75){\line(0,1){0}}
\multiput(16.5,104.5)(.03353659,-.03353659){82}{\line(0,-1){.03353659}}
\put(19.25,101.75){\line(0,1){0}}
\multiput(17.5,108)(.03370787,-.03651685){89}{\line(0,-1){.03651685}}
\put(20.5,104.75){\line(0,1){0}}
\multiput(18.5,112)(.03353659,-.03658537){82}{\line(0,-1){.03658537}}
\put(21.25,109){\line(0,1){0}}
\multiput(19.25,115.5)(.03370787,-.03370787){89}{\line(0,-1){.03370787}}
\put(22.25,112.5){\line(0,1){0}}
\multiput(20.5,118.75)(.03370787,-.03370787){89}{\line(0,-1){.03370787}}
\put(23.5,115.75){\line(0,1){0}}
\multiput(21.5,122.5)(.03370787,-.03651685){89}{\line(0,-1){.03651685}}
\put(24.5,119.25){\line(0,1){0}}
\multiput(23.25,126)(.03333333,-.0375){60}{\line(0,-1){.0375}}
\put(25.25,123.75){\line(0,1){0}}
\multiput(70.75,94.25)(.03932584,-.03370787){89}{\line(1,0){.03932584}}
\multiput(69.25,98.5)(.043650794,-.033730159){126}{\line(1,0){.043650794}}
\put(74.75,94.25){\line(0,1){0}}
\multiput(68,102.75)(.044117647,-.033613445){119}{\line(1,0){.044117647}}
\put(73.25,98.75){\line(0,1){0}}
\multiput(66.25,108.25)(.041044776,-.03358209){134}{\line(1,0){.041044776}}
\put(71.75,103.75){\line(0,1){0}}
\multiput(64.75,113)(.041666667,-.033730159){126}{\line(1,0){.041666667}}
\put(70,108.75){\line(0,1){0}}
\multiput(63.25,118.5)(.042410714,-.033482143){112}{\line(1,0){.042410714}}
\put(68,114.75){\line(0,1){0}}
\put(61.75,123.75){\line(4,-3){4}}
\multiput(89.75,94.25)(.05,-.03333333){45}{\line(1,0){.05}}
\put(92,92.75){\line(0,1){0}}
\multiput(92,92.75)(.05,-.0333333){15}{\line(1,0){.05}}
\multiput(91.5,97.75)(.04326923,-.03365385){52}{\line(1,0){.04326923}}
\put(93.75,96){\line(0,1){0}}
\multiput(93.5,102)(.03731343,-.03358209){67}{\line(1,0){.03731343}}
\put(96,99.75){\line(0,1){0}}
\multiput(95.25,106)(.03353659,-.03353659){82}{\line(0,-1){.03353659}}
\put(98,103.25){\line(0,1){0}}
\multiput(97.75,109.25)(.03358209,-.03358209){67}{\line(0,-1){.03358209}}
\put(100,107){\line(0,1){0}}
\multiput(100.5,111)(.05921053,-.03289474){38}{\line(1,0){.05921053}}
\put(102.75,109.75){\line(0,1){0}}
\multiput(103.75,113.75)(.03846154,-.03365385){52}{\line(1,0){.03846154}}
\put(105.75,112){\line(0,1){0}}
\multiput(106.75,115.5)(.05263158,-.03289474){38}{\line(1,0){.05263158}}
\multiput(107.75,120)(.045731707,-.033536585){164}{\line(1,0){.045731707}}
\put(115.25,114.5){\line(0,1){0}}
\multiput(108,124.25)(.044186047,-.03372093){215}{\line(1,0){.044186047}}
\put(117.5,117){\line(0,1){0}}
\multiput(111.25,127)(.039634146,-.033536585){164}{\line(1,0){.039634146}}
\put(117.75,121.5){\line(0,1){0}}
\multiput(116.5,127)(.03289474,-.03289474){38}{\line(0,-1){.03289474}}
\put(117.75,125.75){\line(0,1){0}}
\multiput(120.5,114.25)(.049751244,-.03358209){201}{\line(1,0){.049751244}}
\multiput(130.5,107.5)(.0326087,-.0326087){23}{\line(0,-1){.0326087}}
\multiput(128.25,105)(.052884615,-.033653846){104}{\line(1,0){.052884615}}
\put(133.75,101.5){\line(0,1){0}}
\multiput(131.5,99.25)(.054123711,-.033505155){97}{\line(1,0){.054123711}}
\put(136.75,96){\line(0,1){0}}
\multiput(134.75,93.25)(.04850746,-.03358209){67}{\line(1,0){.04850746}}
\put(138,91){\line(0,1){0}}
\end{picture}

\end{figure}

The replica $\Delta'$ of the accepted $i$-sector $\Delta$ is constructed. To summarize our effort:
we cut off the $(\theta,a)$-cell and contract up the $a$-edges of $(\theta,q)$-cells for every command $\theta$ of the machine $\hat{\cal S}$ from $\Delta$,
replace the maximal $k_i$- and $k_{i+1}$-bands by their $k_1$- and $k_2$-copies, replace all the labels 
of the remaining edges by their copies from the alphabet ${\cal A}_1$, and
then close up all the holes by pasting in several $H$-cells. The result is the replica $\Delta',$
canonically obtained above. 

\begin{rk} \label{subsector}

{\rm (1) The $i$-sector $\Delta$ is a union of $K+1$ {\it subsectors} $\Gamma_1,\dots,\Gamma_{K+1},$ where 
every $\Gamma_j$ is a trapezia with a base of length $2$ and with height equal to the height of $\Delta.$
The subsector $\Gamma_j$ has a common maximal $q$-band ${\cal C}={\cal C}_{j+1}$ with $\Gamma_{j+1}$ for $j=1,\dots,K.$
If $i$ is odd (even), then $\Gamma_2$ (resp., $\Gamma_K$ ) is the {\it input subsector}. When we construct a replica,  $H$-cells appear only in the input subsector of the replica.
 
 Every $\Gamma_j$ has
its own replica $\Gamma'_j,$  which is a subsector of $\Delta'.$ Similarly, every $q$- or $\theta$-edge of $\cal C$
has a replica in the replica ${\cal C}'$ of ${\cal C}$ in $\Delta'$. Also every vertex $o$ of ${\cal C}$ belongs
to either a $q$-edge or a $\theta$-edge since the boundary of a $(\theta,q)$-cell has no two
consecutive $a$-edges by Lemma
\ref{MtocalS}(5); and so $o$ has a replica $o'$.

(2) The replica of an $i$-sector is not necessarily a minimal diagram.}

\end{rk} 

\subsection{Discs} \label{discs} 

Given a diagram $\Delta,$ one can construct a planar graph whose vertices
are the hubs of this diagram plus one improper vertex outside $\Delta,$ and the
edges are maximal $k$-bands of $\Delta.$  Lemma \ref{extdisc} says that if the 
number $L$ of sectors of a hub is large enough, then
the degree of the improper vertex is positive provided the
graph has at least one proper vertex; and therefore  the homomorphism $H\to G$ turns
out to be injective.

A {\it disc diagram (or a disc)} is a (sub)diagram $\Delta$ such that (1) it has
exactly one hub $\Pi$ (2) there are no $\theta$-edges on the boundary $\partial\Delta$ (3) there are
no $H$-cells of $\Delta$ having an edge on $\partial\Delta.$ In particular,
a hub is a disc diagram.

Let us consider a disc diagram $\Delta$ with hub $\Pi.$ 
 Denote by ${\cal K}_1,\dots,{\cal K}_L$ the 
$k_1,\dots, k_L$-bands starting on
the hub $\Pi.$ Since the hub relation has only one letter $k_i$
for every $i,$ these $k$-bands have to end on $\partial\Delta.$ It
therefore follows from Lemma \ref{NoAnnul} that for every $i,$ the
bands ${\cal K}_i$ and ${\cal K}_{i+1}$ bound, together with 
$\partial\Delta$ and $\partial\Pi,$ either a subdiagram $\Psi_i$
having no cells corresponding to any non-trivial relation of $G$
(in this case  the bands ${\cal K}_i$ and ${\cal K}_{i+1}$ are also trivial),
or $\Psi_i$ is a trapezium, and the boundary label of $\Delta$ has exactly
one letter from ${\cal Q}_i$ for every $i\in [1,\dots,N].$

Similarly, consider two hubs $\Pi_1$ and $\Pi_2$ in a minimal diagram, 
connected by a
$k_i$-band ${\cal K}_i$ and a $k_{i+1}$-band ${\cal K}_{i+1}$, where
$(i,i+1)\ne (1,2)$, and there are no other hubs between these $k$-bands. 
These bands, together with 
$\partial\Pi_1$ and $\partial\Pi_2,$ bound either a subdiagram $\Psi_i$
having no cells,  
or  a trapezium $\Psi_i$. The former case is impossible since in this case
the hubs have a common $k_i$-edge and they are mirror copies of each other
contrary to the reducibility of minimal diagrams. We want to show that
the latter case is not possible either.

Indeed, in the latter case $\Psi_i$ is an accepted trapezium since the $i$-sector
subword of $\Sigma_0$ is $k_iwk_{i+1},$ where $w$ is a (mirror) copy of
the accept word of the machine $\cal S.$ Therefore, according to Subsection
\ref{replicas}, a replica $\Psi_1$ of $\Psi_i$ (as well as the (mirror) copies
$\Psi_j$ for every $j\in [2,\dots,N]$) can be constructed. Then one can
construct a spherical diagram $\Gamma$ from $\Pi_1$, $\Pi_2$, and the diagrams
$\Psi_j$ ($j=1,\dots,N$). There are two subdiagrams $\Gamma'$ and $\Gamma''$ of
$\Gamma$ with common boundary: $\Gamma'$, being a copy of a subdiagram
of the original diagram, is made of $\Pi_1$, $\Pi_2$, and $\Psi=\Psi_i,$ 
and $\Gamma''$ is a union of all $\Psi_j$-s with $j\ne i$. Hence the subdiagram
$\Gamma'$ of the original diagram can be replaced by diagram of lower type with
the same boundary label because $\Gamma''$ has no hubs. This contradicts the
minimality of the original diagram.

Thus, any two hubs of a minimal diagram are connected by at most two $k$-bands,
such that the subdiagram bounded by them contains no other hubs.
This property makes the hub graph of a minimal diagram (where maximal $k$-bands play the role of edges
connecting hubs) hyperbolic (in a sense) since the degree $L$ of every vertex (=hub) is high ($\ge 40$).
Below we give a more precise formulation (proved for diagrams with such a hub graph, in particular, 
in \cite{O}, Lemma 3.2).

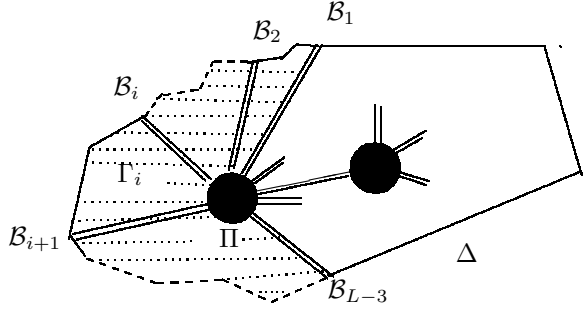
\begin{figure}[h!]

\unitlength 1mm 
\linethickness{0.4pt}
\ifx\plotpoint\undefined\newsavebox{\plotpoint}\fi 
\begin{picture}(94,45)(0,130)
\put(33,146){\circle*{.707}}
\put(30.6389,143.8889){\rule{4.7222\unitlength}{4.7222\unitlength}}
\multiput(31.6964,148.4987)(0,-5.3164){2}{\rule{2.6073\unitlength}{.8191\unitlength}}
\multiput(32.3078,149.2052)(0,-6.2084){2}{\rule{1.3843\unitlength}{.298\unitlength}}
\multiput(33.5797,149.2052)(-1.5821,0){2}{\multiput(0,0)(0,-6.1307){2}{\rule{.4228\unitlength}{.2202\unitlength}}}
\multiput(33.89,149.2052)(-2.0443,0){2}{\multiput(0,0)(0,-6.0805){2}{\rule{.2644\unitlength}{.17\unitlength}}}
\multiput(34.1912,148.4987)(-3.0583,0){2}{\multiput(0,0)(0,-5.0151){2}{\rule{.676\unitlength}{.5178\unitlength}}}
\multiput(34.1912,148.904)(-2.784,0){2}{\multiput(0,0)(0,-5.585){2}{\rule{.4016\unitlength}{.2769\unitlength}}}
\multiput(34.1912,149.0685)(-2.6411,0){2}{\multiput(0,0)(0,-5.8214){2}{\rule{.2587\unitlength}{.1844\unitlength}}}
\multiput(34.4803,148.904)(-3.2123,0){2}{\multiput(0,0)(0,-5.5061){2}{\rule{.2517\unitlength}{.1981\unitlength}}}
\multiput(34.7547,148.4987)(-3.8788,0){2}{\multiput(0,0)(0,-4.8246){2}{\rule{.3694\unitlength}{.3273\unitlength}}}
\multiput(34.7547,148.7135)(-3.7527,0){2}{\multiput(0,0)(0,-5.1378){2}{\rule{.2433\unitlength}{.2109\unitlength}}}
\multiput(35.0116,148.4987)(-4.2569,0){2}{\multiput(0,0)(0,-4.7201){2}{\rule{.2336\unitlength}{.2228\unitlength}}}
\multiput(35.2487,144.9464)(-5.3164,0){2}{\rule{.8191\unitlength}{2.6073\unitlength}}
\multiput(35.2487,147.4412)(-5.0151,0){2}{\multiput(0,0)(0,-3.0583){2}{\rule{.5178\unitlength}{.676\unitlength}}}
\multiput(35.2487,148.0047)(-4.8246,0){2}{\multiput(0,0)(0,-3.8788){2}{\rule{.3273\unitlength}{.3694\unitlength}}}
\multiput(35.2487,148.2616)(-4.7201,0){2}{\multiput(0,0)(0,-4.2569){2}{\rule{.2228\unitlength}{.2336\unitlength}}}
\multiput(35.4635,148.0047)(-5.1378,0){2}{\multiput(0,0)(0,-3.7527){2}{\rule{.2109\unitlength}{.2433\unitlength}}}
\multiput(35.654,147.4412)(-5.585,0){2}{\multiput(0,0)(0,-2.784){2}{\rule{.2769\unitlength}{.4016\unitlength}}}
\multiput(35.654,147.7303)(-5.5061,0){2}{\multiput(0,0)(0,-3.2123){2}{\rule{.1981\unitlength}{.2517\unitlength}}}
\multiput(35.8185,147.4412)(-5.8214,0){2}{\multiput(0,0)(0,-2.6411){2}{\rule{.1844\unitlength}{.2587\unitlength}}}
\multiput(35.9552,145.5578)(-6.2084,0){2}{\rule{.298\unitlength}{1.3843\unitlength}}
\multiput(35.9552,146.8297)(-6.1307,0){2}{\multiput(0,0)(0,-1.5821){2}{\rule{.2202\unitlength}{.4228\unitlength}}}
\multiput(35.9552,147.14)(-6.0805,0){2}{\multiput(0,0)(0,-2.0443){2}{\rule{.17\unitlength}{.2644\unitlength}}}
\put(36.26,146.25){\line(0,1){.3985}}
\put(36.235,146.649){\line(0,1){.1974}}
\put(36.205,146.846){\line(0,1){.1952}}
\put(36.162,147.041){\line(0,1){.1922}}
\put(36.108,147.233){\line(0,1){.1885}}
\put(36.042,147.422){\line(0,1){.1841}}
\put(35.964,147.606){\line(0,1){.179}}
\put(35.876,147.785){\line(0,1){.1732}}
\put(35.776,147.958){\line(0,1){.1668}}
\put(35.666,148.125){\line(0,1){.1598}}
\multiput(35.546,148.285)(-.03235,.03804){4}{\line(0,1){.03804}}
\multiput(35.417,148.437)(-.027697,.028788){5}{\line(0,1){.028788}}
\multiput(35.279,148.581)(-.029409,.027037){5}{\line(-1,0){.029409}}
\multiput(35.132,148.716)(-.03876,.03148){4}{\line(-1,0){.03876}}
\put(34.976,148.842){\line(-1,0){.1625}}
\put(34.814,148.958){\line(-1,0){.1693}}
\put(34.645,149.064){\line(-1,0){.1755}}
\put(34.469,149.16){\line(-1,0){.181}}
\put(34.288,149.244){\line(-1,0){.1858}}
\put(34.103,149.317){\line(-1,0){.19}}
\put(33.913,149.379){\line(-1,0){.1934}}
\put(33.719,149.429){\line(-1,0){.1961}}
\put(33.523,149.467){\line(-1,0){.198}}
\put(33.325,149.493){\line(-1,0){.1993}}
\put(33.126,149.507){\line(-1,0){.3992}}
\put(32.727,149.498){\line(-1,0){.1984}}
\put(32.528,149.475){\line(-1,0){.1967}}
\put(32.331,149.44){\line(-1,0){.1942}}
\put(32.137,149.393){\line(-1,0){.1909}}
\put(31.946,149.335){\line(-1,0){.187}}
\put(31.759,149.264){\line(-1,0){.1823}}
\put(31.577,149.183){\line(-1,0){.177}}
\put(31.4,149.09){\line(-1,0){.1709}}
\put(31.229,148.987){\line(-1,0){.1643}}
\multiput(31.065,148.873)(-.03926,-.03086){4}{\line(-1,0){.03926}}
\multiput(30.908,148.75)(-.03729,-.03321){4}{\line(-1,0){.03729}}
\multiput(30.759,148.617)(-.028152,-.028343){5}{\line(0,-1){.028343}}
\multiput(30.618,148.475)(-.03295,-.03752){4}{\line(0,-1){.03752}}
\multiput(30.486,148.325)(-.03059,-.03947){4}{\line(0,-1){.03947}}
\put(30.364,148.167){\line(0,-1){.1651}}
\put(30.251,148.002){\line(0,-1){.1716}}
\put(30.149,147.83){\line(0,-1){.1776}}
\put(30.058,147.653){\line(0,-1){.1828}}
\put(29.977,147.47){\line(0,-1){.1874}}
\put(29.908,147.283){\line(0,-1){.1913}}
\put(29.851,147.091){\line(0,-1){.1945}}
\put(29.805,146.897){\line(0,-1){.1969}}
\put(29.772,146.7){\line(0,-1){.5978}}
\put(29.744,146.102){\line(0,-1){.1992}}
\put(29.759,145.903){\line(0,-1){.1979}}
\put(29.786,145.705){\line(0,-1){.1958}}
\put(29.826,145.509){\line(0,-1){.193}}
\put(29.877,145.316){\line(0,-1){.1895}}
\put(29.94,145.127){\line(0,-1){.1853}}
\put(30.015,144.941){\line(0,-1){.1804}}
\put(30.1,144.761){\line(0,-1){.1748}}
\put(30.197,144.586){\line(0,-1){.1686}}
\put(30.304,144.418){\line(0,-1){.1617}}
\multiput(30.421,144.256)(.03174,-.03855){4}{\line(0,-1){.03855}}
\multiput(30.548,144.102)(.027236,-.029225){5}{\line(0,-1){.029225}}
\multiput(30.685,143.956)(.028975,-.027502){5}{\line(1,0){.028975}}
\multiput(30.829,143.818)(.03826,-.03209){4}{\line(1,0){.03826}}
\put(30.983,143.69){\line(1,0){.1606}}
\put(31.143,143.571){\line(1,0){.1676}}
\put(31.311,143.462){\line(1,0){.1739}}
\put(31.485,143.364){\line(1,0){.1796}}
\put(31.664,143.277){\line(1,0){.1846}}
\put(31.849,143.2){\line(1,0){.1889}}
\put(32.038,143.136){\line(1,0){.1926}}
\put(32.23,143.083){\line(1,0){.1954}}
\put(32.426,143.041){\line(1,0){.1976}}
\put(32.623,143.012){\line(1,0){.199}}
\put(32.822,142.995){\line(1,0){.3993}}
\put(33.222,142.998){\line(1,0){.1988}}
\put(33.42,143.018){\line(1,0){.1972}}
\put(33.618,143.049){\line(1,0){.1949}}
\put(33.813,143.093){\line(1,0){.1918}}
\put(34.004,143.149){\line(1,0){.188}}
\put(34.192,143.216){\line(1,0){.1836}}
\put(34.376,143.295){\line(1,0){.1784}}
\put(34.554,143.385){\line(1,0){.1726}}
\put(34.727,143.485){\line(1,0){.1661}}
\multiput(34.893,143.596)(.03974,.03023){4}{\line(1,0){.03974}}
\multiput(35.052,143.717)(.03782,.03261){4}{\line(1,0){.03782}}
\multiput(35.203,143.848)(.028599,.027892){5}{\line(1,0){.028599}}
\multiput(35.346,143.987)(.03355,.03699){4}{\line(0,1){.03699}}
\multiput(35.48,144.135)(.03122,.03897){4}{\line(0,1){.03897}}
\put(35.605,144.291){\line(0,1){.1633}}
\put(35.72,144.454){\line(0,1){.17}}
\put(35.825,144.624){\line(0,1){.1761}}
\put(35.92,144.8){\line(0,1){.1815}}
\put(36.003,144.982){\line(0,1){.1863}}
\put(36.075,145.168){\line(0,1){.1904}}
\put(36.135,145.359){\line(0,1){.1937}}
\put(36.184,145.552){\line(0,1){.1963}}
\put(36.221,145.749){\line(0,1){.1982}}
\put(36.245,145.947){\line(0,1){.3031}}
\multiput(34.75,149.25)(.0337171053,.0575657895){304}{\line(0,1){.0575657895}}
\multiput(34.25,149.5)(.0336700337,.0580808081){297}{\line(0,1){.0580808081}}
\multiput(33.25,150.25)(.033505155,.149484536){97}{\line(0,1){.149484536}}
\multiput(35.75,164.5)(-.033505155,-.149484536){97}{\line(0,-1){.149484536}}
\multiput(30.25,148.75)(-.0367063492,.0337301587){252}{\line(-1,0){.0367063492}}
\multiput(21,156.5)(.035714286,-.033613445){238}{\line(1,0){.035714286}}
\multiput(11.75,141.5)(.151260504,.033613445){119}{\line(1,0){.151260504}}
\multiput(12.25,140.75)(.149159664,.033613445){119}{\line(1,0){.149159664}}
\multiput(35.75,144.25)(.0426587302,-.0337301587){252}{\line(1,0){.0426587302}}
\put(46,135.75){\line(0,1){0}}
\multiput(35.5,143.75)(.0418367347,-.0336734694){245}{\line(1,0){.0418367347}}
\multiput(35.5,148.5)(.046391753,.033505155){97}{\line(1,0){.046391753}}
\multiput(36,148.25)(.04573171,.03353659){82}{\line(1,0){.04573171}}
\put(36.25,147.25){\line(5,1){12.5}}
\multiput(36.25,146.75)(.16333333,.03333333){75}{\line(1,0){.16333333}}
\put(36,146.25){\line(1,0){6.25}}
\put(36.25,145.5){\line(1,0){6}}
\put(51.75,150.25){\circle*{.707}}
\put(49.6186,147.8686){\rule{4.7627\unitlength}{4.7627\unitlength}}
\multiput(50.6854,152.5189)(0,-5.3631){2}{\rule{2.6292\unitlength}{.8253\unitlength}}
\multiput(51.3023,153.2317)(0,-6.263){2}{\rule{1.3955\unitlength}{.2996\unitlength}}
\multiput(52.5853,153.2317)(-1.596,0){2}{\multiput(0,0)(0,-6.1846){2}{\rule{.4255\unitlength}{.2212\unitlength}}}
\multiput(52.8983,153.2317)(-2.0623,0){2}{\multiput(0,0)(0,-6.134){2}{\rule{.2657\unitlength}{.1706\unitlength}}}
\multiput(53.2021,152.5189)(-3.0852,0){2}{\multiput(0,0)(0,-5.0592){2}{\rule{.681\unitlength}{.5214\unitlength}}}
\multiput(53.2021,152.9278)(-2.8084,0){2}{\multiput(0,0)(0,-5.634){2}{\rule{.4042\unitlength}{.2784\unitlength}}}
\multiput(53.2021,153.0937)(-2.6643,0){2}{\multiput(0,0)(0,-5.8725){2}{\rule{.26\unitlength}{.185\unitlength}}}
\multiput(53.4938,152.9278)(-3.2406,0){2}{\multiput(0,0)(0,-5.5545){2}{\rule{.2529\unitlength}{.1988\unitlength}}}
\multiput(53.7706,152.5189)(-3.9129,0){2}{\multiput(0,0)(0,-4.867){2}{\rule{.3717\unitlength}{.3292\unitlength}}}
\multiput(53.7706,152.7356)(-3.7856,0){2}{\multiput(0,0)(0,-5.183){2}{\rule{.2444\unitlength}{.2118\unitlength}}}
\multiput(54.0298,152.5189)(-4.2943,0){2}{\multiput(0,0)(0,-4.7616){2}{\rule{.2347\unitlength}{.2238\unitlength}}}
\multiput(54.152,152.5189)(-4.4756,0){2}{\multiput(0,0)(0,-4.7066){2}{\rule{.1716\unitlength}{.1688\unitlength}}}
\multiput(54.2689,148.9354)(-5.3631,0){2}{\rule{.8253\unitlength}{2.6292\unitlength}}
\multiput(54.2689,151.4521)(-5.0592,0){2}{\multiput(0,0)(0,-3.0852){2}{\rule{.5214\unitlength}{.681\unitlength}}}
\multiput(54.2689,152.0206)(-4.867,0){2}{\multiput(0,0)(0,-3.9129){2}{\rule{.3292\unitlength}{.3717\unitlength}}}
\multiput(54.2689,152.2798)(-4.7616,0){2}{\multiput(0,0)(0,-4.2943){2}{\rule{.2238\unitlength}{.2347\unitlength}}}
\multiput(54.2689,152.402)(-4.7066,0){2}{\multiput(0,0)(0,-4.4756){2}{\rule{.1688\unitlength}{.1716\unitlength}}}
\multiput(54.4856,152.0206)(-5.183,0){2}{\multiput(0,0)(0,-3.7856){2}{\rule{.2118\unitlength}{.2444\unitlength}}}
\multiput(54.6778,151.4521)(-5.634,0){2}{\multiput(0,0)(0,-2.8084){2}{\rule{.2784\unitlength}{.4042\unitlength}}}
\multiput(54.6778,151.7438)(-5.5545,0){2}{\multiput(0,0)(0,-3.2406){2}{\rule{.1988\unitlength}{.2529\unitlength}}}
\multiput(54.8437,151.4521)(-5.8725,0){2}{\multiput(0,0)(0,-2.6643){2}{\rule{.185\unitlength}{.26\unitlength}}}
\multiput(54.9817,149.5523)(-6.263,0){2}{\rule{.2996\unitlength}{1.3955\unitlength}}
\multiput(54.9817,150.8353)(-6.1846,0){2}{\multiput(0,0)(0,-1.596){2}{\rule{.2212\unitlength}{.4255\unitlength}}}
\multiput(54.9817,151.1483)(-6.134,0){2}{\multiput(0,0)(0,-2.0623){2}{\rule{.1706\unitlength}{.2657\unitlength}}}
\put(55.288,150.25){\line(0,1){.4018}}
\put(55.264,150.652){\line(0,1){.199}}
\put(55.233,150.851){\line(0,1){.1968}}
\put(55.19,151.048){\line(0,1){.1938}}
\put(55.135,151.241){\line(0,1){.1901}}
\put(55.069,151.431){\line(0,1){.1856}}
\put(54.991,151.617){\line(0,1){.1805}}
\put(54.901,151.798){\line(0,1){.1747}}
\put(54.801,151.972){\line(0,1){.1682}}
\multiput(54.69,152.141)(-.03019,.04029){4}{\line(0,1){.04029}}
\multiput(54.57,152.302)(-.0326,.03836){4}{\line(0,1){.03836}}
\multiput(54.439,152.455)(-.027913,.029036){5}{\line(0,1){.029036}}
\multiput(54.3,152.6)(-.029638,.027273){5}{\line(-1,0){.029638}}
\multiput(54.151,152.737)(-.03906,.03176){4}{\line(-1,0){.03906}}
\put(53.995,152.864){\line(-1,0){.1637}}
\put(53.831,152.981){\line(-1,0){.1706}}
\put(53.661,153.088){\line(-1,0){.1768}}
\put(53.484,153.184){\line(-1,0){.1824}}
\put(53.302,153.27){\line(-1,0){.1873}}
\put(53.114,153.344){\line(-1,0){.1915}}
\put(52.923,153.406){\line(-1,0){.1949}}
\put(52.728,153.457){\line(-1,0){.1977}}
\put(52.53,153.495){\line(-1,0){.1997}}
\put(52.331,153.522){\line(-1,0){.2009}}
\put(52.13,153.536){\line(-1,0){.4025}}
\put(51.727,153.527){\line(-1,0){.2001}}
\put(51.527,153.504){\line(-1,0){.1983}}
\put(51.329,153.469){\line(-1,0){.1958}}
\put(51.133,153.422){\line(-1,0){.1925}}
\put(50.94,153.363){\line(-1,0){.1886}}
\put(50.752,153.292){\line(-1,0){.1839}}
\put(50.568,153.21){\line(-1,0){.1785}}
\put(50.39,153.117){\line(-1,0){.1725}}
\put(50.217,153.013){\line(-1,0){.1658}}
\multiput(50.051,152.899)(-.03962,-.03107){4}{\line(-1,0){.03962}}
\multiput(49.893,152.774)(-.03764,-.03343){4}{\line(-1,0){.03764}}
\multiput(49.742,152.641)(-.028421,-.028539){5}{\line(0,-1){.028539}}
\multiput(49.6,152.498)(-.03328,-.03778){4}{\line(0,-1){.03778}}
\multiput(49.467,152.347)(-.0309,-.03975){4}{\line(0,-1){.03975}}
\put(49.343,152.188){\line(0,-1){.1663}}
\put(49.23,152.021){\line(0,-1){.1729}}
\put(49.127,151.849){\line(0,-1){.1789}}
\put(49.034,151.67){\line(0,-1){.1842}}
\put(48.953,151.485){\line(0,-1){.1889}}
\put(48.883,151.297){\line(0,-1){.1928}}
\put(48.825,151.104){\line(0,-1){.196}}
\put(48.778,150.908){\line(0,-1){.1985}}
\put(48.744,150.709){\line(0,-1){.6027}}
\put(48.715,150.107){\line(0,-1){.2008}}
\put(48.73,149.906){\line(0,-1){.1995}}
\put(48.757,149.706){\line(0,-1){.1975}}
\put(48.796,149.509){\line(0,-1){.1947}}
\put(48.848,149.314){\line(0,-1){.1912}}
\put(48.911,149.123){\line(0,-1){.187}}
\put(48.986,148.936){\line(0,-1){.1821}}
\put(49.072,148.754){\line(0,-1){.1764}}
\put(49.169,148.577){\line(0,-1){.1702}}
\put(49.277,148.407){\line(0,-1){.1633}}
\multiput(49.395,148.244)(.03192,-.03893){4}{\line(0,-1){.03893}}
\multiput(49.522,148.088)(.027396,-.029525){5}{\line(0,-1){.029525}}
\multiput(49.659,147.941)(.029152,-.027792){5}{\line(1,0){.029152}}
\multiput(49.805,147.802)(.0385,-.03244){4}{\line(1,0){.0385}}
\multiput(49.959,147.672)(.04041,-.03003){4}{\line(1,0){.04041}}
\put(50.121,147.552){\line(1,0){.1687}}
\put(50.289,147.442){\line(1,0){.1751}}
\put(50.464,147.342){\line(1,0){.1809}}
\put(50.645,147.254){\line(1,0){.186}}
\put(50.831,147.176){\line(1,0){.1903}}
\put(51.022,147.111){\line(1,0){.194}}
\put(51.216,147.057){\line(1,0){.197}}
\put(51.413,147.015){\line(1,0){.1992}}
\put(51.612,146.985){\line(1,0){.2006}}
\put(51.812,146.967){\line(1,0){.4026}}
\put(52.215,146.969){\line(1,0){.2005}}
\put(52.415,146.988){\line(1,0){.1989}}
\put(52.614,147.02){\line(1,0){.1966}}
\put(52.811,147.063){\line(1,0){.1936}}
\put(53.004,147.119){\line(1,0){.1898}}
\put(53.194,147.186){\line(1,0){.1853}}
\put(53.379,147.265){\line(1,0){.1801}}
\put(53.56,147.355){\line(1,0){.1743}}
\put(53.734,147.456){\line(1,0){.1678}}
\multiput(53.902,147.567)(.04016,.03036){4}{\line(1,0){.04016}}
\multiput(54.062,147.689)(.03823,.03276){4}{\line(1,0){.03823}}
\multiput(54.215,147.82)(.028921,.028033){5}{\line(1,0){.028921}}
\multiput(54.36,147.96)(.02715,.029751){5}{\line(0,1){.029751}}
\multiput(54.496,148.109)(.0316,.0392){4}{\line(0,1){.0392}}
\put(54.622,148.266){\line(0,1){.1642}}
\put(54.739,148.43){\line(0,1){.1711}}
\put(54.845,148.601){\line(0,1){.1772}}
\put(54.94,148.778){\line(0,1){.1828}}
\put(55.025,148.961){\line(0,1){.1876}}
\put(55.098,149.149){\line(0,1){.1917}}
\put(55.16,149.34){\line(0,1){.1951}}
\put(55.21,149.535){\line(0,1){.1978}}
\put(55.247,149.733){\line(0,1){.1998}}
\put(55.273,149.933){\line(0,1){.317}}
\put(52,153.25){\line(0,1){5.5}}
\put(52.75,158.5){\line(0,-1){4.75}}
\multiput(54.5,152.75)(.05666667,.03333333){75}{\line(1,0){.05666667}}
\multiput(54.75,152.25)(.05833333,.03333333){60}{\line(1,0){.05833333}}
\multiput(55,149.75)(.11184211,-.03289474){38}{\line(1,0){.11184211}}
\multiput(55,149.25)(.09868421,-.03289474){38}{\line(1,0){.09868421}}
\multiput(44.5,166.5)(-.375,.03125){8}{\line(-1,0){.375}}
\multiput(41.5,166.75)(-.03333333,-.03333333){60}{\line(0,-1){.03333333}}
\multiput(39.75,165)(-.2065217,-.0326087){23}{\line(-1,0){.2065217}}
\put(34.93,164.43){\line(-1,0){.85}}
\put(33.23,164.43){\line(-1,0){.85}}
\put(31.53,164.43){\line(-1,0){.85}}
\multiput(30.43,164.43)(-.033333,-.070833){10}{\line(0,-1){.070833}}
\multiput(29.763,163.013)(-.033333,-.070833){10}{\line(0,-1){.070833}}
\multiput(29.096,161.596)(-.033333,-.070833){10}{\line(0,-1){.070833}}
\multiput(28.43,160.68)(-.18,-.03){5}{\line(-1,0){.18}}
\multiput(26.63,160.38)(-.18,-.03){5}{\line(-1,0){.18}}
\multiput(24.83,160.08)(-.18,-.03){5}{\line(-1,0){.18}}
\multiput(23.93,159.93)(-.03125,.03125){4}{\line(0,1){.03125}}
\multiput(23.68,160.18)(-.0321429,-.0428571){14}{\line(0,-1){.0428571}}
\multiput(22.78,158.98)(-.0321429,-.0428571){14}{\line(0,-1){.0428571}}
\multiput(21.88,157.78)(-.0321429,-.0428571){14}{\line(0,-1){.0428571}}
\multiput(21.5,157.25)(-.05952381,-.033730159){126}{\line(-1,0){.05952381}}
\multiput(14,153)(-.03353659,-.14939024){82}{\line(0,-1){.14939024}}
\multiput(11.5,141.25)(.03333333,-.04583333){60}{\line(0,-1){.04583333}}
\multiput(13.68,138.18)(.093434,-.032828){9}{\line(1,0){.093434}}
\multiput(15.362,137.589)(.093434,-.032828){9}{\line(1,0){.093434}}
\multiput(17.043,136.998)(.093434,-.032828){9}{\line(1,0){.093434}}
\multiput(18.725,136.407)(.093434,-.032828){9}{\line(1,0){.093434}}
\multiput(20.407,135.816)(.093434,-.032828){9}{\line(1,0){.093434}}
\multiput(22.089,135.225)(.093434,-.032828){9}{\line(1,0){.093434}}
\put(22.93,134.93){\line(1,0){.9722}}
\put(24.874,135.041){\line(1,0){.9722}}
\put(26.819,135.152){\line(1,0){.9722}}
\put(28.763,135.263){\line(1,0){.9722}}
\put(30.707,135.374){\line(1,0){.9722}}
\multiput(31.68,135.43)(.067708,-.03125){12}{\line(1,0){.067708}}
\multiput(33.305,134.68)(.067708,-.03125){12}{\line(1,0){.067708}}
\multiput(34.93,133.93)(.067708,-.03125){12}{\line(1,0){.067708}}
\multiput(36.555,133.18)(.067708,-.03125){12}{\line(1,0){.067708}}
\multiput(38.18,132.43)(.072727,.031818){11}{\line(1,0){.072727}}
\multiput(39.78,133.13)(.072727,.031818){11}{\line(1,0){.072727}}
\multiput(41.38,133.83)(.072727,.031818){11}{\line(1,0){.072727}}
\multiput(42.98,134.53)(.072727,.031818){11}{\line(1,0){.072727}}
\multiput(44.58,135.23)(.072727,.031818){11}{\line(1,0){.072727}}
\put(46.18,135.93){\line(-1,0){.125}}
\multiput(45.93,135.93)(-.03125,.03125){8}{\line(0,1){.03125}}
\multiput(76.75,166.25)(-.03125,.03125){8}{\line(0,1){.03125}}
\put(44.75,166.5){\line(1,0){29.5}}
\multiput(74.5,166.5)(.033557047,-.115771812){149}{\line(0,-1){.115771812}}
\multiput(79.25,149.75)(-.081495098,-.0337009804){408}{\line(-1,0){.081495098}}
\put(31.25,139.75){$\Pi$}
\put(62.75,137.75){$\Delta$}
\put(17.25,160){${\cal B}_i$}
\put(3.25,140){${\cal B}_{i+1}$}
\put(45.5,170){${\cal B}_1$}
\put(35.75,167.5){${\cal B}_2$}
\put(45.75,133){${\cal B}_{L-3}$}
\multiput(30.18,163.18)(.97917,-.02083){13}{{\rule{.4pt}{.4pt}}}
\multiput(28.93,161.18)(.97917,-.02083){13}{{\rule{.4pt}{.4pt}}}
\multiput(22.93,159.43)(.97059,-.02941){18}{{\rule{.4pt}{.4pt}}}
\multiput(21.43,157.18)(.97059,0){18}{{\rule{.4pt}{.4pt}}}
\multiput(17.68,154.93)(.9625,-.0125){21}{{\rule{.4pt}{.4pt}}}
\multiput(14.43,152.68)(.95455,0){23}{{\rule{.4pt}{.4pt}}}
\multiput(14.43,150.93)(.97368,-.02632){20}{{\rule{.4pt}{.4pt}}}
\multiput(24.43,148.18)(.91667,-.04167){7}{{\rule{.4pt}{.4pt}}}
\multiput(12.43,145.68)(.97059,0){18}{{\rule{.4pt}{.4pt}}}
\multiput(11.68,143.18)(.95,.025){21}{{\rule{.4pt}{.4pt}}}
\multiput(12.18,140.18)(.95313,.01563){17}{{\rule{.4pt}{.4pt}}}
\multiput(34.43,140.68)(.95,.05){6}{{\rule{.4pt}{.4pt}}}
\multiput(15.43,137.93)(.96296,.00926){28}{{\rule{.4pt}{.4pt}}}
\multiput(32.43,135.43)(.96154,.07692){14}{{\rule{.4pt}{.4pt}}}
\multiput(36.68,133.68)(.75,0){5}{{\rule{.4pt}{.4pt}}}
\multiput(12.18,172.93)(.125,-.125){3}{{\rule{.4pt}{.4pt}}}
\put(17.75,148.5){$\Gamma_i$}
\end{picture}

\end{figure}

\begin{lemma} \label{extdisc} If a minimal diagram over the group $G$ contains a least one hub,
then there is a hub $\Pi$ in $\Delta$ such that $L-3$ consecutive maximal $k$-bands ${\cal B}_1,\dots
{\cal B}_{L-3} $ start on $\Pi$, end on the boundary $\partial\Delta$, and for any $i\in [1,L-4]$, 
there are no discs in the subdiagram $\Gamma_i$ bounded by ${\cal B}_i$, ${\cal B}_{i+1},$ $\partial\Pi,$ and
$\partial\Delta.$ 
\end{lemma}

\begin{cintro} \label{inject} The canonical homomorphism $H\to G$ given by Lemma \ref{homo} is injective. 
\end{cintro}  

\proof Assume that a word $w$ in the generators ${a_1,...,a_m}$ of the group $H$
is equal to $1$ in $G.$ Then by van Kampen's Lemma, there is a minimal diagram $\Delta$
over $G$ whose boundary label is $w$.
Since $w$ has no $q$-letters, $\Delta$ has no hubs by Lemma \ref{extdisc}. By Lemma \ref{NoAnnul}, $\Delta$ 
contains neither $q$- nor $\theta$-annuli, and so it has neither $(\theta,q)$-cells
nor $(\theta, a)$-cells, because $w$ has neither $q$- nor $\theta$-letters. Hence this diagram can contain $H$-cells only. Since the
boundary labels of $H$-cells are trivial in $H,$ the boundary label $w$ is trivial in $H$
too by van Kampen's Lemma, and so the homomorphism is injective. \endproof.


\section{Comparison of paths in diagrams } \label{compare}

The main lemma of this section is Lemma \ref{cut} which says
that given a diagram $\Delta,$ one can cut off a subdiagram 
having one hub (and cells of lower rank) such that the perimeter
of the remaining diagram does not exceed $|\partial\Delta|.$
This makes it possible to use induction on the number of hubs in the
next section. To carry out the surgery of Lemma \ref{cut} we
must be able to compare the lengths of paths in various 
particular diagrams and to  `finish' some `unfinished' diagrams. 

\subsection{Paths in sectors}\label{paths}

We will consider the words in the generators
of the group $G$ and will modify the length function on this set. This modification is helpful in subsequent subsections,
in particular, we cannot prove an analog of Lemma \ref{rim} using the
standard  length $||\;||.$ 

The
standard length $||\;||$ of a word (of a path) will be called its {\em
combinatorial length}. From now on we use the word {\em length} for
the modified length. We set the length of every
$q$-letter equal 1, and the length of every $a$-letter equal to a small
enough number $\delta>0$ so that

\begin{equation}
\delta < (3N)^{-1}. \label{param}
\end{equation} 

If a word $v$ has $s$ $\theta$-letters, $t$ $a$-letters, and no $q$-letters,
then  $$|v|=s + \delta \max(0, t-s)$$ 
by definition.  For example, the word read between two $q$-letters of
a $(q,\theta)$-relation has length $1$ since it has one $\theta$-letter
and at most one $a$-letter by formulas (\ref{rel1}, \ref{rel2}) and property (3) of Lemma \ref{MtocalS}.

An arbitrary word $w$ is a product $v_0q_1v_1q_2\dots q_m v_{m},$
where $q_1,\dots,q_{m}$ are $q$-letters and the words $v_0,\dots v_{m}$ have
no $q$-letters. Then, by definition, $|w|=m+\sum_{j=0}^{m}|v_j|.$
{\em The length of a path} in a diagram is the
length of its label. The {\em perimeter} $|\partial\Delta|$ of a
van Kampen diagram is similarly defined by a shortest cyclic decomposition
of the boundary $\partial\Delta$. It follows from this definition
that for any product $s=s_1s_2$ of two words or paths, we
have $ |s|\le |s_1|+|s_2|$, and $|s|=|s_1|+|s_2|$ if $s_2$ starts or $s_1$
ends with a $q$-letter.

If a path $p$ starts at a vertex $o$ and ends at $o'$, we will write $o=p_-$ and $o'=p_+.$

\begin{lemma}\label{prohod} Let $\Delta$ be a trapezium bounded by two maximal $q$-bands $\cal C$
and $\cal C'$ and having no $H$-cells. Assume that 
 $\cal C$ has no $a$-edges. 
Let $o$  be vertex lying on both $\cal C$ and the top of the trapezium $\Delta$, and let $o'$
belong to $\cal C'$ and to the bottom of $\Delta.$ Assume that a path $t$ connects $o$ and $o'$ and
has  no $q$-edges.
Then the vertices $o$ and $o'$ can be connected in $\Delta$ by a path $t'$ such that $|t'|\le |t|$
and $t'=t_1t_2t_3$ where $t_1$ and $t_3$ are parts of the sides of $\cal C$ and $\cal C',$
respectively, and $t_2$ consists only of $a$-edges.
\end{lemma}
\proof Since the path $t$ has no $q$-edges, it follows from the assumption of the lemma that every   maximal $\theta$-band $\cal T$ of $\Delta$ has exactly two $(\theta,q)$-cells (the first one and the last one). Therefore
$\cal C$ and $\cal C'$ can be connected along $\cal T$ by a path $x$ consisting of $a$-edges only.

We denote by $t_2$ a shortest  path among such  $x$-s.
Then we define
$t_1$ ($t_3$) as the shortest subpath of the side of $\cal C$ (of $\cal C'$) connecting
$o$ and $(t_2)_-$ ($(t_2)_+$ and $o'$). 

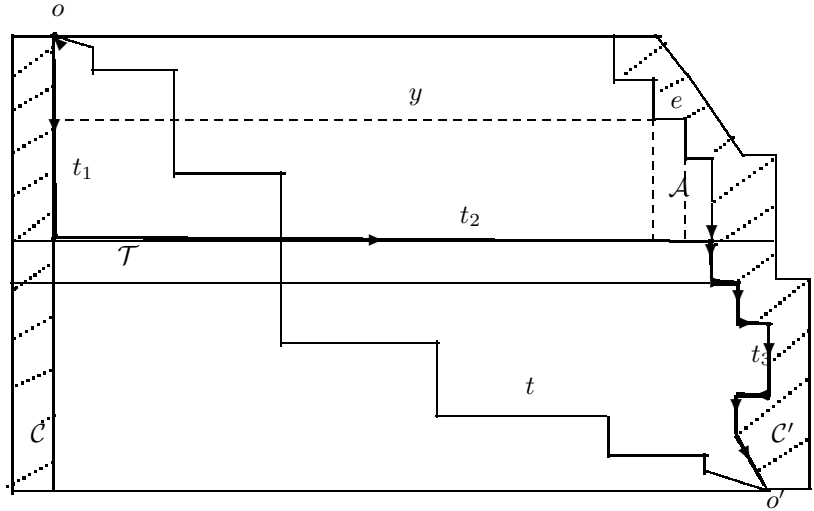
\begin{figure}[h!]
\unitlength 1mm 
\linethickness{0.4pt}
\ifx\plotpoint\undefined\newsavebox{\plotpoint}\fi 
\begin{picture}(248.5,90)(0,90)
\put(12,173.75){\line(0,-1){61.25}}
\put(12,173.5){\line(1,0){86.25}}
\put(92,173.75){\line(0,-1){6.5}}
\multiput(97.5,173.75)(.03358209,-.042910448){134}{\line(0,-1){.042910448}}
\put(92,167.75){\line(1,0){5.75}}
\put(97.25,168){\line(0,-1){5.5}}
\put(97.25,162.5){\line(1,0){4}}
\put(101.5,162.75){\line(0,-1){5.75}}
\put(101.5,157.25){\line(1,0){3.75}}
\put(105,157.25){\line(0,-1){5.75}}
\put(105,151.5){\line(0,1){1}}
\put(105,151.25){\line(0,-1){6.25}}
\put(105,145){\line(0,-1){4.75}}
\put(105,140.5){\line(1,0){3.5}}
\put(108.5,140.5){\line(0,-1){5.25}}
\put(108.5,135.25){\line(1,0){4.5}}
\put(112.5,135){\line(0,-1){9.75}}
\put(112.25,125.5){\line(-1,0){4}}
\put(108.25,125.5){\line(0,-1){5.25}}
\put(12.25,113){\line(1,0){100.5}}
\multiput(101.75,168.5)(.033632287,-.049327354){223}{\line(0,-1){.049327354}}
\put(109.25,157.75){\line(1,0){4.25}}
\put(113.5,157.75){\line(0,-1){5.5}}
\put(113.5,152.25){\line(0,-1){11}}
\put(113.5,141.25){\line(1,0){4.5}}
\put(118,141.25){\line(0,-1){27.75}}
\put(112,113.25){\line(1,0){6.25}}
\multiput(108.25,120.75)(.033613445,-.06092437){119}{\line(0,-1){.06092437}}
\put(17.5,173.5){\line(0,-1){60.25}}
\put(104.75,140.75){\line(-1,0){87.5}}
\put(17.25,140.75){\line(0,1){0}}
\put(17.5,140.75){\line(-1,0){5.75}}
\put(12,146.25){\line(1,0){101.25}}
\put(91.93,167.43){\line(0,1){.125}}
\put(97.18,162.68){\line(0,-1){.9706}}
\put(97.18,160.739){\line(0,-1){.9706}}
\put(97.18,158.797){\line(0,-1){.9706}}
\put(97.18,156.856){\line(0,-1){.9706}}
\put(97.18,154.915){\line(0,-1){.9706}}
\put(97.18,152.974){\line(0,-1){.9706}}
\put(97.18,151.033){\line(0,-1){.9706}}
\put(97.18,149.091){\line(0,-1){.9706}}
\put(97.18,147.15){\line(0,-1){.9706}}
\put(101.43,157.43){\line(0,-1){.9773}}
\put(101.43,155.475){\line(0,-1){.9773}}
\put(101.43,153.521){\line(0,-1){.9773}}
\put(101.43,151.566){\line(0,-1){.9773}}
\put(101.43,149.612){\line(0,-1){.9773}}
\put(101.43,147.657){\line(0,-1){.9773}}
\put(97.43,162.43){\line(-1,0){.9907}}
\put(95.448,162.43){\line(-1,0){.9907}}
\put(93.467,162.43){\line(-1,0){.9907}}
\put(91.485,162.43){\line(-1,0){.9907}}
\put(89.504,162.43){\line(-1,0){.9907}}
\put(87.522,162.43){\line(-1,0){.9907}}
\put(85.541,162.43){\line(-1,0){.9907}}
\put(83.559,162.43){\line(-1,0){.9907}}
\put(81.578,162.43){\line(-1,0){.9907}}
\put(79.596,162.43){\line(-1,0){.9907}}
\put(77.615,162.43){\line(-1,0){.9907}}
\put(75.633,162.43){\line(-1,0){.9907}}
\put(73.652,162.43){\line(-1,0){.9907}}
\put(71.67,162.43){\line(-1,0){.9907}}
\put(69.689,162.43){\line(-1,0){.9907}}
\put(67.707,162.43){\line(-1,0){.9907}}
\put(65.726,162.43){\line(-1,0){.9907}}
\put(63.745,162.43){\line(-1,0){.9907}}
\put(61.763,162.43){\line(-1,0){.9907}}
\put(59.782,162.43){\line(-1,0){.9907}}
\put(57.8,162.43){\line(-1,0){.9907}}
\put(55.819,162.43){\line(-1,0){.9907}}
\put(53.837,162.43){\line(-1,0){.9907}}
\put(51.856,162.43){\line(-1,0){.9907}}
\put(49.874,162.43){\line(-1,0){.9907}}
\put(47.893,162.43){\line(-1,0){.9907}}
\put(45.911,162.43){\line(-1,0){.9907}}
\put(43.93,162.43){\line(-1,0){.9907}}
\put(41.948,162.43){\line(-1,0){.9907}}
\put(39.967,162.43){\line(-1,0){.9907}}
\put(37.985,162.43){\line(-1,0){.9907}}
\put(36.004,162.43){\line(-1,0){.9907}}
\put(34.022,162.43){\line(-1,0){.9907}}
\put(32.041,162.43){\line(-1,0){.9907}}
\put(30.059,162.43){\line(-1,0){.9907}}
\put(28.078,162.43){\line(-1,0){.9907}}
\put(26.096,162.43){\line(-1,0){.9907}}
\put(24.115,162.43){\line(-1,0){.9907}}
\put(22.133,162.43){\line(-1,0){.9907}}
\put(20.152,162.43){\line(-1,0){.9907}}
\put(18.17,162.43){\line(-1,0){.9907}}
\multiput(17.75,173.5)(.11111111,-.03333333){45}{\line(1,0){.11111111}}
\put(22.75,172){\line(0,-1){3.5}}
\put(22.5,169){\line(1,0){11.5}}
\put(33.5,169){\line(0,-1){14.25}}
\put(33.5,155.25){\line(1,0){14.75}}
\put(47.75,155.5){\line(0,-1){23.25}}
\put(47.75,132.75){\line(1,0){21}}
\put(68.5,133){\line(0,-1){10.25}}
\put(68.25,123){\line(1,0){22.75}}
\put(91.25,123){\line(0,-1){5.5}}
\put(91.5,117.75){\line(1,0){13}}
\put(104,118){\line(0,-1){2.75}}
\multiput(104,115.75)(.10365854,-.03353659){82}{\line(1,0){.10365854}}
\multiput(17.5,173.5)(.03125,-.03125){8}{\line(0,-1){.03125}}
\put(99.25,152.25){$\cal A$}
\put(99.5,164){$e$}
\put(64.75,165){$y$}
\put(71.5,148.75){$t_2$}
\put(20,155){$t_1$}
\put(110.25,130.25){$t_3$}
\put(17.25,176.25){$o$}
\put(112.25,110.75){$o'$}
\put(14.5,119.5){$\cal C$}
\put(113,119.5){$\cal C '$}
\put(26,143.25){$\cal T$}
\put(80.25,126){$t$}
\thicklines
\multiput(17.5,173.75)(.03125,.03125){8}{\line(0,1){.03125}}
\put(248,120){\line(-1,0){.25}}
\put(17.25,171.25){\line(0,1){.25}}
\multiput(17.18,171.93)(-.71429,-.53571){8}{{\rule{.8pt}{.8pt}}}
\multiput(17.18,166.68)(-.71429,-.42857){8}{{\rule{.8pt}{.8pt}}}
\multiput(17.18,160.93)(-.75,-.42857){8}{{\rule{.8pt}{.8pt}}}
\multiput(16.43,154.93)(-.75,-.45833){7}{{\rule{.8pt}{.8pt}}}
\multiput(16.68,149.68)(-.64286,-.5){8}{{\rule{.8pt}{.8pt}}}
\multiput(16.68,142.93)(-.60714,-.42857){8}{{\rule{.8pt}{.8pt}}}
\multiput(12.43,139.93)(-.125,-.125){3}{{\rule{.8pt}{.8pt}}}
\multiput(16.93,135.93)(-.79167,-.45833){7}{{\rule{.8pt}{.8pt}}}
\multiput(16.93,128.93)(-.75,-.33333){7}{{\rule{.8pt}{.8pt}}}
\multiput(17.43,122.93)(-.75,-.375){3}{{\rule{.8pt}{.8pt}}}
\multiput(17.18,116.68)(-.75,-.40625){9}{{\rule{.8pt}{.8pt}}}
\multiput(98.68,172.43)(-.65625,-.5){9}{{\rule{.8pt}{.8pt}}}
\multiput(100.93,167.68)(-.6875,-.375){5}{{\rule{.8pt}{.8pt}}}
\multiput(104.43,164.18)(-.65,-.35){6}{{\rule{.8pt}{.8pt}}}
\multiput(107.68,159.93)(-.5625,-.4375){5}{{\rule{.8pt}{.8pt}}}
\multiput(112.93,156.93)(-.775,-.55){11}{{\rule{.8pt}{.8pt}}}
\multiput(112.93,150.68)(-.75,-.5){11}{{\rule{.8pt}{.8pt}}}
\multiput(112.93,144.18)(-.64286,-.46429){8}{{\rule{.8pt}{.8pt}}}
\multiput(117.68,140.93)(-.625,-.5){11}{{\rule{.8pt}{.8pt}}}
\multiput(117.68,134.18)(-.64286,-.53571){8}{{\rule{.8pt}{.8pt}}}
\multiput(117.43,127.68)(-.72917,-.58333){13}{{\rule{.8pt}{.8pt}}}
\multiput(117.43,119.68)(-.69444,-.44444){10}{{\rule{.8pt}{.8pt}}}
\multiput(17.68,173.43)(-.125,0){3}{{\rule{.8pt}{.8pt}}}
\multiput(248.5,127.5)(-.03125,-.03125){8}{\line(0,-1){.03125}}
\put(17.5,173.75){\line(1,0){.25}}
\multiput(18,147)(-.05,-.0333333){15}{\line(-1,0){.05}}
\multiput(17.25,146.5)(.0625,.03125){8}{\line(1,0){.0625}}
\put(17.625,160.5){\vector(0,-1){.07}}\multiput(17.5,173.75)(.03125,-3.3125){8}{\line(0,-1){3.3125}}
\put(17.375,173.625){\vector(-1,1){.07}}\multiput(17.5,173.5)(-.03125,.03125){8}{\line(0,1){.03125}}
\put(61.5,146.5){\vector(1,0){.07}}\multiput(17.75,146.75)(5.8333333,-.0333333){15}{\line(1,0){5.8333333}}
\put(105,146.375){\vector(0,-1){.07}}\put(105,146.5){\line(0,-1){.25}}
\put(104.875,143.875){\vector(0,-1){.07}}\multiput(104.75,146.5)(.03125,-.65625){8}{\line(0,-1){.65625}}
\put(106.875,140.875){\vector(1,0){.07}}\multiput(105,141)(.46875,-.03125){8}{\line(1,0){.46875}}
\put(108.5,137.75){\vector(0,-1){.07}}\put(108.5,140.5){\line(0,-1){5.5}}
\put(110.75,135.375){\vector(1,0){.07}}\multiput(108.75,135.5)(.5,-.03125){8}{\line(1,0){.5}}
\put(112.625,130.625){\vector(0,-1){.07}}\multiput(112.5,135.25)(.03125,-1.15625){8}{\line(0,-1){1.15625}}
\put(110.5,125.875){\vector(-1,0){.07}}\multiput(112.75,125.75)(-.5625,.03125){8}{\line(-1,0){.5625}}
\put(108.25,123.25){\vector(0,-1){.07}}\put(108.25,125.75){\line(0,-1){5}}
\put(110.125,117){\vector(1,-2){.07}}\multiput(108.25,120.25)(.033482143,-.058035714){112}{\line(0,-1){.058035714}}
\end{picture}

\end{figure}

Assume that there is an $a$-band $\cal A$ starting with an $a$-edge of $t_2$ and ending
with an $a$-edge $e$ of $\partial \cal C'.$ Then $e$ belongs to some path
$ye$ where $y$ consists of $a$-edges and connects $\cal C$ and $\cal C'.$
Notice that every maximal $a$-band crossing the path $y$ must cross $t_2$
because it cannot cross $\cal A,$ and $\cal C$ has no $a$-edges.
Hence $|y|_a \le |t_2|_a -1,$ contrary to the minimality in the choice of $t_2.$

Thus every maximal  $a$-band $\cal A$ crossing $t_2$ must connect the top and the bottom
of the trapezium $\Delta$, and therefore the path $t$ must cross every such  $a$-band $\cal A.$
Also $t$ must cross every maximal $a$-band starting on $t_3$ whence $|t|_a\ge |t_2|_a+|t_3|_a = |t'|_a$.
Since the path $|t|$ must cross every maximal $\theta$-band of $\Delta$ we also have the
inequality $|t|_{\theta}\ge |t_1|_{\theta}+|t_3|_{\theta}=|t'|_{\theta}.$ Now it follows
from the definition of path length that $|t'|\le |t|$ as required.
\endproof

\begin{lemma} \label{replica} Let $\Delta$ be an accepted $i$-sector  with index $i\ne 1$
bounded by a maximal $k_i$-band ${\cal C}$ and a $k_{i+1}$-band ${\cal C}'.$ 
Let  $o_1$ and $o_2$ be two vertices lying on ${\cal C}$ and ${\cal C}',$ respectively. 
Assume that $o_1$ and $o_2$ are connected by a path $t,$ in $\Delta.$ Then the replicas
$o'_1$ and $o'_2$ of the vertices $o_1$ and  $o_2$ in the replica $\Delta'$ of $\Delta$ can be
connected by a path $t'$ such that $|t'|\le |t|.$
\end{lemma}

\proof First of all, one may assume that no maximal $q$-band 
${\cal C= C}_1,{\cal C}_2,\dots, {\cal C}_{k+2}={\cal C}'$ is crossed by the path $t$ twice.
Indeed, otherwise $t$ has a subpath $s$ of the form $ezf,$ where $e$ and $f$ are $q$-edges
of some ${\cal C}_j$ separated in this band by $m$  $(\theta,q)$-cells for some $m\ge 0$.
Therefore the path $z$ must cross at least $m$ maximal $\theta$-bands whence $|ezf|\ge m+2.$
But the vertices $e_-$ and $f_+$ can be connected along ${\cal C}_j$ by a path of length $m$
(see the example after the definition of length $|*|$), and so the path $t$ can be shortened.  

Thus the path $t$ is a product $t=t_1\dots t_{k+1},$ where each $t_j$ connects a vertex $o(j)$ lying
on ${\cal C}_j$ with a vertex $o(j+1)$ lying on ${\cal C}_{j+1},$ and for every $j=1,\dots, k$, either
$t_{j+1}$ starts or $t_j$ ends with a $q$-edge, and so $|t|=\sum_{j=1}^{k+1}|t_j|.$ As in the previous paragraph, we have that each of the $t_j$-s crosses every $\theta$-band at most
once. (Consider $ezf$, where $e$ and $f$ are $\theta$-edges of the same $\theta$-band.) Now using the
notation of Remark \ref{subsector}, it suffices to consider the replica $\Gamma'_j$ of the subsector 
$\Gamma_j$ and to find a path $t'_j$ connecting the replicas $o'(j)$ and $o'(j+1),$ with $|t'_j|\le |t_j|.$ 

We may assume that $i$ is odd. (If $i$ is even one should use a mirror argument.) 

We first consider
the path $t_2$ crossing the input subsector $\Gamma_2,$ assuming that $t_2$ has
no $q$-edges, since the $q$-edges (if any) can be attributed to the subpaths $t_1$ and $t_3$. 
By property (6) of Lemma
\ref{MtocalS}, ${\cal C}_2$ has no $a$-edges. Hence, by Lemma
\ref{prohod} applied to a subtrapezium of $\Gamma_2$ containing $t_2$, we may assume that $t_2=s_1s_2s_3,$ where $s_1$ and $s_3$ are the subpaths of the top or
bottom paths of $q$-bands ${\cal C}_2$ and ${\cal C}_3$, respectively, and $s_2$ goes along the top
or bottom of a maximal $\theta$-band $\cal T$. For  both paths $s_1$ and $s_3,$ we 
have paths $s'_1$ and $s'_3$ of the same length lying on
the boundaries of the $q$-bands  ${\cal C}'_2$ and ${\cal C}'_3$ of the replica $\Gamma'_2$ and connecting
the replicas of the vertices $(s_1)_{\pm}$ and $(s_3)_{\pm },$ respectively. The vertices $(s'_1)_+$
and  $(s'_3)_-$ are either connected by a copy $s'_2$ of $s_2$ (if the $\theta$-band $\cal T$ was copied
when we constructed the replica $\Delta'$) or $(s'_1)_+ =(s'_3)_-$ (if the corresponding $\theta$-band
of $\Delta'$ has no $a$-edges). It follows that in any case we have $|t'_2|\le |t_2|$
for $t'_2 = s'_1s'_2s'_3.$ 

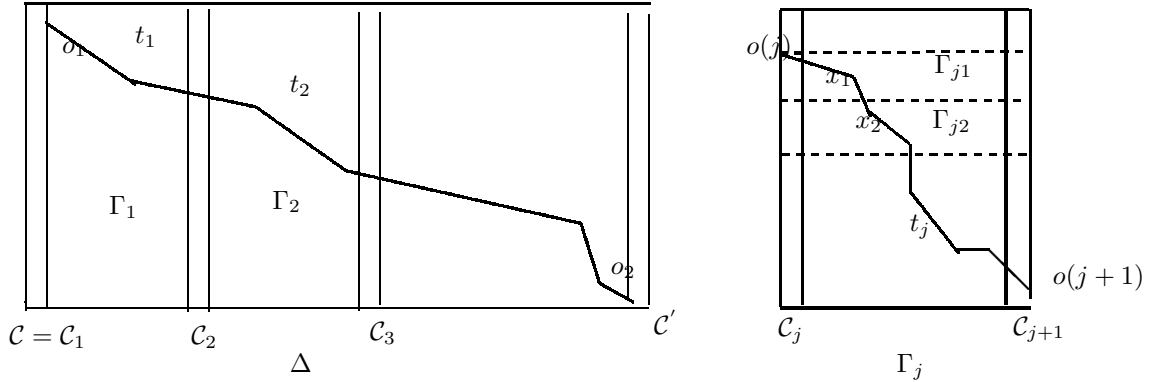
\begin{figure}[h!]
\unitlength 1mm 
\linethickness{0.4pt}
\ifx\plotpoint\undefined\newsavebox{\plotpoint}\fi 
\begin{picture}(133.75,60)(0,110)
\put(-2.75,164.5){\line(0,-1){40.75}}
\put(0,164.25){\line(0,-1){40.25}}
\put(18.75,163.75){\line(0,-1){40.25}}
\put(21.5,163.75){\line(0,-1){40.25}}
\put(41.5,163.25){\line(0,-1){39.5}}
\put(44.25,163.25){\line(0,-1){38.75}}
\put(77.25,163.25){\line(0,-1){38}}
\put(80,163){\line(0,-1){38.5}}
\put(-3,164.5){\line(1,0){83.25}}
\multiput(-2.75,125)(-.03125,-.03125){8}{\line(0,-1){.03125}}
\put(-2.75,124){\line(1,0){82.75}}
\thicklines
\multiput(-.25,162)(.0489795918,-.0336734694){245}{\line(1,0){.0489795918}}
\multiput(11,154.25)(.161057692,-.033653846){104}{\line(1,0){.161057692}}
\multiput(27.75,150.75)(.0476190476,-.0337301587){252}{\line(1,0){.0476190476}}
\multiput(39.75,142.25)(.150240385,-.033653846){208}{\line(1,0){.150240385}}
\multiput(71,135.25)(.03333333,-.10666667){75}{\line(0,-1){.10666667}}
\multiput(73.5,127.25)(.06,-.03333333){75}{\line(1,0){.06}}
\put(97.5,163.5){\line(0,-1){39.25}}
\put(100.5,163.5){\line(0,-1){39}}
\put(127.5,163.5){\line(0,-1){39}}
\put(130.75,163.75){\line(0,-1){38.5}}
\put(97.5,124){\line(1,0){33.25}}
\put(97.5,163.75){\line(1,0){33}}
\put(97.43,157.93){\line(1,0){.9779}}
\put(99.386,157.944){\line(1,0){.9779}}
\put(101.341,157.959){\line(1,0){.9779}}
\put(103.297,157.974){\line(1,0){.9779}}
\put(105.253,157.989){\line(1,0){.9779}}
\put(107.209,158.003){\line(1,0){.9779}}
\put(109.165,158.018){\line(1,0){.9779}}
\put(111.121,158.033){\line(1,0){.9779}}
\put(113.077,158.047){\line(1,0){.9779}}
\put(115.033,158.062){\line(1,0){.9779}}
\put(116.989,158.077){\line(1,0){.9779}}
\put(118.944,158.091){\line(1,0){.9779}}
\put(120.9,158.106){\line(1,0){.9779}}
\put(122.856,158.121){\line(1,0){.9779}}
\put(124.812,158.136){\line(1,0){.9779}}
\put(126.768,158.15){\line(1,0){.9779}}
\put(128.724,158.165){\line(1,0){.9779}}
\put(97.43,151.68){\line(1,0){.9706}}
\put(99.371,151.68){\line(1,0){.9706}}
\put(101.312,151.68){\line(1,0){.9706}}
\put(103.253,151.68){\line(1,0){.9706}}
\put(105.194,151.68){\line(1,0){.9706}}
\put(107.136,151.68){\line(1,0){.9706}}
\put(109.077,151.68){\line(1,0){.9706}}
\put(111.018,151.68){\line(1,0){.9706}}
\put(112.959,151.68){\line(1,0){.9706}}
\put(114.9,151.68){\line(1,0){.9706}}
\put(116.841,151.68){\line(1,0){.9706}}
\put(118.783,151.68){\line(1,0){.9706}}
\put(120.724,151.68){\line(1,0){.9706}}
\put(122.665,151.68){\line(1,0){.9706}}
\put(124.606,151.68){\line(1,0){.9706}}
\put(126.547,151.68){\line(1,0){.9706}}
\put(128.489,151.68){\line(1,0){.9706}}
\put(97.68,144.43){\line(1,0){.9848}}
\put(99.649,144.43){\line(1,0){.9848}}
\put(101.619,144.43){\line(1,0){.9848}}
\put(103.589,144.43){\line(1,0){.9848}}
\put(105.558,144.43){\line(1,0){.9848}}
\put(107.528,144.43){\line(1,0){.9848}}
\put(109.498,144.43){\line(1,0){.9848}}
\put(111.468,144.43){\line(1,0){.9848}}
\put(113.437,144.43){\line(1,0){.9848}}
\put(115.407,144.43){\line(1,0){.9848}}
\put(117.377,144.43){\line(1,0){.9848}}
\put(119.346,144.43){\line(1,0){.9848}}
\put(121.316,144.43){\line(1,0){.9848}}
\put(123.286,144.43){\line(1,0){.9848}}
\put(125.255,144.43){\line(1,0){.9848}}
\put(127.225,144.43){\line(1,0){.9848}}
\put(129.195,144.43){\line(1,0){.9848}}
\multiput(97.5,157.75)(.10955056,-.03370787){89}{\line(1,0){.10955056}}
\multiput(107.25,154.75)(.03358209,-.07835821){67}{\line(0,-1){.07835821}}
\multiput(109.5,150)(.041666667,-.033730159){126}{\line(1,0){.041666667}}
\put(114.75,145.75){\line(0,-1){6.25}}
\multiput(114.75,139.5)(.033678756,-.042746114){193}{\line(0,-1){.042746114}}
\put(121,131.75){\line(1,0){4.25}}
\put(125.25,131.75){\line(1,-1){5.25}}
\put(118,155.25){$\Gamma_{j1}$}
\put(117.75,148){$\Gamma_{j2}$}
\put(114.75,134.25){$t_j$}
\put(97.25,120.25){${\cal C}_j$}
\put(93,157.75){$o(j)$}
\put(133.75,127.25){$o(j+1)$}
\put(32.25,115.5){$\Delta$}
\put(113,115.5){$\Gamma_{j}$}
\put(128.5,120.5){${\cal C}_{j+1}$}
\put(103.5,154){$x_1$}
\put(107.5,148){$x_2$}
\put(-5,119.25){${\cal C}={\cal C}_1$}
\put(19.25,119.25){${\cal C}_2$}
\put(43,119.75){${\cal C}_3$}
\put(8.25,136.5){$\Gamma_1$}
\put(30,137.25){$\Gamma_2$}
\put(80.75,120.5){${\cal C}^{'}$}
\put(11.75,159.25){$t_1$}
\put(32.25,152.75){$t_2$}
\put(2,158){$o_1$}
\put(75,128.5){$o_2$}
\end{picture}

\end{figure}

Assume now that $j\ne 2$. The subsector $\Gamma_j$ has no $H$-cells by Lemma \ref{calAi}, and 
so it is a union of alternating subtrapezia $\Gamma_{j1}, \Gamma_{j2},\dots$ 
whose histories are words either in the alphabet $\Theta$ or in $\hat\Theta.$ Let $t_j=x_1\dots x_d,$
where every $x_k$ belongs to  $\Gamma_{jk}.$ If the history of the trapezia $\Gamma_{jk}$
is a word over $\Theta,$ then we have the copy $\Gamma'_{jk}$ of $\Gamma_{jk}$ in $\Delta'$,
and a subpath $x=x_k$ of $t_j$ lying in $\Gamma_{jk}$ has a copy $x'$ in $\Gamma'_{jk}$. If the history
is a word over $\hat\Theta,$ then $\Gamma'_{jk}$ has no $a$-edges, and for every subpath $x$ of
$t_j$ lying in $\Gamma_{jk}$, we can construct a corresponding subpath $x'$ in $\Gamma'_{jk}$
which copies only $q$- and $\theta$-edges of $x$, but ignores the $a$-edges of $x$. Since
there are no $a$-edges in the common boundaries of the neighbors $\Gamma'_{jk}$ and $\Gamma'_{j,k+1}$,
we have $(x'_k)_+=(x'_{k+1})_-$ for every $k=1,\dots,d-1$, and we obtain $|t'_2|\le |t_2|$ for the
path $t'_2=x'_1\dots x'_d.$ 

Now the required path $t'$ is obtained, and the lemma is proved.
\endproof

 Assume that we have a minimal diagram $\Gamma$ with a cyclically reduced boundary label
 over the group $G({\cal S}\cup \hat{\cal S}, L)$, 
 which is separated by a maximal $k_i$-band (or $k_{i+1}$-band) $\cal C$ in two parts
 $\tilde \Delta$ and $\Delta$ such that $\Delta$ is an accepted $i$-sector. Assume that
 every maximal $\theta$-band ${\cal T}_1,\dots, {\cal T}_m$ of $\Gamma$ crosses $\Delta,
$ the bottom $x$ of ${\cal T}_1$ is a part of the boundary $\partial \Gamma$, and $\Lab(x)$ is a suffix
of the subword $k_{i-1}\dots k_i\dots k_{i+1}$(a prefix of the subword 
$k_{i}\dots k_{i+1}\dots k_{i+2},$ respectively)
of the word $\Sigma_0.$  Also we assume that for every $j\in [2,m]$, the trimmed bottom $\tbott{\cal T}_j$
is a subpath of the top $\topp{\cal T}_{j-1}. $ Below we call such a diagram {\it unfinished}
  if $i,i-1 \ne 1$ (if $i, i+1\ne 1$, respectively).

\begin{lemma} \label{2sectors} The part $\tilde\Delta$ of the unfinished diagram $\Gamma$ 
can be embedded in an $i-1$-sector (respectively, $i+1$-sector) $\nabla$ which is a mirror 
copy of the $i$-sector $\Delta.$
\end{lemma}

\proof Without loss of generality we assume that $\cal C$ is a $k_i$-band in the definition 
of unfinished diagram $\Gamma$. We construct $\nabla$ by induction on the length $m$
of the history 
$h\equiv \theta_1\dots \theta_m$ 
of the $i$-sector $\Delta.$ 

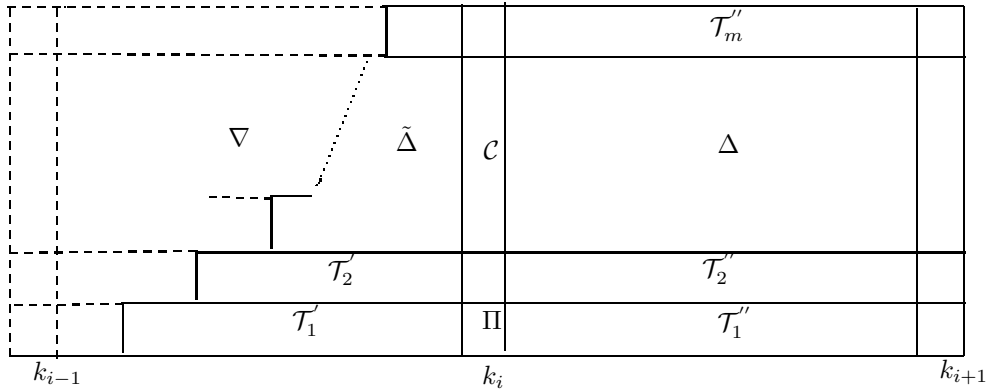
\begin{figure}[h!]
\unitlength 1mm 
\linethickness{0.4pt}
\ifx\plotpoint\undefined\newsavebox{\plotpoint}\fi 
\begin{picture}(139.25,60)(0,120)
\put(12,133){\line(1,0){126.75}}
\put(138.75,133){\line(0,1){46.5}}
\put(72,133.25){\line(0,1){46.5}}
\put(77.75,133.75){\line(0,1){46.25}}
\put(138.75,179.5){\line(-1,0){77.25}}
\put(132.5,179.25){\line(0,-1){46.25}}
\put(61.75,172.75){\line(1,0){77.25}}
\put(62,179.5){\line(0,-1){6.5}}
\put(26.75,140){\line(1,0){112.25}}
\put(27,140){\line(0,-1){6.5}}
\put(139.25,146.75){\line(0,1){0}}
\put(36.75,146.75){\line(1,0){102.25}}
\put(36.75,147){\line(0,-1){6.5}}
\put(46.75,154.25){\line(0,-1){7}}
\put(47,154.25){\line(1,0){5}}
\multiput(53.18,156.43)(0,-.125){3}{{\rule{.4pt}{.4pt}}}
\multiput(52.68,155.68)(.03125,.03125){4}{\line(0,1){.03125}}
\multiput(52.93,155.93)(.35526,.89474){20}{{\rule{.4pt}{.4pt}}}
\put(62.18,179.43){\line(-1,0){.9902}}
\put(60.199,179.43){\line(-1,0){.9902}}
\put(58.219,179.43){\line(-1,0){.9902}}
\put(56.239,179.43){\line(-1,0){.9902}}
\put(54.258,179.43){\line(-1,0){.9902}}
\put(52.278,179.43){\line(-1,0){.9902}}
\put(50.297,179.43){\line(-1,0){.9902}}
\put(48.317,179.43){\line(-1,0){.9902}}
\put(46.337,179.43){\line(-1,0){.9902}}
\put(44.356,179.43){\line(-1,0){.9902}}
\put(42.376,179.43){\line(-1,0){.9902}}
\put(40.395,179.43){\line(-1,0){.9902}}
\put(38.415,179.43){\line(-1,0){.9902}}
\put(36.435,179.43){\line(-1,0){.9902}}
\put(34.454,179.43){\line(-1,0){.9902}}
\put(32.474,179.43){\line(-1,0){.9902}}
\put(30.493,179.43){\line(-1,0){.9902}}
\put(28.513,179.43){\line(-1,0){.9902}}
\put(26.533,179.43){\line(-1,0){.9902}}
\put(24.552,179.43){\line(-1,0){.9902}}
\put(22.572,179.43){\line(-1,0){.9902}}
\put(20.591,179.43){\line(-1,0){.9902}}
\put(18.611,179.43){\line(-1,0){.9902}}
\put(16.631,179.43){\line(-1,0){.9902}}
\put(14.65,179.43){\line(-1,0){.9902}}
\put(12.67,179.43){\line(-1,0){.9902}}
\put(11.93,179.18){\line(0,-1){.984}}
\put(11.93,177.212){\line(0,-1){.984}}
\put(11.93,175.244){\line(0,-1){.984}}
\put(11.93,173.275){\line(0,-1){.984}}
\put(11.93,171.307){\line(0,-1){.984}}
\put(11.93,169.339){\line(0,-1){.984}}
\put(11.93,167.371){\line(0,-1){.984}}
\put(11.93,165.403){\line(0,-1){.984}}
\put(11.93,163.435){\line(0,-1){.984}}
\put(11.93,161.467){\line(0,-1){.984}}
\put(11.93,159.499){\line(0,-1){.984}}
\put(11.93,157.531){\line(0,-1){.984}}
\put(11.93,155.563){\line(0,-1){.984}}
\put(11.93,153.595){\line(0,-1){.984}}
\put(11.93,151.627){\line(0,-1){.984}}
\put(11.93,149.658){\line(0,-1){.984}}
\put(11.93,147.69){\line(0,-1){.984}}
\put(11.93,145.722){\line(0,-1){.984}}
\put(11.93,143.754){\line(0,-1){.984}}
\put(11.93,141.786){\line(0,-1){.984}}
\put(11.93,139.818){\line(0,-1){.984}}
\put(11.93,137.85){\line(0,-1){.984}}
\put(11.93,135.882){\line(0,-1){.984}}
\put(11.93,133.914){\line(0,-1){.984}}
\put(19.18,179.43){\line(1,0){.125}}
\put(18.18,179.18){\line(0,-1){.9894}}
\put(18.18,177.201){\line(0,-1){.9894}}
\put(18.18,175.222){\line(0,-1){.9894}}
\put(18.18,173.244){\line(0,-1){.9894}}
\put(18.18,171.265){\line(0,-1){.9894}}
\put(18.18,169.286){\line(0,-1){.9894}}
\put(18.18,167.307){\line(0,-1){.9894}}
\put(18.18,165.329){\line(0,-1){.9894}}
\put(18.18,163.35){\line(0,-1){.9894}}
\put(18.18,161.371){\line(0,-1){.9894}}
\put(18.18,159.392){\line(0,-1){.9894}}
\put(18.18,157.414){\line(0,-1){.9894}}
\put(18.18,155.435){\line(0,-1){.9894}}
\put(18.18,153.456){\line(0,-1){.9894}}
\put(18.18,151.478){\line(0,-1){.9894}}
\put(18.18,149.499){\line(0,-1){.9894}}
\put(18.18,147.52){\line(0,-1){.9894}}
\put(18.18,145.541){\line(0,-1){.9894}}
\put(18.18,143.563){\line(0,-1){.9894}}
\put(18.18,141.584){\line(0,-1){.9894}}
\put(18.18,139.605){\line(0,-1){.9894}}
\put(18.18,137.627){\line(0,-1){.9894}}
\put(18.18,135.648){\line(0,-1){.9894}}
\put(18.18,133.669){\line(0,-1){.9894}}
\put(61.93,172.93){\line(-1,0){.9804}}
\put(59.969,172.94){\line(-1,0){.9804}}
\put(58.008,172.949){\line(-1,0){.9804}}
\put(56.047,172.959){\line(-1,0){.9804}}
\put(54.087,172.969){\line(-1,0){.9804}}
\put(52.126,172.979){\line(-1,0){.9804}}
\put(50.165,172.989){\line(-1,0){.9804}}
\put(48.204,172.998){\line(-1,0){.9804}}
\put(46.243,173.008){\line(-1,0){.9804}}
\put(44.283,173.018){\line(-1,0){.9804}}
\put(42.322,173.028){\line(-1,0){.9804}}
\put(40.361,173.038){\line(-1,0){.9804}}
\put(38.4,173.047){\line(-1,0){.9804}}
\put(36.44,173.057){\line(-1,0){.9804}}
\put(34.479,173.067){\line(-1,0){.9804}}
\put(32.518,173.077){\line(-1,0){.9804}}
\put(30.557,173.087){\line(-1,0){.9804}}
\put(28.596,173.096){\line(-1,0){.9804}}
\put(26.636,173.106){\line(-1,0){.9804}}
\put(24.675,173.116){\line(-1,0){.9804}}
\put(22.714,173.126){\line(-1,0){.9804}}
\put(20.753,173.136){\line(-1,0){.9804}}
\put(18.792,173.145){\line(-1,0){.9804}}
\put(16.832,173.155){\line(-1,0){.9804}}
\put(14.871,173.165){\line(-1,0){.9804}}
\put(12.91,173.175){\line(-1,0){.9804}}
\put(26.93,139.93){\line(-1,0){.9531}}
\put(25.023,139.898){\line(-1,0){.9531}}
\put(23.117,139.867){\line(-1,0){.9531}}
\put(21.211,139.836){\line(-1,0){.9531}}
\put(19.305,139.805){\line(-1,0){.9531}}
\put(17.398,139.773){\line(-1,0){.9531}}
\put(15.492,139.742){\line(-1,0){.9531}}
\put(13.586,139.711){\line(-1,0){.9531}}
\put(36.68,146.93){\line(-1,0){.99}}
\put(34.7,146.91){\line(-1,0){.99}}
\put(32.72,146.89){\line(-1,0){.99}}
\put(30.74,146.87){\line(-1,0){.99}}
\put(28.76,146.85){\line(-1,0){.99}}
\put(26.78,146.83){\line(-1,0){.99}}
\put(24.8,146.81){\line(-1,0){.99}}
\put(22.82,146.79){\line(-1,0){.99}}
\put(20.84,146.77){\line(-1,0){.99}}
\put(18.86,146.75){\line(-1,0){.99}}
\put(16.88,146.73){\line(-1,0){.99}}
\put(14.9,146.71){\line(-1,0){.99}}
\put(12.92,146.69){\line(-1,0){.99}}
\put(46.68,153.93){\line(-1,0){.9167}}
\put(44.846,153.985){\line(-1,0){.9167}}
\put(43.013,154.041){\line(-1,0){.9167}}
\put(41.18,154.096){\line(-1,0){.9167}}
\put(39.346,154.152){\line(-1,0){.9167}}
\put(75,159.25){$\cal C$}
\put(106,160){$\Delta$}
\put(74.75,136.5){$\Pi$}
\put(74.75,129.25){$k_i$}
\put(15,129.75){$k_{i-1}$}
\put(135.5,130){$k_{i+1}$}
\put(49.5,136.5){${\cal T}_1^{'}$}
\put(54.25,143.25){${\cal T}_2^{'}$}
\put(106,136.25){${\cal T}_1^{''}$}
\put(104,143.5){${\cal T}_2^{''}$}
\put(105,175.75){${\cal T}_m^{''}$}
\put(63.25,160.25){$\tilde\Delta$}
\put(41,161){$\nabla$}
\end{picture}

\end{figure}

 Let $\Pi$  be the intersection cell for ${\cal T}_1$ and $\cal C.$  Then ${\cal T}_1$  
 consists of $\Pi$ and two subbands ${\cal T}'_1$ and ${\cal T}''_1,$  where ${\cal T}'_1$
 belongs to $\tilde\Delta$ and ${\cal T}''_1$ belongs to $\Delta.$ Let $\Pi(1)'$ and $\Pi(1)''$
 be the neighboring cells for $\Pi$  in ${\cal T}'_1$ and in ${\cal T}''_1,$ respectively. Then
 $\Pi(1)'$ (if it exists) is a mirror copy of $\Pi(1)''$ with boundary label in the alphabet ${\cal A}_{i-1}.$
 Indeed, these two cells are determined by the same history $\theta_1$ and mirror $q$-letters of the word 
 $\Sigma_0$ as this follows from Relations (\ref{rel1}, \ref{rel2}). A similar argument shows that if $\Pi'(1)$ has a neighbor $\Pi'(2)$ (where $\Pi'(2)\ne \Pi$)
 in ${\cal T}'_1$,  then the cell $\Pi''(1)$ has a neighbor $\Pi''(2)$ (where $\Pi''(2)\ne \Pi$)
 in ${\cal T}''_1,$ and $\Pi'(2)$ is a mirror copy of $\Pi''(2)$ with boundary label in the alphabet
 ${\cal A}_{i-1}.$  By induction, we obtain that ${\cal T}'_1$ is a mirror copy of a subband of 
 ${\cal T}''_1$ starting with $\Pi.$ Hence one can extend ${\cal T}'_1$ (and the subdiagram $\tilde\Delta$)
 and obtain a $\theta$-band ${\cal T}^{\nabla}_1$ which is a mirror copy of ${\cal T}''_1.$

 Since the trimmed bottom of the band ${\cal T}_2$ is a subpath of the top of ${\cal T}_1,$ one can similarly
 subdivide  ${\cal T}_2$ in $\Pi(2),$ ${\cal T}'_2$, ${\cal T}''_2,$ and prove that ${\cal T}'_2$
 is the mirror copy of a subband of ${\cal T}''_2,$ starting with the cell $\Pi(2)$ and having
 boundary label over the alphabet ${\cal A}_{i-1}.$ (${\cal T'}_2$ and ${\cal T''}_2$ may
 include $(\theta,q)$-cells and also $(\theta,a)$-cells.) Therefore there is an extension ${\cal T}^{\nabla}_2$
 of ${\cal T}'_2$, and this extension is a mirror copy of ${\cal T}''_2.$ Then by induction, we
 construct $\theta$-bands ${\cal T}^{\nabla}_3,\dots, {\cal T}^{\nabla}_m,$ and these $\theta$-bands
 together with the $q$-band $\cal C$ form the required $i-1$-sector $\nabla.$ \endproof

\begin{lemma} \label{sravn}  Let $\Delta$ be a minimal diagram with boundary path $x_1y_1x_2y_2$, where

(1) $\Lab(y_1)$ is a subword $k_j\dots k_i$ of $\Sigma_0$, and this subword
does not contain the letter $k_1;$

(2) $x_1$ and $x_2$ are sides of the maximal $k_j$-band ${\cal K}_j$  and $k_i$-band ${\cal K}_i$
starting on $y_1$;

(3) every cell of $\Delta$ belongs to one of the maximal $\theta$-bands ${\cal T}_1,\dots,{\cal T}_m$ of $\Delta$;

(4) (a) either each of the bands ${\cal T}_1,\dots,{\cal T}_m$ crosses ${\cal K}_j$ or (b) each of them crosses ${\cal K}_i$;

(5) the trimmed bottom path of ${\cal T}_1$ is a subpath of $y_1^{\pm 1}$, and the trimmed bottom path of ${\cal T}_l$ is
a subpath of the top path of ${\cal T}_{l-1}$ for every $l=2,\dots, m;$ 

(6) one can construct a diagram $\bar \Delta$ with boundary $\bar x_1\bar y_1\bar x_2\bar y_2,$
and $\bar\Delta$ satisfies the analogs of properties (1)-(5), but $\Lab(\bar y_1)\equiv k_{j-1}\dots k_i$ in
case 4(a) ($\Lab(\bar y_1)\equiv k_{j}\dots k_{i+1}$ in case 4(b)), and $\Delta$ is embeddable in
$\bar\Delta$ so that the $k_j$-band ${\cal K}_j$  and the $k_i$-band ${\cal K}_i$ remain maximal in $\bar\Delta.$

Then in case (4)(a) (in case 4(b)), there exists a diagram $\Delta'$ over the group $G_1$ with boundary path $x'_1y'_1x'_2y'_2$ such that $\Lab((y'_1)^{-1})$ is the subword $k_{2j-i}\dots k_{j}$ (respectively, 
$k_{i}\dots k_{2i-j}$) of $\Sigma_0$, $x'_1$ and $x'_2$ are sides of the maximal  $k_{j}$-band ${\cal K}_{j}$ 
and $k_{2j-i}$-band ${\cal K}_{2j-i}$ (of the maximal  $k_{i}$-band ${\cal K}_{i}$ 
and $k_{2i-j}$-band ${\cal K}_{2i-j}$ starting on $y'_1$), the labels of $x'_1$ and $x'_2$ are copies of $\phi(x_1)$ and $\phi(x_2)$, resp.,
and $|y'_2|\le|y_2|.$ 

(The subscripts of $k$-bands are taken modulo $L.$)
   \end{lemma}
   
   \proof   We consider the case (4)(b) only. The maximal $k$-bands ${\cal K}_j,{\cal K}_{j+1}\dots {\cal K}_{i+1}$
   subdivide the diagram $\bar\Delta$ in the subdiagrams $\Gamma_l$-s, where $\Gamma_l$ is bounded by ${\cal K}_l$ and
   ${\cal K}_{l+1}$ and $\Gamma_l$ includes these $k$-bands ($l=j,\dots, i$).

Every maximal $\theta$-band of $\Delta$ having a cell in $\Gamma_{i-1},$ must cross both bands ${\cal K}_i$ and ${\cal K}_{i+1}$ of $\bar\Delta.$
Therefore the  parts of these bands in $\Gamma_{i-1}$ and $\Gamma_i$ form an unfinished diagram whose $i$-sector is a subdiagram $\Delta_i$ of $\Gamma_i$.
By Lemma \ref{2sectors}, the subdiagram $\Gamma_{i-1}$ is embedded into the mirror copy $\Delta_{i-1}$ of 
$\Delta_i$. Similarly, $\Gamma_{i-2}$ is embeddable into a mirror copy of $\Delta_{i-1}$ which is the copy of $\Delta_i$ (we denote this copy by $\Delta_{i-2}$)
,..., $\Gamma_{j}$ is embeddable into the copy (or mirror copy) $\Delta_j$ of $\Delta_i.$

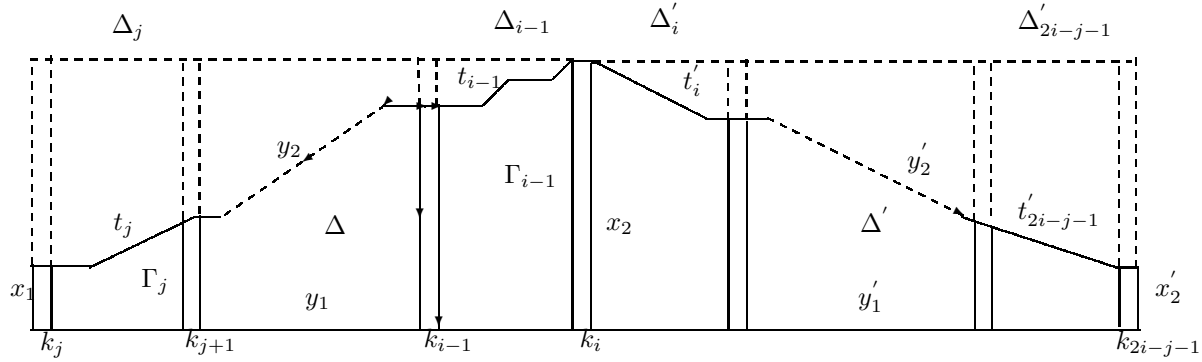
\begin{figure}[h!]
\unitlength 1mm 
\linethickness{0.4pt}
\ifx\plotpoint\undefined\newsavebox{\plotpoint}\fi 
\begin{picture}(162.5,65)(0,100)
\put(13,118.75){\line(1,0){147.5}}
\put(85,119){\line(0,1){35.5}}
\put(87.5,154.25){\line(0,-1){35.25}}
\put(87.75,154.5){\line(-1,0){2.75}}
\multiput(85,154.75)(-.03353659,-.03353659){82}{\line(0,-1){.03353659}}
\put(82.25,152){\line(-1,0){5.75}}
\multiput(76.5,152)(-.033653846,-.033653846){104}{\line(0,-1){.033653846}}
\put(73,148.5){\line(-1,0){13}}
\put(13,127.25){\line(1,0){8.25}}
\multiput(20.75,127)(.069651741,.03358209){201}{\line(1,0){.069651741}}
\put(34.75,133.75){\line(1,0){3.5}}
\put(59.875,148.625){\vector(-1,-1){.07}}\multiput(60,148.75)(-.03125,-.03125){8}{\line(0,-1){.03125}}
\put(49.125,141.125){\vector(-3,-2){.07}}\multiput(60.18,148.68)(-.0467437,-.0320378){17}{\line(-1,0){.0467437}}
\multiput(58.59,147.59)(-.0467437,-.0320378){17}{\line(-1,0){.0467437}}
\multiput(57.001,146.501)(-.0467437,-.0320378){17}{\line(-1,0){.0467437}}
\multiput(55.412,145.412)(-.0467437,-.0320378){17}{\line(-1,0){.0467437}}
\multiput(53.823,144.323)(-.0467437,-.0320378){17}{\line(-1,0){.0467437}}
\multiput(52.233,143.233)(-.0467437,-.0320378){17}{\line(-1,0){.0467437}}
\multiput(50.644,142.144)(-.0467437,-.0320378){17}{\line(-1,0){.0467437}}
\multiput(49.055,141.055)(-.0467437,-.0320378){17}{\line(-1,0){.0467437}}
\multiput(47.465,139.965)(-.0467437,-.0320378){17}{\line(-1,0){.0467437}}
\multiput(45.876,138.876)(-.0467437,-.0320378){17}{\line(-1,0){.0467437}}
\multiput(44.287,137.787)(-.0467437,-.0320378){17}{\line(-1,0){.0467437}}
\multiput(42.698,136.698)(-.0467437,-.0320378){17}{\line(-1,0){.0467437}}
\multiput(41.108,135.608)(-.0467437,-.0320378){17}{\line(-1,0){.0467437}}
\multiput(39.519,134.519)(-.0467437,-.0320378){17}{\line(-1,0){.0467437}}
\put(65.625,148.5){\vector(1,0){.07}}\put(65.43,148.43){\line(1,0){.125}}
\put(64.75,133.625){\vector(0,-1){.07}}\put(64.75,148.25){\line(0,-1){29.25}}
\put(67.625,148.5){\vector(1,0){.07}}\put(67.5,148.5){\line(1,0){.25}}
\put(67.25,119.125){\vector(0,-1){.07}}\put(67.25,119.25){\line(0,-1){.25}}
\put(67.25,119){\line(0,1){29.75}}
\put(13.25,127.5){\line(0,-1){8.5}}
\put(15.75,127.25){\line(0,-1){8.5}}
\put(33.25,133.25){\line(0,-1){14.5}}
\put(35.5,134){\line(0,-1){15.25}}
\put(24,132){$t_j$}
\put(85.18,154.68){\line(-1,0){.9897}}
\put(83.2,154.68){\line(-1,0){.9897}}
\put(81.221,154.68){\line(-1,0){.9897}}
\put(79.241,154.68){\line(-1,0){.9897}}
\put(77.262,154.68){\line(-1,0){.9897}}
\put(75.282,154.68){\line(-1,0){.9897}}
\put(73.303,154.68){\line(-1,0){.9897}}
\put(71.324,154.68){\line(-1,0){.9897}}
\put(69.344,154.68){\line(-1,0){.9897}}
\put(67.365,154.68){\line(-1,0){.9897}}
\put(65.385,154.68){\line(-1,0){.9897}}
\put(63.406,154.68){\line(-1,0){.9897}}
\put(61.426,154.68){\line(-1,0){.9897}}
\put(59.447,154.68){\line(-1,0){.9897}}
\put(57.467,154.68){\line(-1,0){.9897}}
\put(55.488,154.68){\line(-1,0){.9897}}
\put(53.508,154.68){\line(-1,0){.9897}}
\put(51.529,154.68){\line(-1,0){.9897}}
\put(49.55,154.68){\line(-1,0){.9897}}
\put(47.57,154.68){\line(-1,0){.9897}}
\put(45.591,154.68){\line(-1,0){.9897}}
\put(43.611,154.68){\line(-1,0){.9897}}
\put(41.632,154.68){\line(-1,0){.9897}}
\put(39.652,154.68){\line(-1,0){.9897}}
\put(37.673,154.68){\line(-1,0){.9897}}
\put(35.693,154.68){\line(-1,0){.9897}}
\put(33.714,154.68){\line(-1,0){.9897}}
\put(31.735,154.68){\line(-1,0){.9897}}
\put(29.755,154.68){\line(-1,0){.9897}}
\put(27.776,154.68){\line(-1,0){.9897}}
\put(25.796,154.68){\line(-1,0){.9897}}
\put(23.817,154.68){\line(-1,0){.9897}}
\put(21.837,154.68){\line(-1,0){.9897}}
\put(19.858,154.68){\line(-1,0){.9897}}
\put(17.878,154.68){\line(-1,0){.9897}}
\put(15.899,154.68){\line(-1,0){.9897}}
\put(13.919,154.68){\line(-1,0){.9897}}
\put(13.18,127.43){\line(0,1){.9732}}
\put(13.18,129.376){\line(0,1){.9732}}
\put(13.18,131.323){\line(0,1){.9732}}
\put(13.18,133.269){\line(0,1){.9732}}
\put(13.18,135.215){\line(0,1){.9732}}
\put(13.18,137.162){\line(0,1){.9732}}
\put(13.18,139.108){\line(0,1){.9732}}
\put(13.18,141.055){\line(0,1){.9732}}
\put(13.18,143.001){\line(0,1){.9732}}
\put(13.18,144.948){\line(0,1){.9732}}
\put(13.18,146.894){\line(0,1){.9732}}
\put(13.18,148.84){\line(0,1){.9732}}
\put(13.18,150.787){\line(0,1){.9732}}
\put(13.18,152.733){\line(0,1){.9732}}
\put(15.68,126.93){\line(0,1){.9741}}
\put(15.68,128.878){\line(0,1){.9741}}
\put(15.68,130.826){\line(0,1){.9741}}
\put(15.68,132.775){\line(0,1){.9741}}
\put(15.68,134.723){\line(0,1){.9741}}
\put(15.68,136.671){\line(0,1){.9741}}
\put(15.68,138.619){\line(0,1){.9741}}
\put(15.68,140.568){\line(0,1){.9741}}
\put(15.68,142.516){\line(0,1){.9741}}
\put(15.68,144.464){\line(0,1){.9741}}
\put(15.68,146.412){\line(0,1){.9741}}
\put(15.68,148.361){\line(0,1){.9741}}
\put(15.68,150.309){\line(0,1){.9741}}
\put(15.68,152.257){\line(0,1){.9741}}
\put(15.68,154.206){\line(0,1){.9741}}
\put(33.18,133.43){\line(0,1){.9886}}
\put(33.18,135.407){\line(0,1){.9886}}
\put(33.18,137.384){\line(0,1){.9886}}
\put(33.18,139.362){\line(0,1){.9886}}
\put(33.18,141.339){\line(0,1){.9886}}
\put(33.18,143.316){\line(0,1){.9886}}
\put(33.18,145.293){\line(0,1){.9886}}
\put(33.18,147.271){\line(0,1){.9886}}
\put(33.18,149.248){\line(0,1){.9886}}
\put(33.18,151.225){\line(0,1){.9886}}
\put(33.18,153.202){\line(0,1){.9886}}
\put(35.43,134.18){\line(0,1){.9545}}
\put(35.43,136.089){\line(0,1){.9545}}
\put(35.43,137.998){\line(0,1){.9545}}
\put(35.43,139.907){\line(0,1){.9545}}
\put(35.43,141.816){\line(0,1){.9545}}
\put(35.43,143.725){\line(0,1){.9545}}
\put(35.43,145.634){\line(0,1){.9545}}
\put(35.43,147.543){\line(0,1){.9545}}
\put(35.43,149.452){\line(0,1){.9545}}
\put(35.43,151.362){\line(0,1){.9545}}
\put(35.43,153.271){\line(0,1){.9545}}
\put(64.68,148.68){\line(0,1){.8929}}
\put(64.608,150.465){\line(0,1){.8929}}
\put(64.537,152.251){\line(0,1){.8929}}
\put(64.465,154.037){\line(0,1){.8929}}
\put(66.93,154.68){\line(0,-1){.8571}}
\put(66.93,152.965){\line(0,-1){.8571}}
\put(66.93,151.251){\line(0,-1){.8571}}
\put(66.93,149.537){\line(0,-1){.8571}}
\put(69.75,151.5){$t_{i-1}$}
\multiput(87.93,154.18)(-.09375,.03125){4}{\line(-1,0){.09375}}
\multiput(87.75,154.5)(.066304348,-.033695652){230}{\line(1,0){.066304348}}
\put(103,146.75){\line(1,0){8}}
\put(160.25,127.25){\line(0,-1){8.25}}
\put(160,127){\line(-1,0){2.5}}
\multiput(157.5,127)(-.103233831,.03358209){201}{\line(-1,0){.103233831}}
\put(137,134){\vector(2,-1){.07}}\multiput(110.93,146.93)(.0645161,-.0322581){13}{\line(1,0){.0645161}}
\multiput(112.607,146.091)(.0645161,-.0322581){13}{\line(1,0){.0645161}}
\multiput(114.285,145.252)(.0645161,-.0322581){13}{\line(1,0){.0645161}}
\multiput(115.962,144.414)(.0645161,-.0322581){13}{\line(1,0){.0645161}}
\multiput(117.639,143.575)(.0645161,-.0322581){13}{\line(1,0){.0645161}}
\multiput(119.317,142.736)(.0645161,-.0322581){13}{\line(1,0){.0645161}}
\multiput(120.994,141.897)(.0645161,-.0322581){13}{\line(1,0){.0645161}}
\multiput(122.672,141.059)(.0645161,-.0322581){13}{\line(1,0){.0645161}}
\multiput(124.349,140.22)(.0645161,-.0322581){13}{\line(1,0){.0645161}}
\multiput(126.026,139.381)(.0645161,-.0322581){13}{\line(1,0){.0645161}}
\multiput(127.704,138.543)(.0645161,-.0322581){13}{\line(1,0){.0645161}}
\multiput(129.381,137.704)(.0645161,-.0322581){13}{\line(1,0){.0645161}}
\multiput(131.059,136.865)(.0645161,-.0322581){13}{\line(1,0){.0645161}}
\multiput(132.736,136.026)(.0645161,-.0322581){13}{\line(1,0){.0645161}}
\multiput(134.414,135.188)(.0645161,-.0322581){13}{\line(1,0){.0645161}}
\multiput(136.091,134.349)(.0645161,-.0322581){13}{\line(1,0){.0645161}}
\put(105.75,146.75){\line(0,-1){28}}
\put(108.25,146.5){\line(0,-1){27.5}}
\put(157.75,127){\line(0,-1){7.75}}
\put(140.75,118.75){\line(0,1){13.75}}
\put(138.5,133){\line(0,-1){14.25}}
\put(129.5,140.25){$y_2^{'}$}
\put(45.75,142.25){$y_2$}
\put(99.75,150.75){$t_i^{'}$}
\put(10.25,123.25){$x_1$}
\put(87.68,154.43){\line(1,0){.9861}}
\put(89.652,154.43){\line(1,0){.9861}}
\put(91.624,154.43){\line(1,0){.9861}}
\put(93.596,154.43){\line(1,0){.9861}}
\put(95.569,154.43){\line(1,0){.9861}}
\put(97.541,154.43){\line(1,0){.9861}}
\put(99.513,154.43){\line(1,0){.9861}}
\put(101.485,154.43){\line(1,0){.9861}}
\put(103.457,154.43){\line(1,0){.9861}}
\put(105.43,154.43){\line(1,0){.9861}}
\put(107.402,154.43){\line(1,0){.9861}}
\put(109.374,154.43){\line(1,0){.9861}}
\put(111.346,154.43){\line(1,0){.9861}}
\put(113.319,154.43){\line(1,0){.9861}}
\put(115.291,154.43){\line(1,0){.9861}}
\put(117.263,154.43){\line(1,0){.9861}}
\put(119.235,154.43){\line(1,0){.9861}}
\put(121.207,154.43){\line(1,0){.9861}}
\put(123.18,154.43){\line(1,0){.9861}}
\put(125.152,154.43){\line(1,0){.9861}}
\put(127.124,154.43){\line(1,0){.9861}}
\put(129.096,154.43){\line(1,0){.9861}}
\put(131.069,154.43){\line(1,0){.9861}}
\put(133.041,154.43){\line(1,0){.9861}}
\put(135.013,154.43){\line(1,0){.9861}}
\put(136.985,154.43){\line(1,0){.9861}}
\put(138.957,154.43){\line(1,0){.9861}}
\put(140.93,154.43){\line(1,0){.9861}}
\put(142.902,154.43){\line(1,0){.9861}}
\put(144.874,154.43){\line(1,0){.9861}}
\put(146.846,154.43){\line(1,0){.9861}}
\put(148.819,154.43){\line(1,0){.9861}}
\put(150.791,154.43){\line(1,0){.9861}}
\put(152.763,154.43){\line(1,0){.9861}}
\put(154.735,154.43){\line(1,0){.9861}}
\put(156.707,154.43){\line(1,0){.9861}}
\put(159.93,126.93){\line(0,1){.9655}}
\put(159.93,128.861){\line(0,1){.9655}}
\put(159.93,130.792){\line(0,1){.9655}}
\put(159.93,132.723){\line(0,1){.9655}}
\put(159.93,134.654){\line(0,1){.9655}}
\put(159.93,136.585){\line(0,1){.9655}}
\put(159.93,138.516){\line(0,1){.9655}}
\put(159.93,140.447){\line(0,1){.9655}}
\put(159.93,142.378){\line(0,1){.9655}}
\put(159.93,144.309){\line(0,1){.9655}}
\put(159.93,146.24){\line(0,1){.9655}}
\put(159.93,148.171){\line(0,1){.9655}}
\put(159.93,150.102){\line(0,1){.9655}}
\put(159.93,152.033){\line(0,1){.9655}}
\put(159.93,153.964){\line(0,1){.9655}}
\put(157.68,127.18){\line(0,1){.9821}}
\put(157.68,129.144){\line(0,1){.9821}}
\put(157.68,131.108){\line(0,1){.9821}}
\put(157.68,133.073){\line(0,1){.9821}}
\put(157.68,135.037){\line(0,1){.9821}}
\put(157.68,137.001){\line(0,1){.9821}}
\put(157.68,138.965){\line(0,1){.9821}}
\put(157.68,140.93){\line(0,1){.9821}}
\put(157.68,142.894){\line(0,1){.9821}}
\put(157.68,144.858){\line(0,1){.9821}}
\put(157.68,146.823){\line(0,1){.9821}}
\put(157.68,148.787){\line(0,1){.9821}}
\put(157.68,150.751){\line(0,1){.9821}}
\put(157.68,152.715){\line(0,1){.9821}}
\put(158.68,154.43){\line(1,0){.625}}
\put(138.43,133.18){\line(0,1){.9773}}
\put(138.43,135.134){\line(0,1){.9773}}
\put(138.43,137.089){\line(0,1){.9773}}
\put(138.43,139.043){\line(0,1){.9773}}
\put(138.43,140.998){\line(0,1){.9773}}
\put(138.43,142.952){\line(0,1){.9773}}
\put(138.43,144.907){\line(0,1){.9773}}
\put(138.43,146.862){\line(0,1){.9773}}
\put(138.43,148.816){\line(0,1){.9773}}
\put(138.43,150.771){\line(0,1){.9773}}
\put(138.43,152.725){\line(0,1){.9773}}
\put(140.68,154.18){\line(0,-1){.9643}}
\put(140.68,152.251){\line(0,-1){.9643}}
\put(140.68,150.323){\line(0,-1){.9643}}
\put(140.68,148.394){\line(0,-1){.9643}}
\put(140.68,146.465){\line(0,-1){.9643}}
\put(140.68,144.537){\line(0,-1){.9643}}
\put(140.68,142.608){\line(0,-1){.9643}}
\put(140.68,140.68){\line(0,-1){.9643}}
\put(140.68,138.751){\line(0,-1){.9643}}
\put(140.68,136.823){\line(0,-1){.9643}}
\put(140.68,134.894){\line(0,-1){.9643}}
\put(105.68,146.68){\line(0,1){.9688}}
\put(105.68,148.617){\line(0,1){.9688}}
\put(105.68,150.555){\line(0,1){.9688}}
\put(105.68,152.492){\line(0,1){.9688}}
\put(108.18,146.68){\line(0,1){.8889}}
\put(108.18,148.457){\line(0,1){.8889}}
\put(108.18,150.235){\line(0,1){.8889}}
\put(108.18,152.013){\line(0,1){.8889}}
\put(108.18,153.791){\line(0,1){.8889}}
\put(49.5,122){$y_1$}
\put(14.25,116){$k_j$}
\put(86,116.5){$k_i$}
\put(157.25,116.25){$k_{2i-j-1}$}
\put(33.5,116.5){$k_{j+1}$}
\put(65.25,116.5){$k_{i-1}$}
\put(27.75,124.5){$\Gamma_j$}
\put(76,138.25){$\Gamma_{i-1}$}
\put(52,131.5){$\Delta$}
\put(123.25,131.75){$\Delta^{'}$}
\put(144.25,133.5){$t_{2i-j-1}^{'}$}
\put(23.75,158.25){$\Delta_j$}
\put(74.5,158.75){$\Delta_{i-1}$}
\put(95.25,159){$\Delta_i^{'}$}
\put(144.25,158.75){$\Delta_{2i-j-1}^{'}$}
\put(123,122){$y_1^{'}$}
\put(89.5,132){$x_2$}
\put(162.5,123.25){$x_2^{'}$}
\end{picture}

\end{figure}

Let $y_2=t_{i-1}\dots t_j,$ where $t_l$ passes through $\Gamma_l$ (and $\Delta_l$) for $l=i-1,\dots,j.$  
Since every subpath $t_l$  connects two 
 vertices on the $k$-bands of the $l$-sector $\Delta_l,$ we can construct the mirror copy $\Delta'_{2i-l-1}$  (which is a $2i-l+1$-sector) of $\Delta_l$ 
 or the replica of $\Delta_l$ (if $2i-l-1 =1$), 
 and the copies of the vertices $(t_l)_{\pm}$ are connected in $\Delta'_{2i-l-1}$ by a path $t'_{2i-l-1}$ with
 $|t'_{2i-l-1}|\le |t_l|$ by  Lemma \ref{replica}. The desired diagram $\Delta'$ embeds in the union of these
 $\Delta'_{2i-l-1}$-s, and $y'_2 = t'_{i}\dots t'_{2i-j-1}.$ \endproof

\subsection{Shortcuts} \label{sc}

In this subsection we show that Lemma \ref{sravn} helps  cut off a  hub from a diagram using a 
'shortcut'. But at first we consider a few statements about simpler pieces (without hubs) of diagrams over the group $G$.

 \begin{lemma}\label{qband} 
  Suppose a diagram $\Delta$ over $G$ has a $q$-band $\cal C$ starting and ending on $\partial\Delta,$
 and $p$ is a side of $\cal C.$ Assume that no $\theta$-band crosses $\cal C$ twice in $\Delta$
 and that there is a factorization $xy$ of the boundary
 path of $\Delta$ such that $x_-=p_-, x_+=p_+.$ 
 Then $|p|\le |x|,$ $|\partial{\cal C}|\le |\partial \Delta|$,
 and $|xp^{-1}|\le |\partial \Delta|.$
  \end{lemma}
  
 \proof On the one hand, the length $|p|$ of $p$ is equal to the number $m$ of $\theta$-cells in $\cal C,$
 since  every cell of $\cal C$ has one $\theta$-edge and at most one $a$-edge on $p$ by
 condition 3 of Lemma \ref{MtocalS} and the definition of $(\theta,q)$-relations. On the other hand, every maximal $\theta$-band 
 crossing $\cal C$ must terminate on $x,$ since it does not cross $\cal C$ twice. It follows that
 $|x|\ge m,$ and so $|p|\le |x|.$ Similarly, we obtain inequalities $|\partial {\cal C}|\le |\partial \Delta|$,
 and $|xp^{-1}|\le |\partial \Delta|.$ (We take into account the definition of length and the fact 
 that the sides of the band $\cal C$ are separated by two $q$-edges lying on $\partial\Delta$.)    \endproof

A maximal $\theta$-band of a diagram $\Delta$  called a
{\em rim band} if its start and end $\theta$-edges as well as its top or its 
bottom lies on the boundary path of $\Delta$.
\bigskip

\begin{lemma} \label{rim} Suppose $\Delta$ is a diagram over $G$
with a rim band $\ttt$ having at most $N$ $(\theta,q)$-cells. Denote by
$\Delta'$ the subdiagram $\Delta\backslash \ttt$. Then
$|\partial\Delta|-|\partial\Delta'|\ge 1/2$.
\end{lemma}

\proof Let $s$ be the top side of $\ttt$ and
$s\subset\partial\Delta$. Note that by our assumptions the
difference between the number of $a$-edges in the bottom $s'$ of
$\ttt$ and the number of $a$-edges in $s$ cannot be greater than
$2N$ since every $(\theta,q)$-cell has at most two $a$-edges. 
However, $\Delta'$ is obtained
by cutting off $\ttt$ along $s'$, and its boundary contains two fewer
$\theta$-edges than $\Delta$. 
Thus one can  compare  the boundaries of $\Delta$ and $\Delta'$ as follows.
There is a one-to-one correspondence between the $q$-edges of these
boundaries, and to extend this correspondence to the $\theta$- and
$a$-edges 
of the boundaries, one should remove
two $\theta$-edges from $\partial\Delta$ and add at most $2N$ $a$-edges.
Therefore it follows from the definition of length and by Inequality (\ref{param}), that
$|\partial\Delta|-|\partial\Delta'|\ge 2(1-\delta)-2N\delta >1/2$ . \endproof

\begin{lemma} \label{change} Suppose a diagram $\Delta$ has two cells: an $H$-cell $\Pi$
and a $(\theta,a)$-cell $\pi$ which have a common $a$-edge $e$, $ep$ is the boundary of
$\Pi$ and $efe'f'$ is the boundary of $\pi,$ where $\Lab(e)\equiv\Lab(e')^{-1}\equiv a$ and 
$\Lab(f)\equiv\Lab(f')^{-1}\equiv\theta$ for a $\theta$-letter $\theta.$ Then there is a diagram $\Delta'$ with the same 
boundary label as $\Delta$ composed of $\Pi$ and a $\theta$-band $\cal T,$ and $p$
is the side of $\cal T.$  
\end{lemma} 

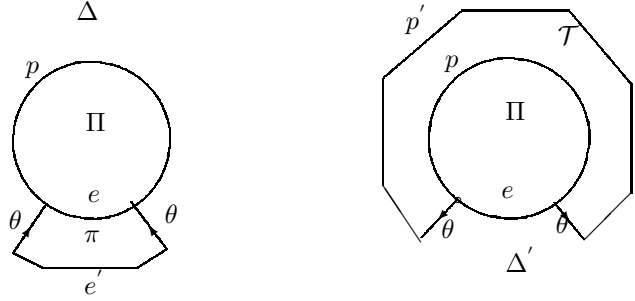
\begin{figure}[h!]
\unitlength 1mm 
\linethickness{0.4pt}
\ifx\plotpoint\undefined\newsavebox{\plotpoint}\fi 
\begin{picture}(104,60)(0,110)
\put(42.693,151.75){\line(0,1){.5611}}
\put(42.678,152.311){\line(0,1){.5595}}
\put(42.633,152.871){\line(0,1){.5562}}
\put(42.558,153.427){\line(0,1){.5514}}
\multiput(42.453,153.978)(-.03361,.13624){4}{\line(0,1){.13624}}
\multiput(42.318,154.523)(-.032709,.107391){5}{\line(0,1){.107391}}
\multiput(42.155,155.06)(-.032026,.087899){6}{\line(0,1){.087899}}
\multiput(41.963,155.588)(-.031459,.073758){7}{\line(0,1){.073758}}
\multiput(41.742,156.104)(-.030955,.062966){8}{\line(0,1){.062966}}
\multiput(41.495,156.608)(-.030483,.054411){9}{\line(0,1){.054411}}
\multiput(41.22,157.097)(-.033362,.052694){9}{\line(0,1){.052694}}
\multiput(40.92,157.572)(-.032531,.045743){10}{\line(0,1){.045743}}
\multiput(40.595,158.029)(-.031765,.039936){11}{\line(0,1){.039936}}
\multiput(40.245,158.468)(-.031043,.03499){12}{\line(0,1){.03499}}
\multiput(39.873,158.888)(-.032878,.033272){12}{\line(0,1){.033272}}
\multiput(39.478,159.287)(-.034618,.031457){12}{\line(-1,0){.034618}}
\multiput(39.063,159.665)(-.039554,.032238){11}{\line(-1,0){.039554}}
\multiput(38.628,160.019)(-.045352,.033073){10}{\line(-1,0){.045352}}
\multiput(38.174,160.35)(-.047064,.030589){10}{\line(-1,0){.047064}}
\multiput(37.704,160.656)(-.054043,.031129){9}{\line(-1,0){.054043}}
\multiput(37.217,160.936)(-.062593,.031703){8}{\line(-1,0){.062593}}
\multiput(36.717,161.19)(-.073378,.032336){7}{\line(-1,0){.073378}}
\multiput(36.203,161.416)(-.087511,.033071){6}{\line(-1,0){.087511}}
\multiput(35.678,161.615)(-.089161,.028322){6}{\line(-1,0){.089161}}
\multiput(35.143,161.785)(-.108665,.028188){5}{\line(-1,0){.108665}}
\put(34.6,161.926){\line(-1,0){.5501}}
\put(34.05,162.037){\line(-1,0){.5553}}
\put(33.494,162.119){\line(-1,0){.5589}}
\put(32.935,162.171){\line(-1,0){.5609}}
\put(32.374,162.193){\line(-1,0){.5612}}
\put(31.813,162.184){\line(-1,0){.56}}
\put(31.253,162.146){\line(-1,0){.5571}}
\put(30.696,162.077){\line(-1,0){.5526}}
\multiput(30.143,161.979)(-.13663,-.03199){4}{\line(-1,0){.13663}}
\multiput(29.597,161.851)(-.107773,-.031427){5}{\line(-1,0){.107773}}
\multiput(29.058,161.694)(-.088274,-.030977){6}{\line(-1,0){.088274}}
\multiput(28.528,161.508)(-.074127,-.030578){7}{\line(-1,0){.074127}}
\multiput(28.01,161.294)(-.06333,-.030202){8}{\line(-1,0){.06333}}
\multiput(27.503,161.052)(-.061616,-.033561){8}{\line(-1,0){.061616}}
\multiput(27.01,160.784)(-.053088,-.032732){9}{\line(-1,0){.053088}}
\multiput(26.532,160.489)(-.046127,-.031983){10}{\line(-1,0){.046127}}
\multiput(26.071,160.169)(-.040311,-.031287){11}{\line(-1,0){.040311}}
\multiput(25.627,159.825)(-.038572,-.033407){11}{\line(-1,0){.038572}}
\multiput(25.203,159.457)(-.033661,-.032479){12}{\line(-1,0){.033661}}
\multiput(24.799,159.068)(-.031867,-.034241){12}{\line(0,-1){.034241}}
\multiput(24.417,158.657)(-.032707,-.039167){11}{\line(0,-1){.039167}}
\multiput(24.057,158.226)(-.033611,-.044955){10}{\line(0,-1){.044955}}
\multiput(23.721,157.776)(-.031147,-.046696){10}{\line(0,-1){.046696}}
\multiput(23.409,157.309)(-.031771,-.053669){9}{\line(0,-1){.053669}}
\multiput(23.124,156.826)(-.032446,-.06221){8}{\line(0,-1){.06221}}
\multiput(22.864,156.329)(-.033208,-.072987){7}{\line(0,-1){.072987}}
\multiput(22.632,155.818)(-.029238,-.074666){7}{\line(0,-1){.074666}}
\multiput(22.427,155.295)(-.029382,-.088817){6}{\line(0,-1){.088817}}
\multiput(22.251,154.762)(-.029481,-.108321){5}{\line(0,-1){.108321}}
\put(22.103,154.221){\line(0,-1){.5487}}
\put(21.985,153.672){\line(0,-1){.5543}}
\put(21.897,153.118){\line(0,-1){.5582}}
\put(21.838,152.559){\line(0,-1){2.2402}}
\put(21.905,150.319){\line(0,-1){.5537}}
\multiput(21.997,149.765)(.03036,-.137){4}{\line(0,-1){.137}}
\multiput(22.118,149.217)(.030141,-.10814){5}{\line(0,-1){.10814}}
\multiput(22.269,148.677)(.029923,-.088637){6}{\line(0,-1){.088637}}
\multiput(22.449,148.145)(.029693,-.074487){7}{\line(0,-1){.074487}}
\multiput(22.657,147.624)(.033652,-.072784){7}{\line(0,-1){.072784}}
\multiput(22.892,147.114)(.032825,-.062012){8}{\line(0,-1){.062012}}
\multiput(23.155,146.618)(.032097,-.053474){9}{\line(0,-1){.053474}}
\multiput(23.444,146.137)(.031431,-.046505){10}{\line(0,-1){.046505}}
\multiput(23.758,145.672)(.030804,-.040681){11}{\line(0,-1){.040681}}
\multiput(24.097,145.224)(.032946,-.038967){11}{\line(0,-1){.038967}}
\multiput(24.459,144.796)(.032076,-.034046){12}{\line(0,-1){.034046}}
\multiput(24.844,144.387)(.033859,-.032273){12}{\line(1,0){.033859}}
\multiput(25.25,144)(.038775,-.033172){11}{\line(1,0){.038775}}
\multiput(25.677,143.635)(.040501,-.03104){11}{\line(1,0){.040501}}
\multiput(26.122,143.293)(.046321,-.031702){10}{\line(1,0){.046321}}
\multiput(26.586,142.976)(.053286,-.032408){9}{\line(1,0){.053286}}
\multiput(27.065,142.685)(.061819,-.033185){8}{\line(1,0){.061819}}
\multiput(27.56,142.419)(.063513,-.029816){8}{\line(1,0){.063513}}
\multiput(28.068,142.181)(.074312,-.030126){7}{\line(1,0){.074312}}
\multiput(28.588,141.97)(.088461,-.030438){6}{\line(1,0){.088461}}
\multiput(29.119,141.787)(.107962,-.030769){5}{\line(1,0){.107962}}
\multiput(29.659,141.633)(.13682,-.03116){4}{\line(1,0){.13682}}
\put(30.206,141.509){\line(1,0){.5532}}
\put(30.759,141.414){\line(1,0){.5575}}
\put(31.317,141.348){\line(1,0){.5602}}
\put(31.877,141.313){\line(1,0){.5613}}
\put(32.438,141.308){\line(1,0){.5607}}
\put(32.999,141.334){\line(1,0){.5586}}
\put(33.557,141.389){\line(1,0){.5548}}
\put(34.112,141.474){\line(1,0){.5494}}
\multiput(34.662,141.589)(.108491,.02885){5}{\line(1,0){.108491}}
\multiput(35.204,141.733)(.088987,.028865){6}{\line(1,0){.088987}}
\multiput(35.738,141.906)(.087307,.033604){6}{\line(1,0){.087307}}
\multiput(36.262,142.108)(.073179,.032783){7}{\line(1,0){.073179}}
\multiput(36.774,142.338)(.062398,.032084){8}{\line(1,0){.062398}}
\multiput(37.273,142.594)(.053853,.031458){9}{\line(1,0){.053853}}
\multiput(37.758,142.877)(.046876,.030875){10}{\line(1,0){.046876}}
\multiput(38.227,143.186)(.04515,.033349){10}{\line(1,0){.04515}}
\multiput(38.678,143.52)(.039357,.032479){11}{\line(1,0){.039357}}
\multiput(39.111,143.877)(.034425,.031668){12}{\line(1,0){.034425}}
\multiput(39.524,144.257)(.032674,.033472){12}{\line(0,1){.033472}}
\multiput(39.916,144.658)(.033631,.038377){11}{\line(0,1){.038377}}
\multiput(40.286,145.081)(.031521,.040128){11}{\line(0,1){.040128}}
\multiput(40.633,145.522)(.032251,.04594){10}{\line(0,1){.04594}}
\multiput(40.956,145.981)(.03304,.052897){9}{\line(0,1){.052897}}
\multiput(41.253,146.458)(.03015,.054595){9}{\line(0,1){.054595}}
\multiput(41.524,146.949)(.03057,.063153){8}{\line(0,1){.063153}}
\multiput(41.769,147.454)(.031009,.073948){7}{\line(0,1){.073948}}
\multiput(41.986,147.972)(.03149,.088092){6}{\line(0,1){.088092}}
\multiput(42.175,148.5)(.032053,.107588){5}{\line(0,1){.107588}}
\multiput(42.335,149.038)(.03278,.13644){4}{\line(0,1){.13644}}
\put(42.466,149.584){\line(0,1){.552}}
\put(42.568,150.136){\line(0,1){.5567}}
\put(42.64,150.693){\line(0,1){1.0573}}
\multiput(21.75,136.75)(.06666667,-.03333333){60}{\line(1,0){.06666667}}
\put(25.75,134.75){\line(1,0){12.5}}
\multiput(38.25,134.75)(.05333333,.03333333){75}{\line(1,0){.05333333}}
\put(98.43,152){\line(0,1){.5715}}
\put(98.415,152.572){\line(0,1){.5699}}
\put(98.369,153.141){\line(0,1){.5666}}
\put(98.293,153.708){\line(0,1){.5617}}
\multiput(98.186,154.27)(-.027282,.111041){5}{\line(0,1){.111041}}
\multiput(98.05,154.825)(-.033185,.109422){5}{\line(0,1){.109422}}
\multiput(97.884,155.372)(-.032494,.089574){6}{\line(0,1){.089574}}
\multiput(97.689,155.909)(-.03192,.075177){7}{\line(0,1){.075177}}
\multiput(97.465,156.436)(-.03141,.064191){8}{\line(0,1){.064191}}
\multiput(97.214,156.949)(-.030934,.055483){9}{\line(0,1){.055483}}
\multiput(96.936,157.449)(-.030473,.048374){10}{\line(0,1){.048374}}
\multiput(96.631,157.932)(-.033017,.046674){10}{\line(0,1){.046674}}
\multiput(96.301,158.399)(-.032243,.040764){11}{\line(0,1){.040764}}
\multiput(95.946,158.847)(-.031514,.035731){12}{\line(0,1){.035731}}
\multiput(95.568,159.276)(-.033381,.033994){12}{\line(0,1){.033994}}
\multiput(95.167,159.684)(-.035152,.032159){12}{\line(-1,0){.035152}}
\multiput(94.745,160.07)(-.04017,.03298){11}{\line(-1,0){.04017}}
\multiput(94.304,160.433)(-.041877,.030783){11}{\line(-1,0){.041877}}
\multiput(93.843,160.771)(-.047811,.031348){10}{\line(-1,0){.047811}}
\multiput(93.365,161.085)(-.054911,.031938){9}{\line(-1,0){.054911}}
\multiput(92.871,161.372)(-.063609,.032574){8}{\line(-1,0){.063609}}
\multiput(92.362,161.633)(-.074584,.033284){7}{\line(-1,0){.074584}}
\multiput(91.84,161.866)(-.076258,.029245){7}{\line(-1,0){.076258}}
\multiput(91.306,162.071)(-.090666,.029309){6}{\line(-1,0){.090666}}
\multiput(90.762,162.247)(-.110526,.029298){5}{\line(-1,0){.110526}}
\put(90.209,162.393){\line(-1,0){.5597}}
\put(89.65,162.51){\line(-1,0){.5651}}
\put(89.084,162.596){\line(-1,0){.5689}}
\put(88.516,162.653){\line(-1,0){.5711}}
\put(87.944,162.678){\line(-1,0){.5717}}
\put(87.373,162.673){\line(-1,0){.5706}}
\put(86.802,162.638){\line(-1,0){.5679}}
\put(86.234,162.572){\line(-1,0){.5636}}
\multiput(85.671,162.476)(-.1394,-.03157){4}{\line(-1,0){.1394}}
\multiput(85.113,162.349)(-.110007,-.031187){5}{\line(-1,0){.110007}}
\multiput(84.563,162.193)(-.090151,-.030858){6}{\line(-1,0){.090151}}
\multiput(84.022,162.008)(-.075746,-.030547){7}{\line(-1,0){.075746}}
\multiput(83.492,161.794)(-.064753,-.030237){8}{\line(-1,0){.064753}}
\multiput(82.974,161.553)(-.063042,-.033659){8}{\line(-1,0){.063042}}
\multiput(82.47,161.283)(-.054356,-.032874){9}{\line(-1,0){.054356}}
\multiput(81.98,160.987)(-.047267,-.032162){10}{\line(-1,0){.047267}}
\multiput(81.508,160.666)(-.041344,-.031496){11}{\line(-1,0){.041344}}
\multiput(81.053,160.319)(-.039599,-.033663){11}{\line(-1,0){.039599}}
\multiput(80.617,159.949)(-.034596,-.032756){12}{\line(-1,0){.034596}}
\multiput(80.202,159.556)(-.032793,-.034561){12}{\line(0,-1){.034561}}
\multiput(79.809,159.141)(-.033706,-.039563){11}{\line(0,-1){.039563}}
\multiput(79.438,158.706)(-.03154,-.04131){11}{\line(0,-1){.04131}}
\multiput(79.091,158.252)(-.032213,-.047232){10}{\line(0,-1){.047232}}
\multiput(78.769,157.779)(-.032933,-.054321){9}{\line(0,-1){.054321}}
\multiput(78.472,157.29)(-.033726,-.063006){8}{\line(0,-1){.063006}}
\multiput(78.203,156.786)(-.030306,-.06472){8}{\line(0,-1){.06472}}
\multiput(77.96,156.269)(-.030628,-.075713){7}{\line(0,-1){.075713}}
\multiput(77.746,155.739)(-.030954,-.090118){6}{\line(0,-1){.090118}}
\multiput(77.56,155.198)(-.031305,-.109974){5}{\line(0,-1){.109974}}
\multiput(77.403,154.648)(-.03172,-.13936){4}{\line(0,-1){.13936}}
\put(77.277,154.091){\line(0,-1){.5634}}
\put(77.18,153.527){\line(0,-1){.5678}}
\put(77.113,152.959){\line(0,-1){1.1423}}
\put(77.072,151.817){\line(0,-1){.5712}}
\put(77.097,151.246){\line(0,-1){.569}}
\put(77.152,150.677){\line(0,-1){.5652}}
\put(77.238,150.112){\line(0,-1){.5598}}
\multiput(77.354,149.552)(.02918,-.110557){5}{\line(0,-1){.110557}}
\multiput(77.5,148.999)(.029212,-.090698){6}{\line(0,-1){.090698}}
\multiput(77.676,148.455)(.029163,-.076289){7}{\line(0,-1){.076289}}
\multiput(77.88,147.921)(.033204,-.07462){7}{\line(0,-1){.07462}}
\multiput(78.112,147.399)(.032505,-.063644){8}{\line(0,-1){.063644}}
\multiput(78.372,146.889)(.03188,-.054945){9}{\line(0,-1){.054945}}
\multiput(78.659,146.395)(.031297,-.047844){10}{\line(0,-1){.047844}}
\multiput(78.972,145.916)(.030738,-.04191){11}{\line(0,-1){.04191}}
\multiput(79.31,145.455)(.032937,-.040205){11}{\line(0,-1){.040205}}
\multiput(79.672,145.013)(.032121,-.035186){12}{\line(0,-1){.035186}}
\multiput(80.058,144.591)(.033958,-.033417){12}{\line(1,0){.033958}}
\multiput(80.465,144.19)(.035698,-.031552){12}{\line(1,0){.035698}}
\multiput(80.894,143.811)(.040729,-.032287){11}{\line(1,0){.040729}}
\multiput(81.342,143.456)(.046638,-.033067){10}{\line(1,0){.046638}}
\multiput(81.808,143.125)(.048341,-.030524){10}{\line(1,0){.048341}}
\multiput(82.292,142.82)(.05545,-.030993){9}{\line(1,0){.05545}}
\multiput(82.791,142.541)(.064158,-.031479){8}{\line(1,0){.064158}}
\multiput(83.304,142.289)(.075143,-.032001){7}{\line(1,0){.075143}}
\multiput(83.83,142.065)(.089539,-.03259){6}{\line(1,0){.089539}}
\multiput(84.367,141.87)(.109386,-.033302){5}{\line(1,0){.109386}}
\multiput(84.914,141.703)(.111011,-.027401){5}{\line(1,0){.111011}}
\put(85.469,141.566){\line(1,0){.5616}}
\put(86.031,141.459){\line(1,0){.5665}}
\put(86.597,141.382){\line(1,0){.5698}}
\put(87.167,141.336){\line(1,0){1.143}}
\put(88.31,141.335){\line(1,0){.5699}}
\put(88.88,141.38){\line(1,0){.5667}}
\put(89.447,141.456){\line(1,0){.5618}}
\multiput(90.009,141.562)(.11107,.027163){5}{\line(1,0){.11107}}
\multiput(90.564,141.697)(.109457,.033067){5}{\line(1,0){.109457}}
\multiput(91.111,141.863)(.089609,.032398){6}{\line(1,0){.089609}}
\multiput(91.649,142.057)(.075212,.03184){7}{\line(1,0){.075212}}
\multiput(92.175,142.28)(.064225,.031342){8}{\line(1,0){.064225}}
\multiput(92.689,142.531)(.055516,.030874){9}{\line(1,0){.055516}}
\multiput(93.189,142.809)(.048406,.030421){10}{\line(1,0){.048406}}
\multiput(93.673,143.113)(.046709,.032967){10}{\line(1,0){.046709}}
\multiput(94.14,143.442)(.040798,.0322){11}{\line(1,0){.040798}}
\multiput(94.589,143.797)(.035765,.031475){12}{\line(1,0){.035765}}
\multiput(95.018,144.174)(.03403,.033344){12}{\line(1,0){.03403}}
\multiput(95.426,144.574)(.032196,.035117){12}{\line(0,1){.035117}}
\multiput(95.813,144.996)(.033023,.040135){11}{\line(0,1){.040135}}
\multiput(96.176,145.437)(.030828,.041844){11}{\line(0,1){.041844}}
\multiput(96.515,145.898)(.031399,.047777){10}{\line(0,1){.047777}}
\multiput(96.829,146.375)(.031997,.054877){9}{\line(0,1){.054877}}
\multiput(97.117,146.869)(.032642,.063574){8}{\line(0,1){.063574}}
\multiput(97.378,147.378)(.033363,.074548){7}{\line(0,1){.074548}}
\multiput(97.612,147.9)(.029326,.076227){7}{\line(0,1){.076227}}
\multiput(97.817,148.433)(.029406,.090635){6}{\line(0,1){.090635}}
\multiput(97.993,148.977)(.029417,.110494){5}{\line(0,1){.110494}}
\put(98.14,149.53){\line(0,1){.5595}}
\put(98.258,150.089){\line(0,1){.565}}
\put(98.345,150.654){\line(0,1){.5689}}
\put(98.402,151.223){\line(0,1){.7769}}
\put(81,143.75){\circle{.707}}
\put(76,138.25){\line(-2,3){5}}
\put(71,145.75){\line(0,1){14.25}}
\multiput(71,160)(.0393258427,.0337078652){267}{\line(1,0){.0393258427}}
\put(81.5,169){\line(1,0){14.25}}
\multiput(95.75,169)(.0336734694,-.0397959184){245}{\line(0,-1){.0397959184}}
\put(104,159.25){\line(0,-1){14.5}}
\put(104,144.75){\line(-1,-1){6.25}}
\put(23.5,160.75){$p$}
\put(31.75,143.5){$e$}
\put(31.5,153){$\Pi$}
\put(31.25,138.25){$\pi$}
\put(31.25,131.25){$e^{'}$}
\put(24.125,140.125){\vector(2,3){.07}}\multiput(22,137)(.033730159,.049603175){126}{\line(0,1){.049603175}}
\put(39.875,140.375){\vector(-3,4){.07}}\multiput(42.25,137.25)(-.033687943,.044326241){141}{\line(0,1){.044326241}}
\put(21.25,140){$\theta$}
\put(42,140.75){$\theta$}
\put(30.25,167.75){$\Delta$}
\put(87.25,134.25){$\Delta^{'}$}
\put(78.375,141.25){\vector(-1,-1){.07}}\multiput(80.75,143.75)(-.033687943,-.035460993){141}{\line(0,-1){.035460993}}
\put(95.75,141){\vector(3,-4){.07}}\multiput(93.75,143.5)(.033613445,-.042016807){119}{\line(0,-1){.042016807}}
\put(78.75,138.75){$\theta$}
\put(94,139.5){$\theta$}
\put(87.25,154.5){$\Pi$}
\put(94.25,164.5){$\cal T$}
\put(79.25,161.5){$p$}
\put(74,166.5){$p^{'}$}
\put(86.75,144){$e$}
\end{picture}

\end{figure}

\proof The letter $a$ commutes with $\theta,$ and therefore every letter of the
boundary label of $\Pi$ commutes with $\theta.$ Hence one can construct a $\theta$-band $\cal T$
with boundary $g'p^{-1}gp',$ where $\Lab(p')\equiv\Lab(p)$, and $\Lab(g)\equiv\Lab(g')^{-1}\equiv\theta^{-1}.$
If we attach the band $\cal T$  to $\Pi$ along the path $p$ and remove $\pi,$ we obtain
the required diagram $\Delta'$. \endproof

 Let $\Pi$ be a hub
of a minimal diagram $\Delta$ given by Lemma \ref{extdisc}. Using the notation of
that lemma, we recall that the subdiagrams $\Gamma_i$ and $\Gamma_{i+1}$ intersect along the 
$k$-band ${\cal B}_{i+1}$ ($i=1,\dots ,L-5$). We denote by $\Psi$ the minimal subdiagram 
containing all the $\Gamma_i$-s for $i=1,\dots, L-4.$ The boundary path of $\Psi$ is 
$x'x'',$ where $x'$ is composed of the sides of ${\cal B}_1,$ ${\cal B}_{L-3},$ and
a subpath of $\partial\Pi,$ while $x''$ is a subpath of $\partial\Delta.$ 

\begin{lemma} \label{svojstva} One can construct a minimal diagram $\Psi'$ over $G_1$ with boundary path 
$x' \bar x$ such that
(1) $\Psi'$ includes the bands ${\cal B}_1$ and ${\cal B}_{L-3}$  
(2) $|\bar x|\le |x''|,$ (3) the subdiagram $\Psi'$ 
has no maximal $q$-bands except for the $q$-bands ${\cal B'}_i$ ($i=1,\dots, L-3$) starting on $x'$,
(4) every maximal $\theta$-band of $\Psi'$ crosses either the band ${\cal B'}_1={\cal B}_1$ 
or the band ${\cal B'}_{L-3}={\cal B}_{L-3}$,
(5) the subdiagram $\Psi'$ has no $H$-cells between the pair of $k$-bands ${\cal B'}_i$ and
${\cal B'}_{i+1},$ unless  this pair is a pair of $k_1$- and $k_2$-bands. 
\end{lemma}

\proof Let $\Psi'$ be a minimal diagram with boundary of the form $x'\bar x$ which satisfies
conditions (1) and (2) and has minimal
$|\bar x|.$ Since $\Psi$ satisfies conditions (1) and (2) with $\bar x=x'$, such $\Psi'$ exists.  
Clearly the path $\bar x$ has no loops.
If the diagram $\Psi'$ has a maximal $q$-band $\cal C$ which does not start/terminate on $\Pi,$
then one can cut off $\cal C$  and shorten $\bar x$ by Lemma \ref{qband}. So $\Psi'$ satisfies
condition (3). 

Assume that the diagram $\Psi'$ does not satisfy condition (4) of the lemma. Then by Lemma \ref{NoAnnul}, we
have a $\theta$-band of $\Psi'$ starting and terminating on $\bar x.$ It follows that there is a $\theta$-band
$\cal T$ starting and terminating on $\bar x,$ such that the subdiagram $\Phi$ bounded by $\cal T$ and a
part $y$ of $\bar x$ has no non-trivial $\theta$-bands, i.e., it contains only $H$-cells. Two $H$-cells of $\Phi$ cannot have a common edge since otherwise they
could be replaced by one $H$-cell contrary to the minimality of $\Psi'.$  It follows that every $H$-cell $\pi$
of $\Phi$ has a common edge with $\cal T$ because the path $\bar x$ has no loops.

Thus we have a series of $H$-cells $\pi_1,\dots,\pi_s$ in $\Phi$ with boundaries $y_iz_i$ ($i=1,\dots,s$),
where $y_i$ is a part of $y$ (or $y_i$ is empty) and $z_i$ belongs to the side $z$ of $\cal T$ . 
If $\sum |\Lab(z_i)|_a\le 2,$ then $\sum |z_i|\le 2\delta$ since every $z_i$ is a product of $a$-edges. Therefore $|z|\le |\bar x|+2\delta.$
It follows from Lemma \ref{rim} that if we remove all $\pi_i$-s and then cut off the band $\cal T$, then
we decrease the length of $\bar x$ since $2\delta<1/2;$ a contradiction.  

Hence  $\sum |\Lab(z_i)|_a\ge 3,$ and so at least one of $\pi_i$-s has a common edge with a $(\theta,a)$-cell
of $\cal T.$ (We recall that $\cal T$ intersects each of the maximal $q$-bands of $\Psi'$ starting on $x''$ at most once by Lemma \ref{NoAnnul}, and so at most
two of the $(\theta,q)$-cells  of $\cal T$ have $a$-edges with $a\in{\cal A}_1$, and each of these two
 $(\theta,q)$-cells can have at most one $a$-edge with label from ${\cal A}_1$ .) Therefore we can apply Lemma \ref{change} to replace
 the $(\theta,a)$-cell by a $\theta$-band passing around the cell $\pi_i$. This modification of the band
 $\cal T$ decreases the number of $H$-cells in $\Phi$, since one of the $H$-cells gets moved over the band $\cal T.$
 
 The modified diagram can be non-minimal, but our surgery preserves $q$-bands and keeps the property
 that every $q$-band and $\theta$-band have at most one common $(\theta,q)$-cell. So sooner or later, this trick makes the inequality $\sum |\Lab(z_i)|_a\le 2$ true, and one can decrease $\bar x,$
 as was explained above. If one replaces the obtained diagram by a minimal one, then  condition
 (1) still holds since the maximal $q$-bands ${\cal B}_1$ and ${\cal B}_{L-3}$ are completely determined
 by the boundary as this follows from Lemma \ref{NoAnnul}.  This contradicts  the assumption
 on the minimality of $|\bar x|.$ Thus   property  (4) holds.
 
 If $\Psi'$ has  $H$-cells in the subdiagram $\Gamma'_i$ between the pair of $k$-bands ${\cal B'}_i$ and
${\cal B'}_{i+1}$ which are not a pair of $k_1$- and $k_2$-bands, then these $H$-cells
(with labels over the alphabet ${\cal A}_1$) cannot have common edges with the maximal $\theta$-bands
of $\Gamma'_i$ since the $\theta$-bands of $\Gamma'_i$ must intersect either ${\cal B'}_i$ or
${\cal B'}_{i+1}$ by (4). This implies that $\Gamma'_i$ has no $H$-cells at all because the path
$\bar x$ has no loops. The lemma is proved. \endproof

\begin{lemma}\label{cut} If a minimal diagram $\Delta$ has a hub, then the (cyclic
shift of the) boundary path of $\Delta$ can be factorized as $pp'$ so that the
subpath $p$ starts and ends with $q$-edges and there is 
a simple path $z$ in $\Delta$ with $z_-=p_-, z_+=p_+$ such that the subdiagram bounded by the 
loop $pz^{-1}$ has exactly one hub $\Pi$ and the label $\Lab(z)$ is equal in the group $G_1$ to a 
word of length $<|p|.$
\end{lemma}

\proof We may assume that a hub $\Pi$ is chosen in $\Delta$ according to Lemma \ref{extdisc}.
Let $x'x''$ be the boundary path of $\Psi$ as in Lemma \ref{svojstva}. We will look for
the path $z$ in the minimal subdiagram $\Delta'$ obtained after removing  $\Psi$ and $\Pi$
from $\Delta.$ Therefore to prove the
lemma, one may assume using the notation of Lemma \ref{svojstva}, that $\Psi'=\Psi$ and 
$\bar x=x''$, i.e., the subdiagram $\Psi$ itself has properties (3), (4), and (5) from
Lemma \ref{svojstva}. (We do not know if the diagram $\Psi'\cup\Pi\cup\Delta'$ is still 
minimal but we will use the minimality of $\Psi'$ only.)
 
There is $d\ge 0$ such that there exist exactly $d$
maximal $\theta$-bands of $\Psi$ crossing each of the $k$-bands ${\cal B}_1,\dots,{\cal B}_{L-3}.$ 
This implies that the initial subbands ${\cal B}_i[d]$ of length $d$ in all  ${\cal B}_i$-s
($i=1,\dots, L-3$) are copies of each other under the shifts of the indices in their
boundary labels.

Let $T$ be the set of remaining maximal $\theta$-bands of $\Psi,$ i.e., every band of $T$ intersects 
exactly one of the bands ${\cal B}_1$, ${\cal B}_{L-3}.$
It follows from Lemma \ref{NoAnnul} for $\Psi$ that there is an integer $l$ ($1\le l<L-3$) such that no $\theta$-band of $T$ crossing 
${\cal B}_1$ crosses ${\cal B}_{l+1}$ and no $\theta$-band of $T$ crossing 
${\cal B}_{L-3}$ crosses ${\cal B}_{l}.$ We have either $(L-3)-l<(L-3)/2 $ or $(l+1)-1<(L-3)/2$
since $L$ is even. Without loss of generality we choose the former inequality, and so $l\ge (L-2)/2.$ 

If ${\cal B}_i$
is $k_1$-band for some $i\le 6$, then we will consider a smaller subdiagram bounded by ${\cal B}_7$ and
${\cal B}_{L-3}$ instead of $\Psi$. (Respectively, we change the complimentary minimal subdiagram $\Delta'$.) If none of ${\cal B}_1,\dots, {\cal B}_6$ is a $k_1$-band, then we 
do not change $\Psi.$  Thus, in any case we can reindex the $k$-bands and assume that
$\Psi$ is bounded by  ${\cal B}_1$ and ${\cal B}_{L-r}$ for some $r\le 9$
and that the bands ${\cal B}_1,\dots ,{\cal B}_{r+3}$ are not $k_1$-bands. Let them be $k_i$-,..., 
$k_{i\pm (r+2)}$-bands for some $i.$ We will assume that ${\cal B}_{r+3}$
is a $k_{i-r-2}$-band. 

Since $L\ge 40, $ after such a reindexing $l\ge (L-2)/2-6 > 12\ge r+3,$ and so no $\theta$-band from the set $T$ crossing the band ${\cal B}_{L-3}$
crosses ${\cal B}_{r+3}.$ We denote by $\Phi$ the part of the diagram $\Psi$ bounded by ${\cal B}_{r+3}$
and ${\cal B}_1.$ Let ${\cal T}^{\Phi}_1,\dots,{\cal T}^{\Phi}_s$ be the maximal $\theta$-bands of $\Phi
$. Then by the choice of the subdiagrams $\Psi$ and $\Phi,$ every cell of $\Phi$ belongs to one of these
$\theta$-bands, and each of these bands crosses the $k_i$-band ${\cal B}_1.$ We will assume that ${\cal T}^{\Phi}_1$
is the closest band to the hub $\Pi$, and so on. 

For every $j\ge 1,$ at least one of the two
$q$-edges of every $(\theta,q)$-cell of ${\cal T}^{\Phi}_i$ belongs to the top of ${\cal T}^{\Phi}_{j-1}$
(to  $\partial\Pi$ if $j=1$) because $\Psi$ has no maximal $q$-bands except for the bands starting on $\Pi.$
Assume that a band ${\cal T}^{\Phi}_j$, starting with a $(\theta,q)$-cell of ${\cal B}_1,$ terminates with
an $(\theta,a)$-cell $\pi$ having no $a$-edges on $\partial {\cal T}^{\Phi}_{j-1}.$ Then an $a$-edge and one
$\theta$-edge of $\pi$ lie on the boundary subpath $x''$ of $\Psi,$ and so if one removes $\pi$ from $\Psi,$
the length of $x''$ does not increase because $\partial\pi$ has two $\theta$-edges and two $a$-edges  
and all properties (3)-(5) hold for the
remaining part of $\Psi.$ Therefore we may assume that every edge $f$ of  $\tbott({\cal T}^{\Phi}_j)$ 
belongs to $\topp({\cal T}^{\Phi}_{j-1})$ (to $\partial\Pi$ for $j=1$). 

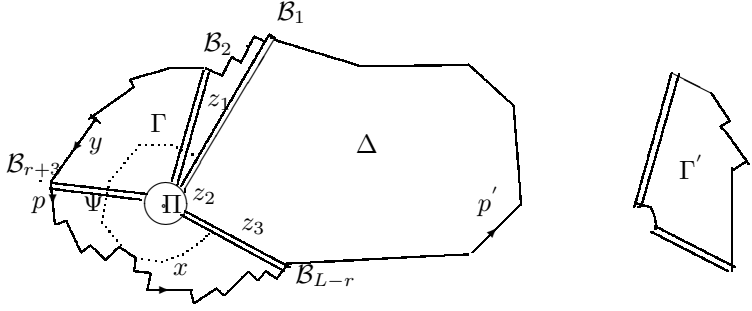
\begin{figure}[h!]
\unitlength 1mm 
\linethickness{0.4pt}
\ifx\plotpoint\undefined\newsavebox{\plotpoint}\fi 
\begin{picture}(98.75,60)(0,125)
\put(20.75,154){$\Pi$}
\put(21,155){\circle{.5}}
\put(21.25,155.25){\circle{5.852}}
\multiput(26.25,173)(-.033730159,-.117063492){126}{\line(0,-1){.117063492}}
\multiput(27,172.75)(-.033730159,-.119047619){126}{\line(0,-1){.119047619}}
\put(4.5,158.25){\line(1,0){.25}}
\multiput(6,158)(.28333333,-.03333333){45}{\line(1,0){.28333333}}
\multiput(5.75,157.25)(.27222222,-.03333333){45}{\line(1,0){.27222222}}
\multiput(23.75,154.25)(.062780269,-.033632287){223}{\line(1,0){.062780269}}
\put(23.5,153.75){\line(1,0){.25}}
\multiput(23.25,153.75)(.061659193,-.033632287){223}{\line(1,0){.061659193}}
\put(32.25,176.25){\line(-1,0){.25}}
\put(23.75,156.75){\line(0,1){.5}}
\multiput(23.25,157.5)(.0336676218,.0580229226){349}{\line(0,1){.0580229226}}
\multiput(35,177.75)(.05,-.0333333){15}{\line(1,0){.05}}
\put(23.5,157.25){\line(0,1){0}}
\put(23.75,157.25){\line(3,5){12}}
\multiput(26.5,173.25)(.075,-.0333333){30}{\line(1,0){.075}}
\multiput(28.75,172.25)(.03289474,.06578947){38}{\line(0,1){.06578947}}
\multiput(30,174.75)(.0543478,-.0326087){23}{\line(1,0){.0543478}}
\multiput(31.25,174)(.03365385,.06730769){52}{\line(0,1){.06730769}}
\multiput(33,177.5)(.0652174,-.0326087){23}{\line(1,0){.0652174}}
\put(26.25,173.25){\line(-1,0){4.5}}
\multiput(21.75,173.25)(-.11666667,-.03333333){45}{\line(-1,0){.11666667}}
\multiput(16.5,171.75)(.0333333,-.0666667){15}{\line(0,-1){.0666667}}
\multiput(17,170.75)(-.047619048,-.033730159){126}{\line(-1,0){.047619048}}
\multiput(11.25,167)(.05,-.0333333){15}{\line(1,0){.05}}
\multiput(12,166.5)(-.09375,.03125){8}{\line(-1,0){.09375}}
\multiput(11.5,166.75)(.0434783,-.0326087){23}{\line(1,0){.0434783}}
\multiput(6.75,152.5)(.0833333,.0333333){15}{\line(1,0){.0833333}}
\multiput(8,153)(.03333333,-.07083333){60}{\line(0,-1){.07083333}}
\multiput(10,148.75)(.0416667,.0333333){30}{\line(1,0){.0416667}}
\multiput(11.25,149.75)(.03333333,-.03333333){60}{\line(0,-1){.03333333}}
\multiput(13.25,147.75)(.03365385,.03365385){52}{\line(0,1){.03365385}}
\multiput(15,149.5)(.05,-.03333333){45}{\line(1,0){.05}}
\multiput(17.25,148)(-.03125,-.21875){8}{\line(0,-1){.21875}}
\multiput(17,146.25)(.28125,-.03125){8}{\line(1,0){.28125}}
\multiput(19.25,146)(-.0333333,-.15){15}{\line(0,-1){.15}}
\multiput(24.5,143.75)(.0333333,-.1166667){15}{\line(0,-1){.1166667}}
\multiput(25,142)(.05,.03333333){60}{\line(1,0){.05}}
\multiput(28,144)(.0333333,-.0333333){30}{\line(0,-1){.0333333}}
\multiput(29,143)(.033505155,.033505155){97}{\line(0,1){.033505155}}
\multiput(32.25,146.25)(.03289474,-.03289474){38}{\line(0,-1){.03289474}}
\multiput(33.5,145)(.0333333,.0416667){30}{\line(0,1){.0416667}}
\multiput(34.5,146.25)(.05,-.0333333){30}{\line(1,0){.05}}
\multiput(36,145.25)(.03333333,.04444444){45}{\line(0,1){.04444444}}
\multiput(15.18,149.68)(-.5,.6){6}{{\rule{.4pt}{.4pt}}}
\multiput(12.68,152.68)(.14286,.78571){8}{{\rule{.4pt}{.4pt}}}
\multiput(13.68,158.18)(.5,.59375){9}{{\rule{.4pt}{.4pt}}}
\multiput(17.68,162.93)(.85,0){6}{{\rule{.4pt}{.4pt}}}
\multiput(21.93,162.93)(.7,-.3){6}{{\rule{.4pt}{.4pt}}}
\multiput(25.43,161.43)(0,0){3}{{\rule{.4pt}{.4pt}}}
\multiput(17.18,147.93)(.875,-.125){5}{{\rule{.4pt}{.4pt}}}
\multiput(20.68,147.43)(.6875,.25){5}{{\rule{.4pt}{.4pt}}}
\multiput(23.43,148.43)(.625,.4375){5}{{\rule{.4pt}{.4pt}}}
\multiput(25.93,150.18)(.3333,.3333){4}{{\rule{.4pt}{.4pt}}}
\multiput(35.43,177.18)(-.125,0){3}{{\rule{.4pt}{.4pt}}}
\multiput(35.75,177)(.110576923,-.033653846){104}{\line(1,0){.110576923}}
\multiput(47.25,173.5)(1.78125,.03125){8}{\line(1,0){1.78125}}
\multiput(61.5,173.75)(.036585366,-.033536585){164}{\line(1,0){.036585366}}
\multiput(67.5,168.25)(.0333333,-.4416667){30}{\line(0,-1){.4416667}}
\multiput(61.5,148.5)(-.63815789,-.03289474){38}{\line(-1,0){.63815789}}
\put(26.75,168.25){$z_1$}
\put(24.75,155.75){$z_2$}
\put(31.25,151.75){$z_3$}
\put(46.5,162){$\Delta$}
\put(19.25,164.75){$\Gamma$}
\put(11,162.5){$y$}
\put(8.875,162.25){\vector(-3,-4){.07}}\multiput(11.75,166.25)(-.033625731,-.046783626){171}{\line(0,-1){.046783626}}
\put(22.25,145.75){$x$}
\put(21.75,143.75){\vector(1,0){.07}}\put(18.75,143.75){\line(1,0){6}}
\put(6.25,155.125){\vector(0,-1){.07}}\multiput(6,157.75)(.0333333,-.35){15}{\line(0,-1){.35}}
\put(3.5,154.75){$p$}
\put(62.75,154.25){$p^{'}$}
\put(65.25,152){\vector(1,1){.07}}\multiput(62,148.75)(.033678756,.033678756){193}{\line(0,1){.033678756}}
\put(36,179.75){${\cal B}_1$}
\put(26.25,175.5){${\cal B}_2$}
\put(0,159.5){${\cal B}_{r+3}$}
\put(38.5,144.75){${\cal B}_{L-r}$}
\multiput(88.5,172.75)(-.033557047,-.120805369){149}{\line(0,-1){.120805369}}
\multiput(84.75,155)(-.09375,-.03125){8}{\line(-1,0){.09375}}
\multiput(84.5,155)(.033557047,.117449664){149}{\line(0,1){.117449664}}
\multiput(86.25,152.25)(.06554878,-.033536585){164}{\line(1,0){.06554878}}
\multiput(85.75,151.5)(.065705128,-.033653846){156}{\line(1,0){.065705128}}
\multiput(83.75,155.25)(.1333333,-.0333333){15}{\line(1,0){.1333333}}
\multiput(85.75,154.75)(.0333333,-.0833333){15}{\line(0,-1){.0833333}}
\multiput(86.25,153.5)(.03125,-.15625){8}{\line(0,-1){.15625}}
\multiput(86.5,152.25)(-.03125,-.0625){8}{\line(0,-1){.0625}}
\put(88.75,172.5){\line(2,-1){5}}
\multiput(93.75,170)(.03333333,-.05){75}{\line(0,-1){.05}}
\multiput(96.25,166.25)(-.0543478,-.0326087){23}{\line(-1,0){.0543478}}
\multiput(95,165.5)(.0326087,-.0326087){23}{\line(0,-1){.0326087}}
\multiput(95.75,164.75)(.03370787,-.04775281){89}{\line(0,-1){.04775281}}
\multiput(98.25,160.75)(-.05,-.0333333){30}{\line(-1,0){.05}}
\put(97,147.25){\line(0,1){0}}
\put(96.5,146.25){\line(0,1){13.75}}
\put(89.5,158.5){$\Gamma^{'}$}
\put(10.25,154.25){$\Psi$}
\end{picture}

\end{figure}

To construct the path $z$, we 
go along the side of the band ${\cal B}_2$ which is closer to ${\cal B}_1$, then go along the part
of the boundary of the hub which is not part of $\partial\Psi,$ and finally go to $\partial\Delta$
along the side of ${\cal B}_{L-r}$ which is closer to ${\cal B}_{L-r+1}$. Thus we have $z=z_1z_2z_3$ 
according to this definition, and respectively, $\Lab(z) \equiv Z\equiv Z_1Z_2Z_3.$  Let $p$ be the subpath
of $\partial\Delta$ and $\partial\Psi$ such that $p_-=z_-, p_+=z_+$. Then the boundary path of 
$\Delta$ is of the form $pp',$ for an appropriate $p'$.

We denote by $\Gamma$ the part of the diagram $\Phi$ bounded by ${\cal B}_{r+3}$
and ${\cal B}_2$.  Let $y$ be the maximal common subpath of $\partial\Gamma$ and $p$  with $y_-=p_-=z_-$ that does not cross the band ${\cal B}_{r+3}$. 
We may apply Lemma \ref{sravn} to the pair $(\Gamma, \Phi)$ and  obtain a new diagram $\Gamma'$ over $G_1$
given by that lemma. One of the four boundary sections of $\Gamma'$ is a subword
$k_{i-1}...k_{i+r+1}$ of the word $\Sigma_0.$

Note that ${\cal B}_{L-r}$ is a $k_{i+r+1}$-band since we take the indices of the $k$-letters
 modulo $L.$ By Lemma \ref{sravn}, $\Gamma'$ has a loop with label of the form
$Z_1Z_2Z'Y'$, where $Z'$ is the copy of the word written along the band ${\cal B}_{r+3},$ and
$|Y'|\le |y|.$  Hence $Z= (Y')^{-1}(Z')^{-1}Z_3$ in $G_1.$ Therefore to prove that $Z$ is equal
in $G_1$ to a word of length $\le|p|,$ it suffices to prove that the word  $(Z')^{-1}Z_3$ is
equal in $G_1$ to a word of length $<|x|,$ where $p=yx.$

We observe that the words $Z'$ and $Z_3$ have equal prefixes of length $d$  since ${\cal B}_{L-r}[d]$
is a copy of ${\cal B}_{r+3}[d].$ Therefore the word $(Z')^{-1}Z_3$ is equal to $(\bar Z')^{-1}\bar Z_3,$ where
$(\bar Z')^{-1}$ copies the label of the part of the side of ${\cal B}_{r+3}$ 
crossed by the $\theta$-bands from the set $T$ only, and $\bar Z_3$ is the label of the part of $z_3$ crossed
by the $\theta$-bands of $T$ only. We denote the union of these two subsets of $T$ by $\bar T$. The length of  $(\bar Z')^{-1}\bar Z_3$ does not exceed the number $|\bar T|$ of
bands in $\bar T$ since neither of the bands from $T$ crosses both ${\cal B}_{r+3}$ and ${\cal B}_{L-r}$.
By Lemma \ref{NoAnnul}, every band of $\bar T$ must end on the subpath $x.$ Hence 
$|(\bar Z')^{-1}\bar Z_3|\le |x|.$ In fact this inequality is strict since $L-r > r+3$ (as $r\le 9$) and
so $x$ must include some $q$-edges as well.  The lemma is proved. \endproof

 \section{Spaces of words} \label{sw} 
 
 Let $\cal W$ be a set of vanishing in the group $G$ words. An easy observation (see Lemma \ref{kuski} below) shows that under a natural
 condition, spaces of  words from  $\cal W$
 can be bounded from above by the same (up to equivalence) function as the spaces 
 of  words from a smaller set $\cal W'.$ In this section, we define several sets
 of words (as boundary labels of several sets of diagrams) and estimate from
 above the spaces of words from larger sets using the upper bounds for the spaces
 of words from smaller sets. To apply Lemma \ref{kuski}, we will
 cut up a diagram $\Delta$ into two pieces such that one of the pieces belongs to
 a smaller set of diagrams while the perimeter of the second piece is bounded
 by $|\partial\Delta|.$ The most difficult statement is Lemma \ref{final} (separating
 of a one-hub subdiagram) whose proof is based on Lemma \ref{cut}.
 
 
 \subsection{Spaces of boundary labels of some diagrams.}

We call a disc  {\it simple} if it has no $H$-cells and
either all its $\theta$-edges have labels from $\Theta$ or all of them have labels
from $\hat\Theta.$ 

\begin{lemma} \label{simple} Let $\Delta$ be a diagram having exactly one hub $\Pi$. Then there is
a diagram $\bar\Delta$ with the same boundary label as $\Delta$ such that $\bar\Delta$ has a 
simple disc subdiagram $D$, and the annular diagram $\Gamma = \Delta\backslash D$ is a minimal
annular diagram without $\theta$-annuli. 

Moreover, one may assume that the boundary label
of $D$ is of the form  $(k_1W_1k_2W_2\dots k_L W_L)^{\pm 1},$ where $k_1W_1k_2W_2\dots k_L W_L$
is accepted by either the machine ${\cal S}(L)$ or by $\hat{\cal S}(L),$ and the lengths  of
$\theta$-annuli in $D$ do not exceed $N+L (space_{{\cal S}\cup\hat{\cal S}}(W_2))$.

\end{lemma}

\proof We may assume that $\Delta$ is a minimal diagram. Let $D_1$ be a maximal
disc in $\Delta,$ and denote by ${\cal K}_1,\dots,{\cal K}_L$ the 
maximal $k_1,\dots, k_L$-bands of $D_1$ starting on the hub $\Pi.$ We denote by
$\Gamma_i$ the maximal accepted $i$-sector of $D_1$ bounded by $k_i$ and $k_{i+1}$
($i=1,\dots,L$). By Lemmas \ref{calAi} and \ref{simul}(1), these sectors have
no $H$-cells for $i\ne 1$ and each of them is a copy or a mirror copy of the $2$-sector
$\Gamma_2.$ 

Now we replace the $1$-sector $\Gamma_1$ by the replica $\Gamma'_2$ of $\Gamma_2$ in $D_1.$
(To achieve this, one can make a cut along ${\cal K}_1$, ${\cal K}_2$ and the part of the
boundary of $\Pi$ between these $k$-bands, and insert two mirror copies of $\Gamma'_2$ 
along this cut.) By the definition of  replica, we obtain a modification $D_2$ of the disc 
diagram $D_1$ with boundary label of the form $k_1W_1k_2W_2\dots k_L W_L,$ where $W_1$ is 
a mirror copy of $W_2$ or $W_1$ has no $a$-letters. Since $k_2W_2k_3$ is an accepted $2$-sector word, 
the word  $k_1W_1k_2W_2\dots k_L W_L$ is accepted by either the machine ${\cal S}(L)$ or the machine $\hat{\cal S}(L)$
by Lemma \ref{inH}. Moreover the length  of this computation does not exceed the length of
the computation of ${\cal S}\cup\hat{\cal S}$ by the same lemma. Therefore by Lemma \ref{simul}(2),
 the disc $D_2$ can be replaced by a disc $D_3$ which has no $H$-cells, and whose labels of $\theta$-edges 
 either all belong to $\Theta$ or all belong to $\hat\Theta$, and  
 whose number of $\theta$-annuli does not exceed that number for $D_2.$  
 
 Now, if necessary, the annular diagram $\Delta\backslash D_3$ can be replaced by
 a minimal diagram $\Gamma$ over the group $G_1.$ Assume that $\Gamma$ has a $\theta$-annulus ${\cal T.}$
 By Lemma  \ref{NoAnnul}, ${\cal T}$ surrounds the disc diagram $D_3,$ and so $D_3$ can be included in
 a larger disc subdiagram $D_4$ which contains more $(k_i,\theta)$-cells for $i\ne 1,2$ since the extensions 
 of the ${\cal K}_i$-s have to cross $\cal T.$ Then one can make the surgery as above and replace
 $D_3$ by a larger simple disc $D_4.$ This procedure terminates, because we do not change
 the number of $(k_i,\theta)$-cells ($i\ne 1,2$) in the compliment of the discs when passing from $D_1$ to $D_2$ and
 from $D_2$ to $D_3$ and we reduce this number when passing from $D_3$ to $D_4.$  (Recall
 that the rank of such cells is higher than the ranks of other cells in diagrams over $G_1$.) 
 Thus the procedure terminates with a desired diagram $\Gamma.$ 
 
 Finally, one can replace the disc by a disc corresponding to a computation of minimal  space and use
 Lemma \ref{inH} to make 
 the second claim of the lemma true.  
 
 \endproof
 
 \begin{lemma}\label{podisku} There are positive constants $c_1$ and $c_2$ with the
 following property. For the boundary label $w\equiv k_1W_1\dots k_2W_l$ of the simple disc $D$ from Lemma
 \ref{simple}, there is a derivation
 $w\equiv w_0\to w_1\to\dots\to w_t=1$ using the relations of $G({\cal S}\cup\hat{\cal S},L)$ 
 and the hub relation (\ref{rel3}), such that $|w_i|\le c_1 S'_{\cal S}(|W_2|)+c_2$ ($i=0,1,\dots,t$).
 \end{lemma}
 
 \proof Let us start with the words $w\equiv w(0), w(1),\dots, w(m)\equiv\Sigma_0,$ written on the boundaries
 of the $\theta$-annuli of $D.$ By Lemma \ref{simple} and \ref{inH}, there is $c_1>0$ such that $|w(i)| \le c_1   
 S'_{\cal S}(|W_2|)+N.$ Notice that $w(i)$ is written on the top of a $\theta$ annulus $\cal T$ of $D$
 and $w(i+1)$ is written on its bottom. The band $\cal T$ has $N$ $(\theta,q)$-cells. For these cells
 the combinatorial length of their tops and bottoms differ by at most $\pm 1$. The remaining cells are
 $(\theta,a)$-cells, and their tops and bottoms have one $a$-edge. Therefore one can insert several
 elementary transformations between $w(i)$ and $w(i+1)$ corresponding to a sequential removal of the  cells of $\cal T$
 so that the combinatorial length of the words obtained after the refinement of our sequence does 
 not exceed  $|w(i)|+N +2.$ To complete the proof, it suffices to set $c_2=2N+2.$\endproof
 
 \begin{lemma}\label{klyaksy} 
 Assume that (1) a minimal diagram $\Delta$ has no hubs and no $q$-bands or (2) $\Delta$ 
 is a union of a simple disc  $D$ and a minimal annular diagram $\Gamma$ over $G_1$ surrounding
 the disc subdiagram $D$ and having no $\theta$-annuli, and every maximal $q$-band of 
 $\Delta$ starts on the hub of $D.$ Then the sum of perimeters $\sigma_H$ of all $H$-cells in $\Delta$ does
 not exceed $c_3|\partial\Delta|$ for some constant $c_3$ independent of $\Delta.$
 \end{lemma}  
 
 \proof We consider the condition (2) of the lemma only since a simplified argument works
 if $\Delta$ satisfies condition (1).
 
 Two different $H$-cells cannot be connected by an $a$-band in a minimal diagram since
 otherwise this subdiagram could be replaced by a diagram with one $H$-cell and several 
 $(\theta,a)$-cells contrary to the minimality of the diagram. (See Lemma 3.12 (2) in \cite{OS2}.)
 A maximal $a$-band $\cal A $ cannot connect $a$-edges of the same $H$-cell $\pi$ in a disc subdiagram by
 Lemma 3.12 (3) in \cite{OS2}, and also $\cal A$ and $\pi$ cannot surround the disc $D$ since
their boundaries have no $q$-edges. 
 Therefore every maximal $a$-band starting on an $H$-cell ends either on $\partial\Delta$
 or on $\partial D$, or on a $(\theta,q)$-cell of the annular diagram $\Gamma.$ Therefore
 to estimate  $\sigma_H,$ we should give an estimate for the number $n_{\theta,q}$
 of $(\theta,q)$-cells in $\Gamma$ and for the number of $a$-letters in the word $W_1,$
 where $k_1W_1\dots k_LW_L$ is the boundary label of $D.$ (Recall that the edges of $H$-cells
 are labeled by $a$-letters from the alphabet ${\cal A}_1$ and so cannot be connected by
 $a$-bands with subpaths of $\partial D$ labeled by $W_2,\dots, W_L.$)
 
 Let $\cal T$ be a maximal $\theta$-band of $\Gamma.$ Since both the start and the end $\theta$-edges of $\cal T$ must belong to $\partial\Delta$, it follows from Lemma \ref{NoAnnul} that $\cal T$ 
 crosses every
 $q$-band of $\Gamma$ at most once, and so has at most $N$ $(\theta,q)$-cells, because every
 maximal $q$-band of $\Gamma$ starts on the disc $D$. The number of maximal $\theta$-bands of
 $\Gamma$ does not exceed $\partial\Delta$ since $\Gamma$ has no $\theta$-annuli. It follows
 that $n_{\theta,q}\le N |\partial\Delta|.$ 
 
 Since we have $|W_1|_a\le |W_2|_a=|W_3|_a=\dots $ for the simple disk $D,$ we will look for an upper estimate
 for $|W_2|_a.$   Every maximal $a$-band $\cal A$ starting on the subpath $p$ of $\partial D$ labeled by $W_2$
 cannot end on $p$ by Lemma \ref{NoAnnul} since the word $W_2$ is reduced and has no $\theta$-letters. Therefore $\cal A$
 terminates on one of the two closest maximal $q$-bands $Q$ and $Q'$ starting on $D$ or on 
 $\partial\Delta.$ The lengths of $Q$ and $Q'$ are at most $|\partial\Delta|$ as was explained
 in the previous paragraph, and so each of the sides of these $q$-bands has at most $|\partial\Delta|$ 
 $a$-edges. Therefore $|W_2|_a\le 3|\partial\Delta|$.
 
 Thus, $\sigma_H\le c_3|\partial\Delta|$ for the constant $c_3=1+3+2N.$ \endproof
 
Assume that a word $w$ vanishes in the group $G$ given
 by Relations (\ref{rel1}), (\ref{rel2}), and (\ref{rel3}).
 Then we denote by $Space_G(w)$ the minimal number $m$ such that there is an elementary reduction
 of $w$ to the empty word such that at every step $i$, we have a tuple of words $(w_{i,1},\dots, w_{i, s(i)})$
with $|w_{i,1}|+\dots + |w_{i,s(i)}|\le m.$ 
For a set $\cal W$ of words vanishing in $G$ (call such
a set {\it vanishing}), we
define the space function $f_{G, \cal W}(x)$ as the maximum of $Space_G(w)$ over the words $w\in \cal W$
with $|w|\le x.$ (Thus we use the length $|\;\;|$ here unlike the length
$||\;\;||$ used in the Introduction.)

Now we denote by ${\cal W}_1$ the set of words read on the boundaries of simple discs,
 and we say that $w\in {\cal W}_2$ if $w$ can be read on the boundary of an $H$-cell. 
 
 \begin{lemma}\label{Hcell}  The function $f_{G,{\cal W}_2}$ is bounded from above by a function
 equivalent to  $f_{G,{\cal W}_1}.$
 \end{lemma}
 
 \proof Assume that $|w|=n>0$ and $w\in {\cal W}_2.$  Since the set of generators
 ${a_1,...,a_m}$ of $H$ is symmetric, for every $i\le m$, we have a positive relation of the form
 $a_ia_{i'}$ for some $i'\le m.$ By Lemma \ref{homo} these 2-letter relations are consequences 
 of the relations of $G$. Since the set of  $2$-letter relations is finite, there is a constant $c_4$
 such that one can convert $w$ letter-for-letter to
 a positive word $u$ of the same length, and the space of the corresponding derivation is $\le n+c_4.$
 The word $u$ is a product of cyclic shifts of the words $\Sigma(u)$ and $\hat\Sigma(u)$
 as  was explained at the end of Subsection \ref{conemb}. Here both words $\Sigma(u)$ and $\hat\Sigma(u)$
 belong to ${\cal W}_1$ and their lengths are at most $Ln+N\le(L+N)n.$ Therefore for any 
 $w\in {\cal W}_2,$ we have
 $space_G(w)\le 2f_{G,{\cal W}_1}((L+N)n)+c_4$ which implies the statement of the lemma. \endproof

 \begin{lemma} \label{kuski} Let $C>0$. Assume that for every $w$ from a vanishing set $ \cal W$, 
 there is a derivation $w\equiv w_0\to w_1\to\dots\to w_t\equiv 1$
 such that for every $i=0,\dots,t-1,$ a cyclic shift of the word $w_i$ is freely equal 
 to a product of a cyclic shift of $w_{i+1}$ and a word $v_i$ from a vanishing set $\cal W'$, where
 $\max(|v_i|,|w_i|)\le C|w|$ ($i=0,\dots,t-1$). Then the function $f_{G,\cal W}$ is bounded from
 above by a function equivalent to $f_{G,\cal W'}.$   
 \end{lemma}
 \proof The hypothesis of the lemma implies that we can apply the following series of elementary transformations
 to $w_i$. A first series of transformations replaces the word by its cyclic shift, a second series
 deletes/inserts mutual inverse letters operating with the words of length $\le 2C|w|$; then we split the obtained word into a product of
 a cyclic shift of $w_{i+1}$ and the word $v_i;$ then we keep $w_{i+1}$ unchanged and use the appropriate
 procedure reducing the word $v_i$ to the empty word, and finally obtain the word $w_{i+1}$
 using cyclic shifts. Clearly, we have that the space of this procedure is at most
 $2C|w|+f_{G,W'}(C|w|).$  Thus, by induction on $i$, we have $Space_G(w)\le 2Cn+ f_{G,W'}(Cn)$ for arbitrary
 word $w\in \cal W$ of length at most $n$. The lemma is proved. \endproof

 Now we introduce the set ${\cal W}_3$ of  boundary labels of diagrams $\Delta$ satisfying the condition 
 of Lemma \ref{klyaksy}, i.e., either (1) $\Delta$ is a minimal diagram $\Delta$ having no hubs and 
 no $q$-bands or (2) $\Delta$  is a union of a simple disc  $D$ and a minimal annular diagram $\Gamma$ over 
 $G_1$ surrounding the disc subdiagram $D$ and having no $\theta$-annuli, and every maximal $q$-band of 
 $\Delta$ starts on the hub of $D.$
 
 \begin{lemma} \label{bezq} The function $f_{G,{\cal W}_3}$ is bounded from above by a function
 equivalent to  $f_{G,{\cal W}_1\cup{\cal W}_2}.$ 
 \end{lemma}
 
 \proof We will assume that a word $w$ of length $n\ge \delta$ is the boundary label of a diagram $\Delta$ 
 satisfying condition (2) in the definition of the set ${\cal W}_3.$ Every cell of the annular subdiagram
 $\Gamma$ is either an $H$-cell or a $\theta$-cell. Therefore there is a sequence of diagrams 
 $\Delta=\Delta_0, \Delta_1,\dots,\Delta_t=D$ such that for $i=1,\dots, t$, the diagram $\Delta_i$
 results from $\Delta_{i-1}$ after one cuts off either (a) an $H$-cell or  (b) a rim $\theta$-band,
 or an edge $e$ such that $ee^{-1}$ belongs to the boundary path of $\Delta_{i-1}.$
 Surgery of types (b) and (c) decreases the perimeter by Lemma \ref{rim}. Although  surgery of type
 (a) can increase the perimeter, it follows from Lemma \ref{klyaksy} that the perimeter of
 every diagram $\Delta_i$ is at most $(1+c_4)n.$
 
 Now we have a sequence  $w\equiv w_0, w_1,\dots, w_t, 1$ where $w_0,\dots,w_t$ are the boundary 
 labels of $\Delta_0,\Delta_1,
 \dots,  \Delta_t=D,$ of lengths at most $(1+c_4)n,$ such that, for every $i=0,\dots,t,$ a cyclic shift
 of the word $w_i$ is a product of a cyclic shift of the word $w_{i+1}$ and a word $v_i,$ where
 $v_i$ is either a boundary label of an $H$-cell of $\Delta$ or a $2$-letter word $aa^{-1}$, or the boundary label of the simple disc $D,$
 or the boundary label of the the rim $\theta$-band of $\Delta_i.$ In all of these cases $|v_i|\le 2(1+c_4)n;$
 in the case of the rim band, we obviously have $Space_G |v_i|\le 2(1+c_4)n,$ and in the other cases 
 $Space_G |v_i|\le f_{G,{\cal W}_1\cup{\cal W}_2}((1+c_4)n).$ Therefore one can apply Lemma
 \ref{kuski} and complete the proof. \endproof
 
 By definition, the set of words ${\cal W}_4$ contains the set ${\cal W}_3$ and consists of 
 boundary labels of diagrams $\Delta,$ where either (1) $\Delta$ is a minimal diagram having no hubs or
 (2) $\Delta$  is a union of a simple disc  $D$ and a minimal annular diagram $\Gamma$ over 
 $G_1$ surrounding the disc subdiagram $D$ and having no $\theta$-annuli.
 
 \begin{lemma} \label{onehub} The function $f_{G,{\cal W}_4}$ is bounded from above by a function
 equivalent to  $f_{G,{\cal W}_3}.$ 
 \end{lemma}
 
 \proof Again we assume that a word $w$ of length $n\ge \delta$ is the boundary label of a diagram $\Delta$ 
 satisfying condition (2) in the definition of the set ${\cal W}_4.$ 
 Assume that $\Delta$ has a maximal $q$-band $\cal C$ which does not start or terminate on the simple disc $D.$
 Then $\Delta$ is separated into $3$ subdiagrams: $\Gamma_1$ contains the disc $D$, $\Gamma_2=\cal C$, and
 $\Gamma_3$ is the remaining part of $\Delta.$ On the one hand, the lengths of the top and the bottom of 
 $\cal C$ are equal to the number of cells $m$ in $\cal C$ since every cell of $\cal C$ has one $\theta$-edge 
 and at most one $a$-edge on each of the sides of $\cal C.$ On the other hand, each maximal $\theta$-band
 of $\Delta$ crossing $\cal C$ must start and terminate on $\partial\Delta.$ This implies that
 the perimeters of  each of the subdiagrams $\Gamma_1, \Gamma_2,$ and $\Gamma_3$ are at most $2n$.
 (Here we use the definition of length and take into account that the band $\cal C$ starts and ends on
 $q$-edges.) 
 
 Therefore there is a sequence of diagrams 
 $\Delta=\Delta_0, \Delta_1,\dots,\Delta_t$ of perimeters $\le 2n$ such that for $i=1,\dots, t-1$, the diagram $\Delta_i$
 results from $\Delta_{i-1}$ after one cut off either (a) a subdiagram without $q$-bands or (b) a $q$-band,
 and $\Delta_t$ has no maximal $q$-bands except for those starting/terminating on the hub of $D.$ 
 Let $w\equiv w_0, w_1,\dots, w_t$ be the boundary 
 labels of these diagrams. Every word $w_i$ ($i=0,\dots,t$) of the series $w\equiv w_0, w_1,\dots, w_t, w_{t+1}=1$ 
 (or its cyclic shift) is a product of a cyclic shift of the word $w_{i+1}$ and a word $v_i,$ where
 $v_i$ either belongs to ${\cal W}_3$ or it is
 the boundary label of a $q$-band. In all of these cases $|v_i|\le 2n,$ and
 in the later case we obviously have $Space_G |v_i|\le 2n.$ In the former case we have 
 $space_G |v_i|\le f_{G,{\cal W}_3}(2n).$ To complete the proof, we apply Lemma
 \ref{kuski}. \endproof

 \begin{lemma} \label{final} Let ${\cal W}_5$ be the set of all words vanishing in $G$. The  function $f_{G,{\cal W}_5}(n)$ is bounded 
 from above by a function
 equivalent to  $f_{G,{\cal W}_3\cup{\cal W}_4}(n).$ 
 \end{lemma}
 \proof Let $w=1$ in $G$ and $|w|=n>0.$ If $w=1$ in $G_1,$ then $Space_G(w)\le Space_{G,{\cal W}_3}(w)$ 
 by Lemma \ref{onehub}. Otherwise the minimal diagram $\Delta$ with boundary
 label $w$ has $t\ge 1$ hubs, and by Lemma \ref{cut}, a cyclic shift of the word $w\equiv w_0$ is a product of
 a word $v_1$ written on the boundary of a diagram $\Gamma$ having one hub, and a word $w_1$ which
 is a boundary label of a diagram with $t-1$ hubs, and $|w_1|\le n, |v_1|\le 2n.$  By Lemma \ref{simple},
 $v_1\in {\cal W}_4.$  Now induction on $t$ gives a series of words $w_0,w_1,\dots, w_t=1$ and words
 $v_1,\dots,v_t$ satisfying the conditions of Lemma \ref{kuski} with $C=2,$ and our statement follows
 from that lemma. \endproof 
 
 \subsection{Proofs of main statements.}\label{vse}
 
 {\bf Proof of Theorem \ref{main1}.} Suppose  $H$ is a group with a finite generating set
 $\{a_1,\dots,a_m\}$ and the word problem for $H$ is solvable by a DTM with space complexity $f(n).$ As we mentioned in the beginning of Subsection \ref{conemb}, one can
double the set of generators so that every $a_i^{-1}$ becomes equal in $H$ to some $a_j.$ 
Hence the set of positive words vanishing in $H$ is a set of defining relations for $H$,
and it is recognizable by a DTM $M$ with space complexity $f(n).$
 By Lemma 2.10 (a), there is
 an $S$-machine $\cal S$ recognizing the same language, and the generalized space complexity
 $S'_{\cal S}(n)$ of $\cal S$
 is equivalent to $f(n).$    The group $G$ constructed on the basis of $\cal S$ 
 in Subsection \ref{conemb} is finitely presented and contains $H$ as a subgroup by 
 Corollary \ref{inject}. The consecutive application of Lemmas \ref{final}, \ref{onehub},
 \ref{bezq}, \ref{Hcell}, and \ref{podisku}
 results in the inequality $f_{G,{\cal W}_5}(n)\preceq  S'_{\cal S}(n).$ 
 Note that the space functions $f_{G,{\cal W}_5}(n)$ and $s_G(n)$ of $G$ are equivalent since
 the lengths functions $||*||$ and $|*|$ satisfy inequalities $\delta||w||\le|w|\le ||w||$
 for every word $w$.  Hence  $s_G(n)\preceq  S'_{\cal S}(n).$ Since $S'_{\cal S}(n)\sim f(n)$ by Lemma
 \ref{MtocalS} (2), we have $s_G(n)\preceq f(n).$ 
 
 To invert this inequality and to obtain the second claim of Theorem \ref{main1}, we first  note that $s_G(n)\succeq \log d_G(n),$ where $d_G(n)$ is the Dehn function for $G.$ (Indeed, up to equivalence, the length $t$ of a rewriting $W_0=(w_0)\to\dots
 \to W_t=(1)$ without repetitions 
 does not exceed $\exp(\max_{i=0}^t ||W_i||);$ see also Theorem C in \cite{BR}).  Therefore
 it suffices to show that $f(n)\preceq \log d_G(n).$
 By Lemma \ref{generspace}(3), $S'_{M'}(n) \sim f(n)$, 
 and by Lemma \ref{spacecalS} (4),  $T'_{\cal S}(n) \succeq \exp (S'_{M'}(n)).$ Thus it remains to
 explain that $d_G(n)\succeq T'_{\cal S}(n).$
 
 Let $W$ be a word accepted by $\cal S,$ such that $1\le |W|_a\le n$ and
 $time_{\cal S}(W)=T'_{\cal S}(n).$  Then the word $V\equiv k_1W_1\dots k_LW_L$ 
 (where $W_i$-s are copies or mirror copies of $W$) is accepted by
 ${\cal S}(L)$ and therefore it is conjugate to the
 word $\Sigma_0$ (see subsection \ref{conemb}). Hence $V=1$ in $G$.
 Let $\Delta$ be a minimal diagram over $G$ with boundary label $V.$
 Since $\partial\Delta$ has only one $k_1$-edge, the maximal $k_1$-band starting
 on this $k_1$-edge must end on the boundary of a hub. Hence $\Delta$
 has a hub $\Pi$ satisfying the condition of Lemma \ref{extdisc}. If one removes
 $\Pi$ together with the bands ${\cal B}_1,\dots , {\cal B}_{L-3}$ and
 subdiagrams $\Gamma_i$ ($i=1,\dots,L-4$), then the remaining diagram
 $\Delta'$ has at most $L-(L-3)+3=6$$\;\;$ $k$-edges. It follows from Lemma
 \ref{extdisc} that $\Delta'$ has no hubs since $6<L-3.$ Thus $\Delta$ has
 exactly one hub.
 
 Obviously, every $k_i$-band starting on $\partial\Pi$ ends on $\partial\Delta$
 ($i=1,\dots,L$), and so we can consider the $2$-sector $\Gamma$ of $\Delta$ bounded
 by $k_2$- and $k_3$-bands. By Lemma \ref{simul}(1), the number of maximal $\theta$-bands
 of $\Gamma$ is equal to the length of a computation of ${\cal S}(L)\cup{\hat{\cal S}}(L)$
 accepting the word $W.$ By Lemma \ref{conemb}, we can replace ${\cal S}(L)\cup{\hat{\cal S}}(L)$ by $\cal S$ in the previous phrase, and so $area(\Delta)> T'_{\cal S}(n).$
Since $||V||\le C|W|_a$ for a constant $C,$ we have $d_G(Cn)>T'_{\cal S}(n),$
and the required lower bound is obtained. $\Box$

\medskip

{\bf Proof of Theorem \ref{realiz}.} Let $M$ be a $DTM$ with space complexity $f(n).$
Then, as in the proof of Theorem \ref{main1}, we can construct an $S$-machine $\cal S$
with $S'_{\cal S}(n)\sim f(n)$ and $T'_{\cal S}(n) \sim \exp f(n).$
 But now we simplify the construction of the group $G$:
It is defined by the relations associated with the machine ${\cal S}(L)$ only and the
hub relation (there is no $\hat{\cal S}(L)$ now, and so we have no $H$-relations at all). The
(simplified) proof of Theorem \ref{main1} works in this setting and
therefore $s_G(n)\sim f(n)$ and $\log d_G(n) \sim  f(n).$ 
Finally, by Lemma \ref{simul} (1,2) and by the argument exploited in the proof of Theorem
\ref{main1}, an input word $u$ of $M$ is accepted iff
the word $\Sigma(u)$ vanishes in $G;$ this completes the proof.  
$\Box$

\medskip

{\bf Proof of Corollary \ref{main}.} Let an $NTM$ have  $FSC$ space complexity $f(n)$
and solve the word problem in a finitely generated group $H$.
Then by Savitch's Theorem (see \cite{DK}, Theorem 1.30), there is a $DTM$ which solves
the same problem with space $\sim f(n)^2.$ It remains to refer to Theorem \ref{main1}. $\Box$

\begin{rk} \label{extend} Corollary \ref{main} has a problem-reduction version. If 
a $DTM$ $M$ has an $FSC$ space complexity function $f(n),$ then one can construct  a finitely presented 
group $G=G(M)$ according to Theorem \ref{realiz}, and so $G$ has a 
space function equivalent to $f(n)$. By Proposition \ref{propos}, 
the word problem of $G$ is solvable by an $NTM$ of space complexity $f(n)$,
and so it is solvable by a $DTM$ with space $f(n)^2.$ Therefore Corollary
\ref{complete} can be extended to say that many space complexity classes 
have space complete word problem for some finitely presented group 
(depending on the class) with 
space function equivalent to the given space complexity $f(n).$
\end{rk}

\medskip

We need one more lemma to prove Corollary \ref{alpha}. This is a version of
Savitch's theorem (see Theorem 1.30 in \cite{DK}), but instead of  simulation of the work of an
$NTM$ by a $DTM,$ now we need a $DTM$ computing the space complexity of a given $NTM.$

\begin{lemma}\label{savitch}
Let $M$ be an $NTM$ with space complexity $S(n)$ bounded from above by
an $FSC$ function $f(n).$ Then there is a $DTM$ $M_0$ such that (1) $M_0$ computes
$S(n)$, i.e., for any input $n\in \mathbb{N}$ given in binary, it computes $S(n)$;
(2) the space complexity $S_{M_0}(n)$ is $O(f(n)^2).$ 
\end{lemma}

\proof Without loss of generality, we may assume that $M$ satisfies the $\vec s_{10}$
condition. 

If $u$ is an accepted input word for $M$ and $||u||\le n,$ then the time of any
computation (without repetitions) of space $\le f(n)$ accepting $u$ is at most $\le 2^{cf(n)},$
for some integer $c>0$ since  all configurations in this computation are of length $\le  f(n)+c_0$ for a constant $c_0.$ So the goal
for the $DTM$ $M'$ we want to define,  is  to find a computation $C$ of $M$ of minimal space and of length at most $2^{cf(n)}$ which connects the input configuration  and the accept configuration of $M,$ and then to compute the space of this computation. (If such a
computation exists; otherwise $M'$ says that $u\notin {\cal L}_M.$)
Indeed, the required machine $M_0$ will examine all words $u$ with $||u||\le n$ in
lexicographical order, it will switch on $M'$ for every such input word $u,$
and it will compare the spaces $space_M(u)$ of the $u$-s on an additional tape keeping
only the maximal one after every return.

For arbitrary words $w$ and $w'$ and $k\in \mathbb{N}$, define 
predicates $reach_n(w,w',k)$ to mean that  $w$ and $w'$ are configurations of $M$
and there is a computation $w\to\dots \to w'$ with time $\le k$ and with
space $\le f(n)+c_0.$ By this definition, $M$ accepts an input word $u$  of
combinatorial length $\le n$ iff $reach_n(w_0, w_f, 2^{cf(n)}),$ where
$w_0\equiv w(u)$ is the unique input configuration on input $u$ and $w_f$ is
the unique accept configuration of $M.$ 

Note that $reach_n(w,w',k+j)$ iff $(\exists w'')(reach_n(w,w'',k)\;\; and \;\; reach_n(w'',w',j))$ and the minimal space of computations $w\to\dots\to w'$ is
the minimum over all $w''$ of the maximums of minimal spaces for $w\to\dots \to w''$
and for $w''\to\dots\to w'.$  This observation leads to the following (slight)
modification of Savitch's machine.

\medskip 

\begin{itemize}

\item Given $n\in \mathbb{N},$ compute $2^{cf(n)}$ (in binary) using that $f(n)$ is an $FSC$
function.

\item Then, $reach_n(w,w',1)$  is true if $w\equiv w'$ or $w$ transforms into $w'$ under
the application of a single command of $M.$ The  space of the computation $w\to w'$
is $\max(|w|_a, |w'|_a).$  

\item If $k\ge 2$, then for all possible configurations $w''$ of $M$ with length $\le f(n)+c_0,$
compute whether it is true that $reach_n(w,w'',[(k+1)/2])$ and $reach_n(w'',w',[(k+1)/2)].$
Set $reach_n(w,w',k)$ to be true iff such $w''$ exists. Find the minimal space of
computations $w\to\dots \to w'$ of length $\le k$ using the information on
the minimal space of computations $w\to\dots\to w''$ and $w''\to\dots\to w'$
of length $\le[(k+1)/2].$

\end{itemize}
\medskip

It is easy to see that the constructed machine $M'$ computes, in particular,
the space of any accepted input word $u$ of length $\le n.$ Passing from
$k$ to $[(k+1)/2],$ we need additional space to store the information on
$k,$ $w$, $w',$ on the current $w'',$ and afterwards,  on the minimal space of computations $w\to\dots\to w'$ of length $\le k.$ Clearly, this additional space is $O(f(n))$,
and since we start with $k=2^{cf(n)}$, we divide $k$ by $2\;\;$  $cf(n)$ times, so
the total space used by $M'$ and by $M_0$ is $O(f(n)^2)$.

\endproof 

\medskip

{\bf Proof of Corollary \ref{alpha}.}
Assume that $\alpha$ is computable with space $\le 2^{2^m}$.
It follows that for $m= [\log_2\log_2 n]$ we can recursively compute 
binary rationals $\alpha_m$ such that 
\begin{equation}\label{nalpha}
|\alpha-\alpha_m|=O(2^{-m})= O((\log_2 n)^{-1})
\end{equation}
and the space of the computation of $\alpha_m$ is at most $n.$
In addition, one may assume that the number of digits in the
binary expansion of $\alpha_m$ is $O(m).$ Therefore the computation
of $[\log_2n]$ (in binary) and of the product $\alpha_m[\log_2 n]$ needs  space at most $O((\log_2 n)^2).$  
Then we rewrite the binary presentation of $[\alpha_m[\log_2 n]]$
in unary (as a sequence of $1$-s). This well-known rewriting (e.g., see p.352 in \cite{SBR}) has 
space complexity of the form $[\alpha_m[\log_2 n]]+O(1).$ One more rewriting
of this type applied to the unary presentation of $[\alpha_m[\log_2 n]],$
will have space complexity of the form $2^{[\alpha_m[\log_2 n]]}+O(1).$
Using (\ref{nalpha}), we can present this function as $$2^{\alpha[\log_2 n]+O(1)}+O(1)
\sim 2^{\alpha[\log_2 n]} \sim [n^{\alpha}]$$

Thus the subsequent application of the above mentioned $DTM$-s has space complexity equivalent to $n^{\alpha}$,
and we can apply Corollary \ref{realiz} to obtain a finitely presented group with space function equivalent to $n^{\alpha}.$

   Now assume that a function $[n^{\alpha}]$ is equivalent to a space function of a
   finitely presented group $G.$ Then by Proposition \ref{propos}, there is an $NTM$
   $M$ whose space complexity $S_M(n)$ is equivalent to $[n^{\alpha}],$ that is
   
   \begin{equation}\label{equival}
   c_1n^{\alpha} < S_M(n)< c_2n^{\alpha}
   \end{equation} 
   for some positive $c_1,$ positive integer $c_2,$ and
   every sufficiently large $n$. In particular, we have $S_M(n) < c_2 n^d$ for some integer $d$ 
   and every $n.$  Since $c_2n^d$ is an $FSC$ function, we may apply Lemma \ref{savitch},
   and obtain a $DTM$ $M_0$ computing the function $S_M(n)$ with space $O(n^{2d})$.
   Hence $M_0$ computes the function  $S_M(2^{2^m})$ of $m$ with space  $O((2^{2^m})^{2d}).$
   This space is less than $2^{2^{m+c_3}}$ for some $c_3.$ Hence for some  $c_4\in \mathbb{N}$,
   $M_0$ computes the function $S_M(2^{2^{m-c_4}})$ with space at most $ 2^{2^{m-1}}.$
   
   Let us plug $n=2^{2^{m-c_4}}$ into Inequalities (\ref{equival}) and then take 
   $\log_2$ of the terms. We obtain
   
   $$ \lambda_1+\alpha 2^{m-c_4}\le \log_2 S_M(2^{2^{m-c_4}})\le \lambda_2+\alpha 2^{m-c_4}$$
   where $\lambda_i=\log_2c_i$ ($i=1,2$). It follows that  
   \begin{equation}\label{near}
   |\alpha- 2^{-m+c_4}\log_2S_M( 2^{2^{m-c_4}})|< c2^{-m+c_4}= O(2^{-m})
   \end{equation}
    where $c=\max(|\lambda_1|,|\lambda_2|).$  Recall that $S_M(2^{2^{m-c_4}})\le c_2(2^{2^{m-c_4}})^d\le 2^{2^{m-1}}$ and so this number has at most $2^{m-1}+1$ binary digits. Therefore the real numbers
   $\log_2 (S_M(2^{2^{m-c_4}})2^{-m+c_4}) $ are computable with error $O(2^{-m})$ and space  $O(2^{2^{m-1}}).$
   Now it follows from (\ref{near}) that the real number $\alpha$ is  computable with space
   $2^{2^{m}}.$
$\Box$
\medskip

{\bf Acknowledgment.} 

The author is grateful to J.-C. Birget, M. Bridson, T.Davis, M. Elder, S.V.Ivanov,
T.Riley, M. Sapir, and the anonymous referee for useful discussions, comments, and criticism.

\end{document}